\documentclass[fleqn]{elsarticle}
\usepackage{amsmath,amsfonts,graphicx,subfigure,epstopdf,float,amssymb,bm,amsthm,multirow,appendix,bbm,color,enumitem}

\usepackage{geometry}
\numberwithin{equation}{section}
\newtheorem{thm}{Theorem}[section]

\newtheorem{rmk}{Remark}[section]
\newtheorem{lem}{Lemma}[section]
\newtheorem{prop}{Proposition}[section]
\newtheorem*{prf}{Proof}

\allowdisplaybreaks[4] 

\begin{document}

\begin{frontmatter}

	\title{An explicit and practically invariants-preserving method for conservative systems}
	
	\author{Wenjun Cai$^{a}$, Yuezheng Gong$^{b}$, Yushun Wang$^{a,*}$}
	\address[1]{Jiangsu Key Laboratory for NSLSCS, \\
		School of Mathematical Sciences, Nanjing Normal University, Nanjing 210023, China}
	\address[2]{College of Science, Nanjing University of Aeronautics and Astronautics, Nanjing 210016, China}
	
	\begin{abstract}
		An explicit numerical strategy that practically preserves invariants is derived for conservative systems by combining an explicit high-order Runge-Kutta (RK) scheme with a simple modification of the standard projection approach, which is named the explicit invariants-preserving (EIP) method. The proposed approach is shown to have the same order as the underlying RK method, while the error of invariants is analyzed in the order of $\mathcal{O}\left(h^{2(p+1)}\right),$ where $h$ is the time step and $p$ represents the order of the method. When $p$ is appropriately large, the EIP method is practically invariants-conserving because the error of invariants can reach the machine accuracy. The method is illustrated for the cases of single and multiple invariants, with regard to both ODEs and high-dimensional PDEs.	Extensive numerical experiments are presented to verify our theoretical results and demonstrate the superior behaviors of the proposed method in a long time numerical simulation. Numerical results suggest that the fourth-order EIP method preserves much better the qualitative properties of the flow than the standard fourth-order RK method and it is more efficient in practice than the fully implicit integrators.
	\end{abstract}
	
	\begin{keyword}
		conservative systems;
		explicit Runge-Kutta method;
		explicit invariants-preserving method;
		practically invariants-conserving.
	\end{keyword}
	
\end{frontmatter}

\begin{figure}[b]
	\small \baselineskip=10pt
	\rule[2mm]{1.8cm}{0.2mm} \par
	$^{*}$Corresponding author.\\
	E-mail address: wangyushun@njnu.edu.cn (Y. Wang).
\end{figure}

\section{Introduction}

For an autonomous ordinary differential equation system with the initial condition
\begin{equation}\label{eq-1-1}
	\dot{y}=f(y), \quad y(0)=y_0, \quad y(t)\in \mathbb{R}^d,
\end{equation}
where $f:\mathcal{D}\subset\mathbb{R}^d\rightarrow \mathbb{R}^d$ is sufficiently smooth, a scalar function $I(y)$ is called a first integral or strong invariant if $$\nabla I(y)^T f(y) = 0, \quad y\in \mathcal{D}.$$ A function $I(y)$ is called a weak invariant of \eqref{eq-1-1} if $y_0\in \mathcal{M}$ implies the exact solution $y(t)\in \mathcal{M}$ and
$$\frac{\textrm{d}}{\textrm{d} t}I(y) = \nabla I(y)^T f(y) = 0.$$
Note that strong invariants are independent of initial values, while weak invariants depend on appropriate initial conditions. These invariants play an important role both in qualitative and quantitative studies of the flow of \eqref{eq-1-1}. Nowadays the ability to preserve some invariant properties is a criterion to judge the success of a numerical integrator \cite{li-95-FD-KG}.

In the past decades, there has been an increasing interest in invariant-preserving integrators that preserve as much as possible the invariants of underlying differential systems. It is well known that all Runge-Kutta (RK) methods preserve linear invariants, while only those that satisfy the symplectic condition conserve all quadratic invariants \cite{cooper-87-RK}. However, for $k \geq 3$ no RK method can conserve all polynomial invariants of degree $k$ \cite{hlw06}. In order to develop numerical methods that conserve general invariants, McLachlan et al. rewrote the initial problem with a first integral into a gradient system and then proposed discrete gradient (DG) methods, which are based upon the definition of a discrete counterpart of the gradient operator \cite{mcLachlan-99-DGM}. There are several choices of discrete gradients, including the coordinate increment DG \cite{itoh-88-DGM}, the midpoint DG \cite{gonzalez-96-DGM} and the mean value DG \cite{harten-83-DGM} (also called the averaged vector field (AVF) method in \cite{quispel-08-EP-AVF}) as well as their extensions \cite{wu-13-AAVF,cai-18-PAVF}. Later, these ideas are generalized for conservative PDEs by Furihata, Matsuo and collaborators using the concept of discrete variational derivatives (DVD) \cite{furihata-99-EP-DVD,matsuo-01-EP-linearly-implicit,furihata-11-DVDM}. In recent years, high-order integrators have also been widely developed, including the high-order AVF methods \cite{H10,CH11bit,LWQ16}, Hamiltonian boundary value methods (HBVMs) \cite{BIT10,BCMR12,brugnano-16-LIM}, time finite element methods \cite{BS00,TS12} and so on. Unfortunately, the structure-preserving algorithms mentioned above are often fully implicit for general conservative systems, which require a nonlinear iteration to solve them.

In contrast to fully implicit schemes, linearly implicit methods only involve to solve a linear system at each time step, and therefore have attracted continuous attention in recent years. To design linear-implicit conservative schemes, Matsuo and Furihata proposed the multiple points DVD method, which was a generalization of the DVD approach \cite{matsuo-01-EP-linearly-implicit}. More generally, Dahlby and Owren utilized the polarization technique to develop a general framework for deriving linearly implicit conservative algorithms for PDEs with polynomial invariants \cite{dahlby-11-general-SP}. More recently, the invariant energy quadratization (IEQ) \cite{yang-16-IEQ,yang-17-IEQ,gong-18-IEQ-binary-fluid} and the scalar auxiliary variable (SAV) approaches \cite{shen-18-SAV,shen-19-review,qiao-19-SAV}, originally proposed for dissipative gradient flow models, have been successfully applied to various conservative systems \cite{cai-19-BD-SG,cai-20-linearly-implicit-EP-MS,jiang-20-high-order-IEQ-CH,jiang-20-multiple-SAV-CH,chen-20-IEQ-wave,zhang-20-IEQ-explicit-RK}. However, these linearly implicit methods either conserve the polarized invariant or the modified quadratic energy.

Although linearly implicit methods can already reduce the computational cost, it is more preferable to design explicit invariant-preserving algorithms for conservative problems. So far, there are only a few related works existed in the literature.  When considering special second-order systems $\ddot{y}=f(y)$, the symplectic Runge-Kutta-Nystr\"om (RKN) methods can naturally preserve quadratic invariants of the form $y^\top D\dot{y}$ where $D$ is a matrix of appropriate dimensions, and importantly various explicit RKN methods with high order were already constructed (see \cite{qin-91-RKN,calvo-93-RKN} and references therein). For more general invariants of second-order problems, two kinds of explicit schemes based on the DVD method were proposed in \cite{furihata-01-explicit-EP-DVD,matsuo-07-explicit-EP-DVD} to preserve the energy of nonlinear wave equations. For quadratic invariants of general systems, a rational four-stage explicit RK method was derived by del Buono and Mastroserio \cite{buono-02-explicit-RK-quadratic-invariant}, which can be viewed as an incremental direction projection method and has been further generalized by Calvo et al. \cite{calvo-06-explicit-RK-quadratic-invariant} to any high-order explicit RK method. Combining the IEQ approach, general invariants can be reformulated into quadratic forms. Based on this fact, Zhang et al. \cite{zhang-20-IEQ-explicit-RK} extended the incremental direction projection method to construct explicit conservative schemes for general Hamiltonian ODEs and PDEs. A similar but more efficient strategy was recently proposed by Jiang et al. \cite{jiang-20-IEQ-explicit-RK} that incorporates high-order explicit RK methods with standard orthogonal projection techniques. The key of the above projection methods that make the resulting schemes explicit is attributed to the quadratization of the invariant, and as a consequence, the related unknown parameter (also called the Lagrangian multiplier) can be explicitly solved.  Otherwise, for non-quadratic invariant, nonlinear iterations are inevitable to solve this parameter in the standard projection methods and the resulting schemes are no long strictly explicit, although usually only one nonlinear equation is involved when considering to preserve single invariant. However, it should be noticed that the modified quadratic invariant by the IEQ approach is equivalent to the original invariant in the continuous case but will be different after discretizations. Therefore, either the increment direction projection method \cite{zhang-20-IEQ-explicit-RK} or the orthogonal projection method \cite{jiang-20-IEQ-explicit-RK} actually preserve a modified invariant. Moreover, these two kinds of methods are only suitable for single invariant because there are no explicit solutions for unknown parameters in multiple invariants cases.

The above mentioned structure-preserving algorithms can be rigorously proved to conserve some invariant properties exactly in theory, but they tend to achieve a practical conservation in numerical computing due to the following reasons. Firstly, nonlinear iterations are required for fully implicit methods whereas linear solvers are usually employed for linearly implicit methods, where an iteration tolerance must be given that may induce accumulated errors in the invariant. Secondly, as all the computations are done on a finite precision computer, even the explicit methods will inevitably suffer round-off errors. Based on this fact, in this paper we propose a novel class of explicit and practically invariants-preserving methods to capture the preservation of original invariants with general forms, which can also be recast into the framework of projection methods. Consider the nonlinear iterations for the unknown parameters in the standard orthogonal projection approach, the basic idea of the proposed method is to stop the Newton iteration after only one step, and thus obtain an explicit form of the unknown parameters. Similar idea can be found in \cite{brugnano-12-EP}, without theoretical analysis, to correct the drift of energy for Hamiltonian ODEs where the underlying methods are chosen as conservative HBVMs. In the present work, we show that any explicit high-order RK methods, even non-conservative, can be taken into our framework to achieve a practical invariant preservation for both conservative ODEs and PDEs. More specifically, we give a rigorous proof that such simplification does not affect the accuracy of the projection method, which retains the same order as the underlying RK method. Meanwhile, the invariants are also preserved to round-off errors as long as high-order RK methods are employed. The proposed schemes are illustrated by many practical applications, including the perturbed Kepler problem, the solar system, the charged particle dynamics and the rotating Gross-Pitaevskii equation. Numerical results suggest that the proposed methods perform the superior behaviors in a long time simulation and essentially improve the numerical performance of the standard fourth-order RK method.

The rest of this paper is organized as follows. In section 2, we present the explicit and practically invariants-preserving method and analyze its order of accuracy as well as the corresponding error in the invariant. In section 3, extensive numerical tests are reported to confirm the theoretical results and demonstrate the efficiency of the proposed method. Finally, concluding remarks are drawn in section 4.

\section{Explicit and practically invariants-preserving method}

In this section, we first present the algorithm of the explicit invariants-preserving (EIP) method. Afterward, the analysis on its accuracy and invariants-preserving property are carried out. Some implementation issues are provided in the last of this section.

Suppose we have an $(d-l)$-dimensional submanifold of $\mathbb{R}^d$,
\begin{equation}
	\mathcal{M} = \{y; g(y) = 0\}
\end{equation}
$(g: \mathbb{R}^d\longrightarrow \mathbb{R}^l),$ and the system \eqref{eq-1-1} with the property that
\begin{equation}
	y_0\in \mathcal{M} \qquad \textrm{implies} \qquad y(t)\in\mathcal{M} ~ \textrm{for~all~} t.
\end{equation}
Note that all components $g_i(y)$ of $g(y)$ are regarded as weak invariants of \eqref{eq-1-1}. Denote $y_n$ be the approximation of $y(t)$ at $t_n=t_0+nh$ with a step size $h$, $n=0, 1,2,\cdots$. Let $\Phi_h$ define an explicit RK method of order $p$ with time step $h$. For given $y_n\in \mathcal{M}$, the EIP method is defined by:
\begin{itemize}
	\item[1.] Compute $\widehat{y}_{n+1}$ such that $\widehat{y}_{n+1}=\Phi_h(y_n)$;
	\item[2.] Compute the value of vector $\widehat{\lambda}\in \mathbb{R}^l$ by
	\begin{equation}\label{2-1}
		\widehat{\lambda}=-\left(\nabla g(\widehat{y}_{n+1})^\top \nabla g(\widehat{y}_{n+1})\right)^{-1}g(\widehat{y}_{n+1}),~~\mbox{with}~~\nabla g(\widehat{y}_{n+1})\in\mathbb{R}^{d\times l};
	\end{equation}
	\item[3.] Update $y_{n+1}=\widehat{y}_{n+1}+\nabla g(\widehat{y}_{n+1})\widehat{\lambda}$.
\end{itemize}

\begin{prop}\label{prop-1}
	The EIP method is equivalent to the standard orthogonal projection approach \cite{hlw06} for system \eqref{eq-1-1} with only one Newton iteration to solve the corresponding nonlinear system.
\end{prop}

\begin{prf}
	
	The algorithm of the standard orthogonal projection approach can be written as follows:
	
	\begin{itemize}
		\item[1.] Compute $\widehat{y}_{n+1}$ by $\widehat{y}_{n+1}=\Phi_h(y_n)$;
		\item[2.] Project the value $\widehat{y}_{n+1}$ onto the manifold $\mathcal{M}$ through
		\begin{equation}
			y_{n+1}=\widehat{y}_{n+1}+\nabla g(\widehat{y}_{n+1})\lambda\quad \mbox{and}\quad g(y_{n+1})=0,
		\end{equation}
	\end{itemize}
	where $\nabla g(\widehat{y}_{n+1})$ is the projection direction. At each time step, we have to solve a nonlinear system of $\lambda$ as
	\begin{equation}\label{2-3}
		g\left(\widehat{y}_{n+1}+\nabla g(\widehat{y}_{n+1})\lambda\right)=0.
	\end{equation}
	Denote $F(\lambda)=g\left(\widehat{y}_{n+1}+\nabla g(\widehat{y}_{n+1})\lambda\right)$, then a Newton iteration for \eqref{2-3} yields
	\begin{equation}\label{2-4} \lambda_{k+1}=\lambda_k-\big[\nabla F(\lambda_{k})\big]^{-1}F(\lambda_{k}),~~k=0,1,2,\cdots.
	\end{equation}
	Since $\lambda$ is actually of small magnitude as stated in \cite{hlw06} and will be further proved below, we let the starting value $\lambda_0=0$, then one Newton iteration step of \eqref{2-4} yields the formula \eqref{2-1} exactly. \qed
	
\end{prf}

\subsection{Accuracy and invariants preservation}

Although the EIP method is a further simplification of the standard orthogonal projection approach, we will show that its accuracy is still the same as the underlying RK method. Moreover, the invariant error can be significantly reduced to a magnitude of $\mathcal{O}(h^{2(p+1)})$, that is, utilizing a high-order RK method the EIP method can preserve the invariant to the machine accuracy in practical computation.

As stated in Proposition \ref{prop-1}, the EIP method is related to the standard projection approach. Therefore, we can analyze the EIP method from the aspect of the simplified projection approach, specifically regarding to the solution of the parameter $\lambda$ which is the key difference between the two kinds of methods.

\begin{lem}\label{lem-1}
	Assume $\nabla g(y_n)$ has full column rank, there exists $h_\star>0$ such that the nonlinear system \eqref{2-3} for the standard projection approach has a unique solution $\lambda_\star=\lambda_\star(h)$ for $h\in[0,h_\star]$ and  $\lambda_\star$ is of size $\mathcal{O}(h^{p+1})$. Therefore, the standard projection approach is of order $p$.
\end{lem}

\begin{prf}
	For simplicity, we denote
	\[
	\widehat{F}(\lambda,h)=g\left(\widehat{y}_{n+1}+\nabla g(\widehat{y}_{n+1})\lambda\right)=0.
	\]
	Notice that
	\[
	\widehat{F}(0,0)=g(y_n)=0,~~\mbox{and}~~\nabla_\lambda\widehat{F}(0,0)=\nabla g(y_n)^\top \nabla g(y_n)
	\]
	is invertible as $\nabla g(y_n)$ has full rank, then the implicit function theorem ensures the existence of a neighborhood $[0,h_\star]$ and a unique smooth solution $\lambda_\star=\lambda_\star(h)$ such that $\lambda_\star(0)=0$ and $\widehat{F}\big(\lambda_\star(h),h\big)=0$ for all $h\in[0,h_\star]$. Moreover, we can expand
	\[
	\widehat{F}(\lambda_\star,h)=\widehat{F}(0,h)+\nabla_\lambda\widehat{F}(0,h)\lambda_\star+\mathcal{O}(\lambda_\star^2),
	\]
	with
	\[
	\widehat{F}(0,h)=g(\widehat{y}_{n+1})=\mathcal{O}(h^{p+1}),\quad	\nabla_\lambda\widehat{F}(0,h)=\nabla_\lambda\widehat{F}(0,0)+\mathcal{O}(h).
	\]
	Then, we can derive $\lambda_\star=\lambda_\star(h)=\mathcal{O}(h^{p+1})$, which implies the standard projection approach is of order $p$.\qed
	
\end{prf}

\begin{rmk}
	The assumption that $\nabla g(y_n)$ has full column rank can be satisfied as long as the $l$ invariants $g_1, g_2,\cdots, g_l$ are functional independent, and thus the results of Lemma \ref{lem-1} hold for most cases of multiple invariants.
\end{rmk}

\begin{lem}\label{lem-2}
	Assume $\lambda_\star$ is the exact solution of the nonlinear system \eqref{2-3}, then the $k$-th iteration solution $\lambda_k$ obtained by \eqref{2-4} satisfies $\lambda_{k}-\lambda_\star=\mathcal{O}(h^{2^k(p+1)})$. Moreover, we have $F(\lambda_k) = \mathcal{O}(h^{2^k(p+1)})$.
\end{lem}

\begin{prf}
	A direct calculation yields
	\[
	\begin{aligned}
	\lambda_{k+1}-\lambda_\star&=\lambda_k-\big[\nabla F(\lambda_k)\big]^{-1}F(\lambda_k)-\lambda_\star\\
	&=\lambda_k-\lambda_\star-\big[\nabla F(\lambda_k)\big]^{-1}\left(F(\lambda_\star)-\nabla F(\lambda_k)(\lambda_\star-\lambda_k)-\mathcal{O}(\|\lambda_\star-\lambda_k\|^2)\right)\\
	&=\lambda_k-\lambda_\star+\big[\nabla F(\lambda_k)\big]^{-1}\left(\nabla F(\lambda_k)(\lambda_\star-\lambda_k)+\mathcal{O}(\|\lambda_\star-\lambda_k\|^2)\right)\\
	&=\mathcal{O}(\|\lambda_\star-\lambda_k\|^2),
	\end{aligned}
	\]
	where we have used $F(\lambda_\star)=0$ and $\|\cdot\|$ is the Euclidean norm. Let the starting value $\lambda_0=0$, we can obtain $\lambda_1-\lambda_\star=\mathcal{O}(\|\lambda_\star\|^2)=\mathcal{O}(h^{2(p+1)})$ by Lemma \ref{lem-1}. Thus we can recursively derive $\lambda_{k}-\lambda_\star=\mathcal{O}(h^{2^k(p+1)})$, which leads to $F(\lambda_k) = \mathcal{O}(h^{2^k(p+1)})$. This completes the proof.\qed
\end{prf}

As a direct consequence of Lemma \ref{lem-2}, we have the following main results about the accuracy and invariant error of the EIP method.

\begin{thm}\label{thm-2}
	The EIP method retains the same order of accuracy as the underlying RK method, and the invariant error is of size $\mathcal{O}(h^{2(p+1)})$.
\end{thm}

\begin{prf}
	
	From Lemma \ref{lem-2}, the parameter $\widehat{\lambda}$ in the EIP method satisfies $\widehat{\lambda}=\lambda_1=\lambda_\star+\mathcal{O}(h^{2(p+1)})=\mathcal{O}(h^{p+1})$.
	Since $\Phi_h$ corresponds to a RK method of order $p$, by Taylor's Theorem we have $\widehat{y}_{n+1}=y(t_{n+1})+\mathcal{O}(h^{p+1})$ and
	\[
	\nabla g(\widehat{y}_{n+1})=\nabla g(y(t_{n+1}))+\mathcal{O}(h^{p+1}).
	\]
	Hence the local error of the EIP method yields
	\[
	y_{n+1}=\widehat{y}_{n+1}+\nabla g(\widehat{y}_{n+1})\widehat{\lambda}=y(t_{n+1})+\mathcal{O}(h^{p+1}),
	\]
	that is, the EIP method has an accuracy with the same order as the underlying RK method. While the invariant error is a straightforward result of Lemma \ref{lem-2} by setting $k=1$.\qed

\end{prf}

\begin{rmk}
	Given a time step $h$, if we choose a high-order underlying RK method in the first step of the EIP method, the invariant can reach to the round-off error of double precision machines. In this sense, we call it a practically invariant-preserving method.
\end{rmk}

\begin{rmk}
	Due to the equivalence between the EIP method and the projection method, $\nabla g(\widehat{y}_{n+1})$ in the last step of the EIP method actually corresponds to a projection direction.  Therefore, we can replace this term with any other suitable projection directions. Of course, the steepest descent direction for convergence is $\nabla g(y_{n+1})$. However, to ensure that the resulting method is fully explicit, the direction can only involve some known solutions, such as $y_n$ and $\widehat{y}_{n+1}$. For example, if we take $\nabla g(y_n)$ or $\nabla g\left(\frac{1}{2}(y_n+\widehat{y}_{n+1})\right)$ as the projection direction, we can also obtain the same results about the accuracy and the invariant error. Meanwhile, it is not necessary to restrict the direction to forms of $\nabla g(y)$. Instead, any other effective projection directions can be utilized to construct the EIP method. However, among those directions, which one has the fastest convergence turns into an optimization problem and is worthy of further study.
\end{rmk}

To illustrate the above results about the accuracy and invariant error of the EIP method, we consider the harmonic oscillator as a concrete example
\[
\dot{y}=\left(\begin{array}{cc}
0 & \omega \\
-\omega & 0
\end{array} \right)y,
\]
where $y=(p,q)^\top$ and by definition $\omega>0$. This is a Hamiltonian system with quadratic Hamiltonian $H(y)=\frac{\omega}{2} y^\top y$. For a given initial data $y_0$, the invariant in this example becomes $g(y)=H(y)-H(y_0)$. After one step of a $p$-th order RK method $y_n\rightarrow \widehat{y}_{n+1}$, we first apply the standard projection approach and obtain the corresponding nonlinear system \eqref{2-3} which is actually a quadratic equation of $\lambda$
\[
a_2\lambda^2+2a_1\lambda+a_0=H(y_0),
\]
where $a_k=\frac{\omega^{k+1}}{2}\widehat{y}_{n+1}^\top\widehat{y}_{n+1}$ and thereby can be solved analytically, which is exactly the basic ideas in \cite{calvo-06-explicit-RK-quadratic-invariant,zhang-20-IEQ-explicit-RK,jiang-20-IEQ-explicit-RK} for quadratic invariant. Of the two solutions that exist for $\lambda$ we choose the one with the smallest absolute value since we expect $\lambda$ to be close to zero for small $h$. Hence, we get
\begin{equation}\label{e-1}
\lambda_\star=-\frac{1}{\omega}\left(1-\frac{\|y_0\|}{\|\widehat{y}_{n+1}\|}\right)=\mathcal{O}(h^{p+1}).
\end{equation}
While for the EIP method, $\widehat{\lambda}$ can be directly calculated by
\[
\widehat{\lambda}=-\frac{1}{2\omega}\left(1-\frac{\|y_0\|^2}{\|\widehat{y}_{n+1}\|^2}\right)=\frac{1}{2}\left(1+\frac{\|y_0\|}{\|\widehat{y}_{n+1}\|}\right)\lambda_\star=\lambda_\star+\mathcal{O}(h^{2(p+1)}).
\]
Obviously, the local errors of $\widehat{\lambda}$ and $\lambda_\star$ coincide with the theoretical results in Lemma \ref{lem-1} and \ref{lem-2}, respectively. Subsequently, we can also derive the invariant error as $g(y_{n+1})=\mathcal{O}(h^{2(p+1)})$ by a direct substitution and Taylor's series expansion.

We numerically simulate the harmonic oscillator by the EIP methods with the underlying RK methods being the classic explicit ones of order 1 to 4, whose Butcher tabular are listed as follows:
\[
\begin{array}{cc}
& \\
& \\
& \\
\multicolumn{1}{c|}{0}	&  \\ \hline
\multicolumn{1}{c|}{}& 1
\end{array} \qquad
\begin{array}{ccc}
& &\\
& &\\
\multicolumn{1}{c|}{0}&  &  \\
\multicolumn{1}{c|}{1/2}	& 1/2 &  \\ \hline
\multicolumn{1}{c|}{}& 0 & 1
\end{array} \qquad
\begin{array}{cccc}
& & & \\
\multicolumn{1}{c|}{0}	&  &  &  \\
\multicolumn{1}{c|}{1/3}	& 1/3 &  &  \\
\multicolumn{1}{c|}{2/3}	& 0 & 2/3 &  \\ \hline
\multicolumn{1}{c|}{}& 1/4 & 0 & 3/4
\end{array} \qquad
\begin{array}{ccccc}
\multicolumn{1}{c|}{0}&  &  &  &  \\
\multicolumn{1}{c|}{1/2}& 1/2 &  &  &  \\
\multicolumn{1}{c|}{1/2}& 0 & 1/2 &  &  \\
\multicolumn{1}{c|}{1}& 0 & 0 & 1 &  \\ \hline
\multicolumn{1}{c|}{}& 1/6 & 2/6 & 2/6 & 1/6
\end{array}.
\]
For convenience, we denote the above four RK methods as RK1, RK2, RK3, RK4, respectively. In Table.~\ref{ex0-tab1}, the errors of $|\widehat{\lambda}-\lambda_\star|$ and the corresponding orders are presented. For RK1 and RK3, the results are coincided with the analytical derivation, whereas for RK2 and RK4, we can observe a surprising result that the convergence order becomes $2(p+2)$ instead of $2(p+1)$ which is mainly because the simulated model a linear equation. As a consequence, the related energy errors also behave a similar result in Table.~\ref{ex0-tab2}, where the results labeled iter\#1 correspond to the EIP method that can be viewed as a projection method with one step iteration. From Lemma~\ref{lem-2}, the invariant error of the standard projection method is of size $\mathcal{O}(h^{2^k(p+1)})$ with $k$ being the iteration number for the nonlinear system \eqref{2-3}. Therefore, we also
provide results with $k=2$ for RK1 and RK2, which again conform the theoretical analysis. While for RK3 and RK4, iteration twice will make the energy errors reach the machine accuracy quickly and consequently we can hardly get the convergence orders.

\begin{table}[H]
	\renewcommand\arraystretch{1.2}
	\centering
	\begin{tabular*}{0.9\textwidth}[h]{@{\extracolsep{\fill}}ccccccccc} \hline
		\multirow{2}{*}{$h$}	& \multicolumn{2}{c}{RK1} &\multicolumn{2}{c}{RK2} & \multicolumn{2}{c}{RK3}  & \multicolumn{2}{c}{RK4}\\ \cline{2-3}\cline{4-5}\cline{6-7}\cline{8-9}
		& error & order & error & order & error & order & error & order \\ \hline
		0.1/1 & 6.3506e-03 &  & 6.2524e-04 &  &  4.0849e-05 &  & 1.8688e-06 & \\[1ex]
		0.1/2 & 6.2524e-04  & 3.3444   &   3.0048e-06 &  7.7010 &   2.8629e-07 & 7.1567 & 5.5259e-10 & 11.7237\\[1ex]
		0.1/4 & 4.5956e-05  & 3.7661 &    1.1909e-08 &    7.9790 &  1.2703e-09 & 7.8161 &  1.4148e-13 &  11.9314 \\[1ex]
		0.1/8 & 3.0048e-06    &  3.9349  &  4.6563e-11 &    7.9987 &    5.1204e-12 &  7.9547 &  3.3307e-17 & 12.0525\\ [1ex] \hline
	\end{tabular*}
	\caption{Errors of $|\widehat{\lambda}-\lambda_\star|$ and the corresponding orders for EIP methods with different underlying RK methods for the harmonic oscillator problem with $\omega=10$, $y_0=[1,0]^\top$ at $t=1$.}\label{ex0-tab1}
\end{table}

\begin{table}[H]
	\renewcommand\arraystretch{1.2}
	\centering
	\begin{tabular*}{0.9\textwidth}[h]{@{\extracolsep{\fill}}ccccccc} \hline
		& & & $h_0/1$ & $h_0/2$ & $h_0/4$ & $h_0/8$\\\hline
		\multirow{4}{*}{RK1}&\multirow{2}{*}{iter\#1} & error & 1.0354 & 7.0644e-02 &4.7404e-03  & 3.0283e-04\\
		& & order & - & 3.8735 & 3.8975 & 3.9684\\\cline{2-7}
	    & \multirow{2}{*}{iter\#2} & error&  1.7712e-02 &1.9303e-04 & 1.0550e-06 & 4.5142e-09\\
		& & order & - & 6.5197 &7.5154  &  7.8686\\\hline
		\multirow{4}{*}{RK2} & \multirow{2}{*}{iter\#1} &error & 3.1922 & 7.0644e-02 & 3.0283e-04 & 1.1915e-06\\
		& &order & -& 8.8198 & 7.8659 & 7.9896\\\cline{2-7}
		&\multirow{2}{*}{iter\#2}  &error &5.6576e-01  &1.9303e-04 & 4.5142e-09 &7.1054e-14\\
		& &order &- & 11.5171 & 15.3840 & 15.9552\\\hline
		\multirow{2}{*}{RK3}& \multirow{2}{*}{iter\#1} & error & 2.1230e-01 &  3.9722e-03 &2.8561e-05 & 1.2701e-07\\
		& & order& -& 5.7400&7.1198 & 7.8129\\\hline
		\multirow{2}{*}{RK4}& \multirow{2}{*}{iter\#1} & error &3.4710e-01 &1.8575e-04 &5.5253e-08 &1.4149e-11\\
		& & order&- & 10.8678&11.7150 &11.9312\\
		 [1ex] \hline	
	\end{tabular*}
	\caption{Energy errors and the corresponding orders for EIP methods with different underlying RK methods for the harmonic oscillator problem with $\omega=10$, $y_0=[1,0]^\top$ till $t=1$. $h_0=0.1$ for RK1 and $h_0=0.2$ for the rest.}\label{ex0-tab2}
\end{table}

\begin{rmk}
	According to the above simple example, we actually have illustrate the key strategy in \cite{calvo-06-explicit-RK-quadratic-invariant,zhang-20-IEQ-explicit-RK,jiang-20-IEQ-explicit-RK} to make their methods explicit for quadratic invariant or general invariant after quadratization. Nevertheless, this strategy is only feasible for preservation of single invariant because there will be no analytical solutions of vector-valued $\lambda$ for multiple invariants.
	
	Although the standard projection method can handle multiple invariants, it inevitably require iterations to solve the nonlinear system \eqref{2-3} for vector-valued $\lambda$, which makes the standard projection method less efficient than the EIP method.
\end{rmk}

\subsection{Implementation issue}

As demonstrated in the algorithm of the EIP method, the implementation is completely explicit and the main effort is paid to solve the parameter $\lambda$. In this subsection, with respect to single and multiple invariants, we will provide the solution of $\widehat{\lambda}$ in details for both ODEs and PDEs circumstances.

\begin{itemize}[leftmargin=*]
	\item Case I : single invariant
\end{itemize}

Suppose we have obtained $\widehat{y}_{n+1}$ through the explicit RK method $\Phi_h$. For ODEs or one-dimensional PDEs after semi-discretization, $g(y):\mathbb{R}^N\rightarrow\mathbb{R}$ and $\widehat{\lambda}$ is a scalar. Here, we use $N$ as the notation of dimension instead of the aforementioned $d$ for convenience. Consequently, the expression of $\widehat{\lambda}$ \eqref{2-1} can be simplified as
\begin{equation}\label{2-2-1}
\widehat{\lambda}=-\dfrac{g(\widehat{y})}{\|\nabla g(\widehat{y})\|^2},
\end{equation}
where we have omitted the subscript of $\widehat{y}_{n+1}$ for simplicity. Specifically, when considering the canonical Hamiltonian system and let the Hamiltonian energy $H(y)$ be the targeted invariant, the formula \eqref{2-2-1} becomes
\begin{equation}\label{2-2-2}
\widehat{\lambda}=-\dfrac{H(\widehat{y})-H(y_0)}{\|\nabla H(\widehat{y})\|^2},
\end{equation}
and subsequently $y=\widehat{y}+\lambda\nabla H(\widehat{y})$.

\begin{rmk}
	In Ref. \cite{brugnano-12-EP}, the authors proposed a correction technique to prevent the accumulation of round-off errors for the HBVMs. Actually, such technique belongs to the framework of the EIP method with the parameter calculated exactly as \eqref{2-2-2} and the HBVMs as the one-step method $\Phi_h$. As analyzed in the above section, the EIP method can make a non-conservative RK method invariants-preserving. Therefore, the additional employment of the EIP method as a correction technique for the HBVMs may confuse the performance in the preservation of invariants. Nevertheless, for general conservative methods, the accumulation of round-off errors often renders the invariant errors growing linearly. How to design a technique to prevent it is still worthy of further study.
\end{rmk}

Due to the explicit form, the EIP method will be more competitive for PDEs of high dimensions. However, it is not straightforward to extend the formula \eqref{2-1} for such cases, even if only one invariant needs to be preserved. Consider a two-dimensional problem with $\widehat{y}_{n+1}\in\mathbb{R}^{N\times N}$ for example, where $N$ represents the grid size in both $x$ and $y$ directions, then the gradient of single invariant $\nabla g(\widehat{y}_{n+1})$ is of dimension $N\times N$ in general. A direct substitution into the formula \eqref{2-1} will suffer the disagreement of dimension. To overcome this obstacle, we define a function ``{vec}" that represents the vectorization of matrix, i.e., vec($A_{N\times N}$) results a vector of dimension $N^2\times 1$, composed of the elements of $A$. Then we can obtain
\begin{equation}\label{2-2-3}
\widehat{\lambda}=-\dfrac{g(\widehat{y})}{\|\mbox{vec}(\nabla g(\widehat{y}))\|^2},
\end{equation}
which can be generalized to three-dimensional cases in a straightforward manner.

\begin{itemize}[leftmargin=*]
	\item Case II : multiple invariants
\end{itemize}

For ODEs and one-dimensional semi-discrete PDEs with $l$ invariants, the express of $\lambda$ retains the same as that in \eqref{2-1} with $\nabla g(\widehat{y})\in\mathbb{R}^{N\times l}$. Different from \eqref{2-2-1} and \eqref{2-2-3} for single invariant where the denominators are scalars, one has to compute the inverse of a matrix of dimension $l\times l$. For high-dimensional PDEs after semi-discretization, the technique in \eqref{2-2-3} should also be employed, and the calculation of $\widehat{\lambda}$ yields
\begin{equation}\label{2-2-4}
\widehat{\lambda}=-\left(\mathcal{G}^\top\mathcal{G}\right)^{-1}g(\widehat{y}_{n+1}),
\end{equation}
where $\mathcal{G}=\big(\mbox{vec}\left(\nabla g_1(\widehat{y})\right), \mbox{vec}\left(\nabla g_2(\widehat{y})\right),\cdots, \mbox{vec}\left(\nabla g_l(\widehat{y})\right)\big)\in\mathbb{R}^{N^2\times l}$ for two-dimensional problems. Accordingly, the implementation of the EIP method for multiple invariants is extremely simple, whereas most of existing conservative methods fail to simultaneously preserve two or more invariants.

We define the inner product and the corresponding norms for 2D problems as follows:
\[
(\bm u,\bm v)=h_xh_y\sum_{i=1}^{N_x}\sum_{j=1}^{N_y}\bm u_{jk}\bm u_{jk},\quad \|\bm u\|_h=(\bm u,\bm u)^{\frac{1}{2}}.
\]
Similar definitions can be obtained for 1D and 3D cases.

\section{Numerical experiments}

In this section,  we will present various numerical examples for both conservative ODEs and PDEs to verify the theoretical results and demonstrate the efficiency of the EIP method. The underlying one step method is uniformly taken as RK4 unless otherwise stated. All the experiments are performed in Matlab, version 2016b, on a computer equipped with the Intel Core i7-6700 processor running at 3.40GHz. For the part of accuracy tests, when the exact solution is unknown, we take the following strategy to compute the convergence order.

Consider a $p$th-order numerical approximation of an exact solution $y$ (a scalar variable for simplicity). The approximation depends on a small parameter $h$, either as the grid size or time step, and we denote it by $y_h$. By the $h$-expansion \cite{hnw08} of the global error, we have
\[
y_h=y+Ch^p+\mathcal{O}(h^{p+1}),
\]
where the number $C$ is independent of $h$ and typically depends on the exact solution. Next, we halve the step size $h$ successively and obtain
\[
y_{h/2}=y+C(h/2)^p+\mathcal{O}(h^{p+1}),\quad \mbox{and}\quad y_{h/4}=y+C(h/4)^p+\mathcal{O}(h^{p+1}).
\]
The ratio of the difference between the above three numerical solutions is calculated by
\begin{equation}\label{accuracy}
\frac{y_h-y_{h/2}}{y_{h/2}-y_{h/4}}=\frac{1-2^{-p}+O(h)}{2^{-p}-2^{-2p}+O(h)}=2^p+\mathcal{O}(h).
\end{equation}
Hence, after continually decreasing $h$ and compute the logarithm of the left side of \eqref{accuracy}, we can  get an estimate of $p$. For vector-valued solutions, the difference should be measured under some appropriate norms.

\subsection{Experiment I: perturbed Kepler system}

We first consider the perturbed Kepler system with Hamiltonian
\begin{equation}\label{ex1-eq1}
	H(q_1,q_2,p_1,p_2)=\frac{1}{2}(p_1^2+p_2^2)-\frac{1}{\sqrt{q_1^2+q_2^2}}-\frac{0.005}{2\sqrt{(q_1^2+q_2^2)^3}},
\end{equation}
which describes the motion of a planet in the Schwarzschild potential for Einstein's general relatively theory. Besides the above Hamiltonian, this problem has another invariant, i.e., the angular momentum
\begin{equation}\label{ex1-eq2}
L(q_1,q_2,p_1,p_2)=q_1p_2-q_2p_1.
\end{equation}
The most of existing invariant-preserving methods usually preserve only one of the above two invariants. For examples, the AVF method can be used to preserve the Hamiltonian \eqref{ex1-eq1} while the symplectic partitioned RK method can handle with the angular momentum \eqref{ex1-eq2} in general. In this experiment, we will not only verify the effectiveness of the proposed method through the accuracy test, but also show its flexibility in the conservation of either one or two invariants.

Let the initial conditions $q_1(0)=1-e$, $q_2(0)=0$, $p_1(0)=0$, $p_2(0)=\sqrt{(1+e)/(1-e)}$ where $e$ represents the eccentricity and is taken as $e=0.6$. Although one can derive the exact solution of the perturbed Kepler system based on the two invariants \cite{hlw06}, we here choose the formula \eqref{accuracy} instead to compute the convergence order of the EIP methods. Denote \textbf{EIP-H}, \textbf{EIP-L} and \textbf{EIP-HL} as the methods designed to preserve $H$, $L$ and both, respectively. The results of accuracy tests for those methods are listed in Table.~\ref{ex1-tab1}, which uniformly show a convergence order of 4 as expected.

\begin{table}[H]
	\renewcommand\arraystretch{1.2}
	\centering
	\begin{tabular*}{0.9\textwidth}[h]{@{\extracolsep{\fill}}ccccccc} \hline
		\multirow{2}{*}{$h$}	& \multicolumn{2}{c}{\textbf{EIP-H}} &\multicolumn{2}{c}{\textbf{EIP-L}} & \multicolumn{2}{c}{\textbf{EIP-HL}} \\ \cline{2-3}\cline{4-5}\cline{6-7}
		& error & order & error & order & error & order \\ \hline
		0.02/1  &  &  &  & & & \\
		0.02/2 & 1.3935e-06  &    &   2.2607e-06 &   &    4.7011e-07 &  \\
		0.02/4 &  8.6286e-08  & 4.0134 &    1.3647e-07 &    4.0501 &  3.2365e-08 & 3.8605 \\
		0.02/8 & 5.3636e-09    &  4.0079  &  8.3757e-09 &   4.0262 &   2.1111e-09 &  3.9384 \\
		0.02/16 &  3.3424e-10    &  4.0042 &   5.1865e-10 &  4.0134 &   1.3462e-10 &   3.9710 \\  \hline			
	\end{tabular*}
\caption{Errors between two adjacent time steps and the corresponding orders for different EIP methods.}\label{ex1-tab1}
\end{table}

Different from the harmonic oscillator model, the Kepler system is fully nonlinear. We also present the convergence test for the invariant errors of different EIP methods to verify the result in Lemma~\ref{lem-2}. In Table~\ref{ex1-tab2}, due to the barrier of machine accuracy, we only give the corresponding errors and orders when the underlying methods are chosen as RK1 and RK2, which is clearly coincided with our analysis. In addition, there is no super convergence result that occurs for RK2 in the harmonic oscillator case.
\begin{table}[H]
	\renewcommand\arraystretch{1.2}
	\centering
	\begin{tabular*}{0.9\textwidth}[h]{@{\extracolsep{\fill}}cccccccc} \hline
	&	\multirow{2}{*}{$h$}	& \multicolumn{2}{c}{\textbf{EIP-H}} &\multicolumn{2}{c}{\textbf{EIP-L}} & \multicolumn{2}{c}{\textbf{EIP-HL}} \\ \cline{3-4}\cline{5-6}\cline{7-8}
	&	& error & order & error & order & error & order \\ \hline
\multirow{4}{*}{RK1}	&	0.03/1  & 7.6344e-08 &  & 1.1704e-07 & & 1.0251e-07& \\
	&	0.03/2 & 3.3388e-09  &  4.5151  &  7.3733e-09 & 3.9886  &    5.1644e-09 & 4.3110\\
	&	0.03/3 &  6.3275e-10 & 4.1021 &    1.5123e-09 &   3.9072 &  9.8071e-10 & 4.0972\\
	&	0.03/4 & 1.9642e-10    &  4.0664  &  4.8734e-10 &   3.9364 &   3.0466e-10 &   4.0638 \\ \hline	
	\multirow{4}{*}{RK2}	&	0.03/1  & 1.5467e-12 &   & 2.3863e-11 & &5.7380e-11 & \\
	&	0.03/2 &  1.9984e-14  &  6.2742  &   3.2674e-13 & 6.1905  &    7.7061e-13&  6.2184\\
	&	0.03/3 &  1.8874e-15  & 5.8198 &   2.9088e-14 &    5.9656 &  6.8501e-14 & 5.9693\\
	&	0.03/4 & 3.3307e-16    &  6.0296  &  5.1070e-15 &   6.0473 &   1.2434e-14 &   5.9314 \\ \hline			
	\end{tabular*}
	\caption{Energy errors and the corresponding orders for EIP methods with RK1 and RK2 as the underlying methods till $t=1$.}\label{ex1-tab2}
\end{table}

\begin{figure}[H]
\centering
\includegraphics[width=0.3\linewidth]{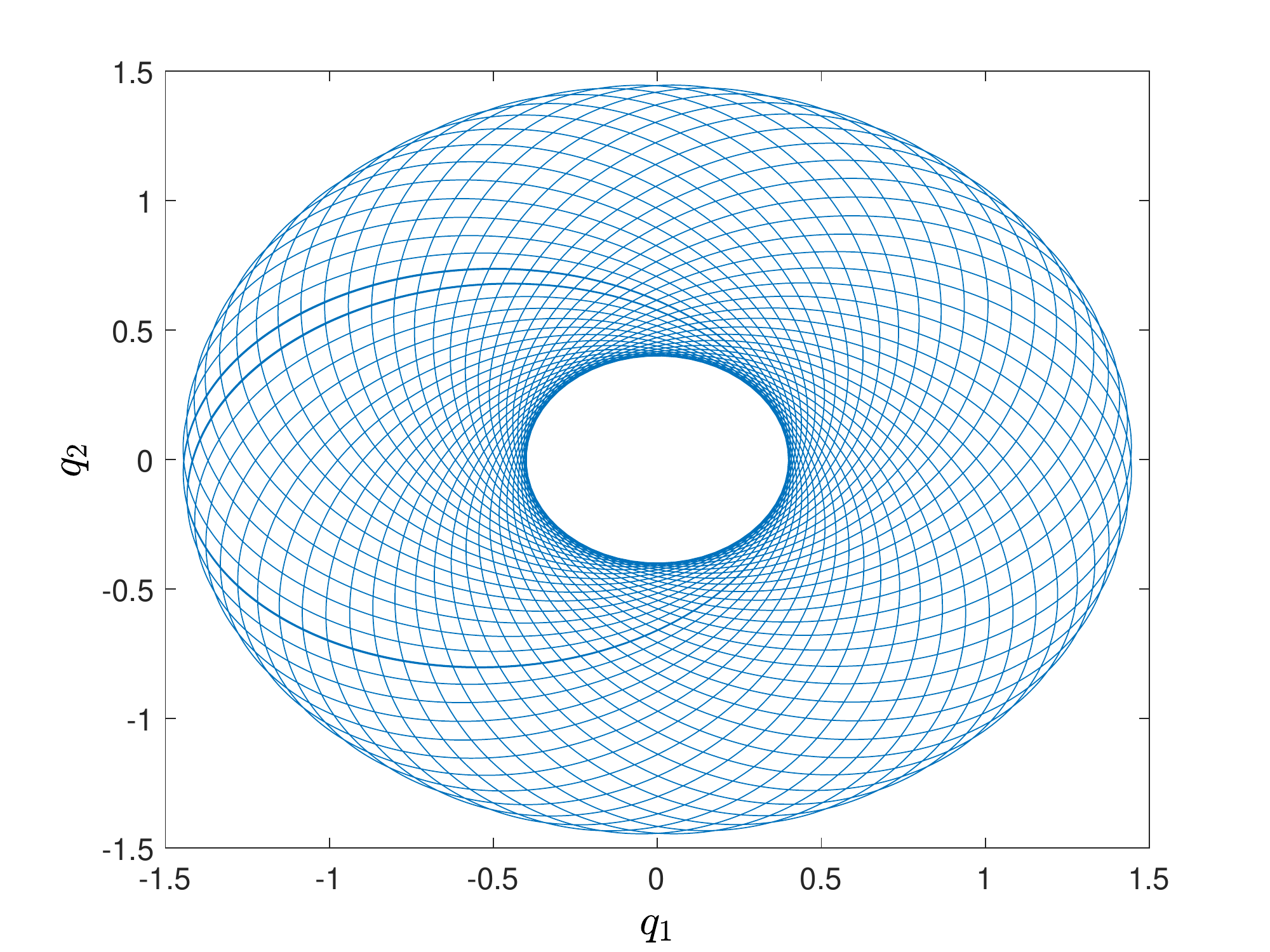}
\includegraphics[width=0.3\linewidth]{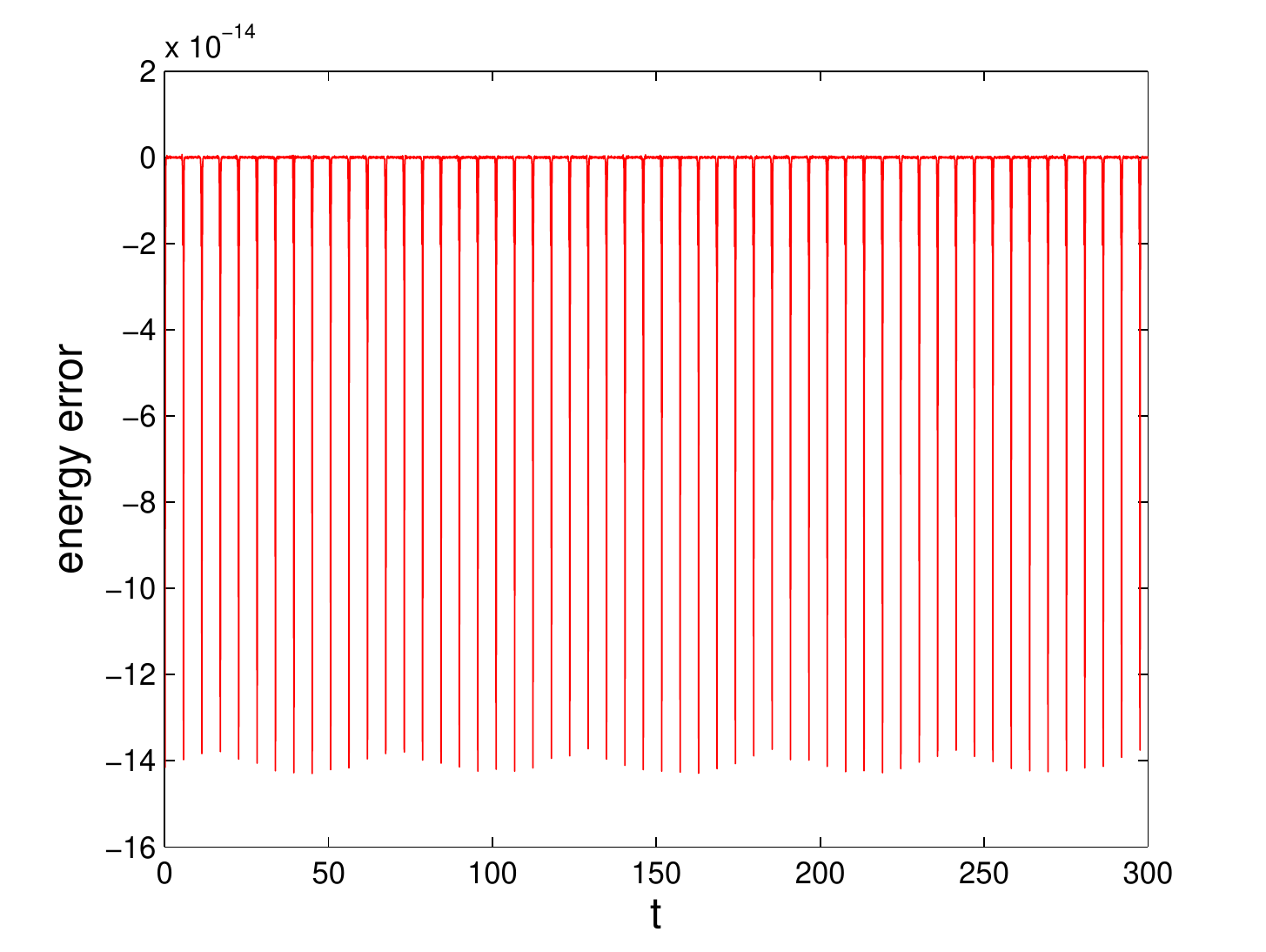}
\includegraphics[width=0.3\linewidth]{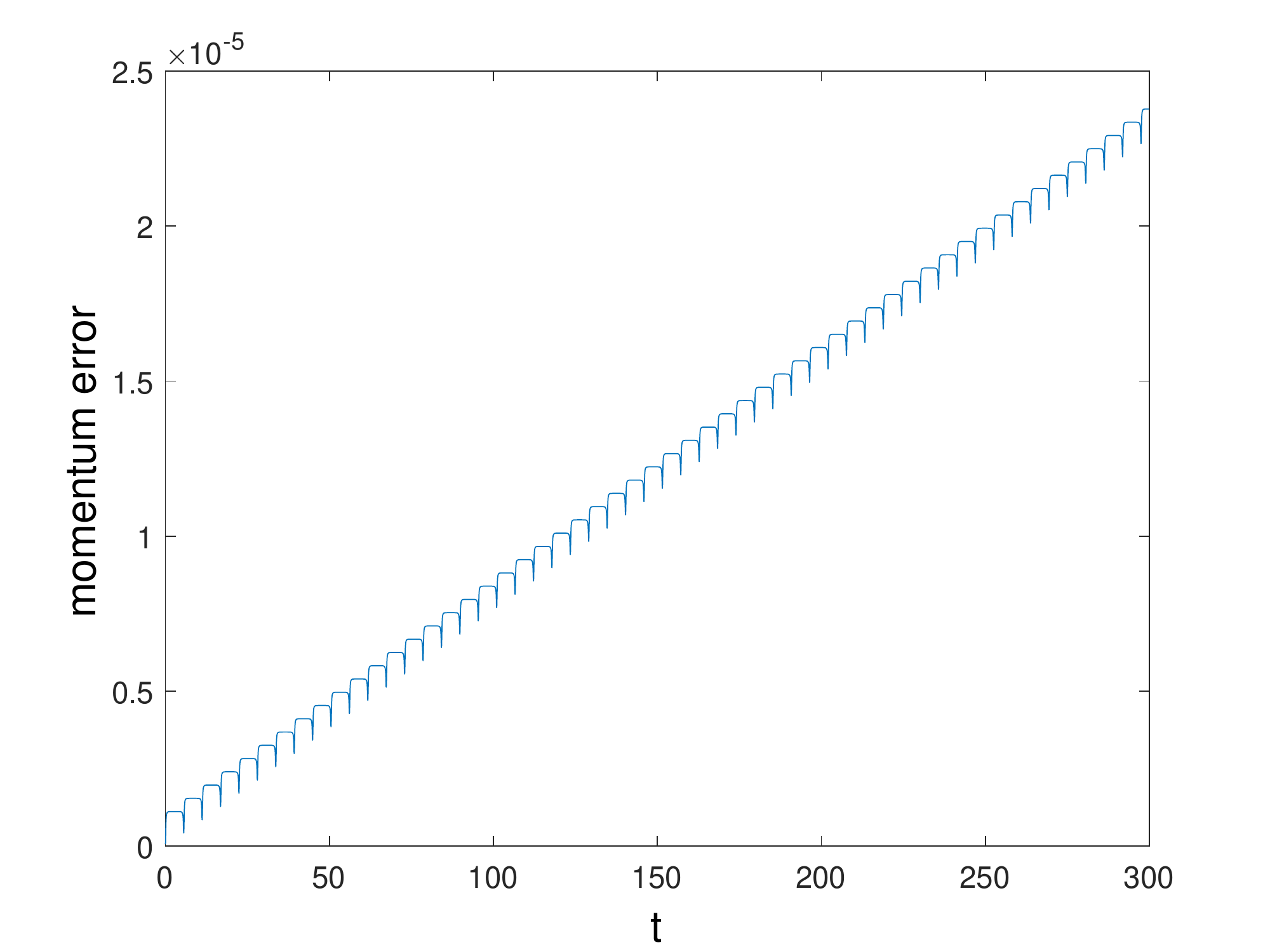}\\
\includegraphics[width=0.3\linewidth]{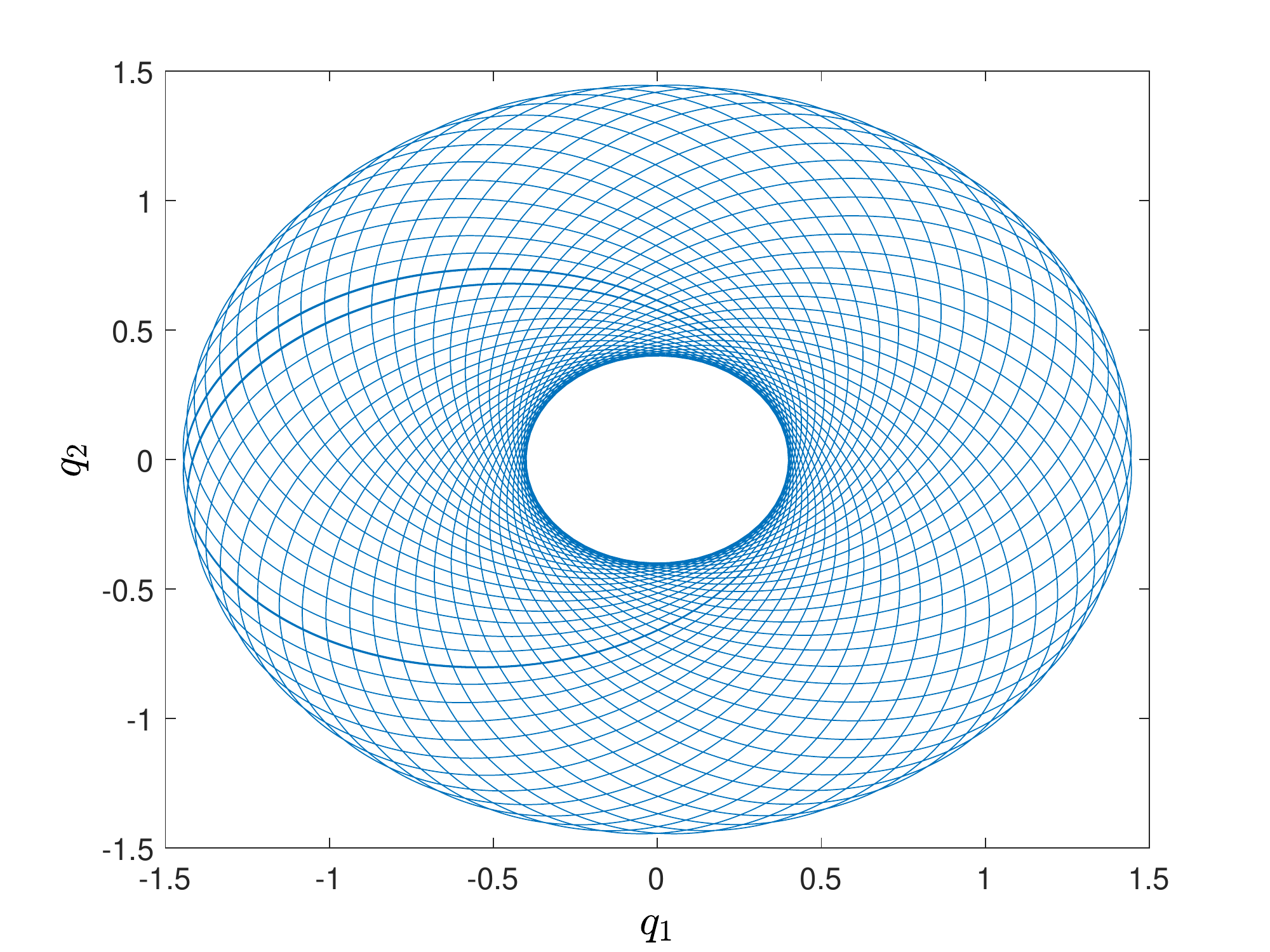}
\includegraphics[width=0.3\linewidth]{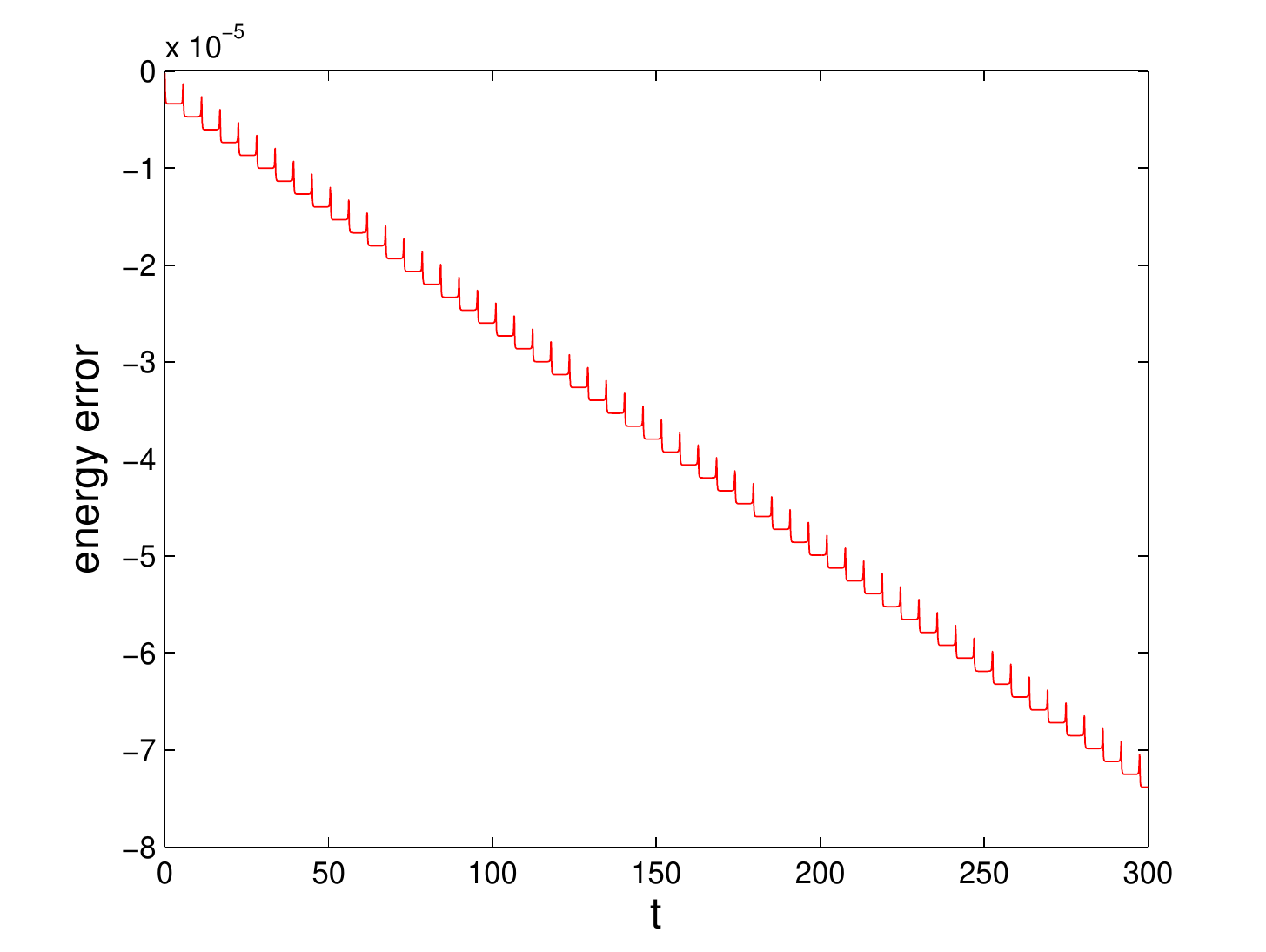}
\includegraphics[width=0.3\linewidth]{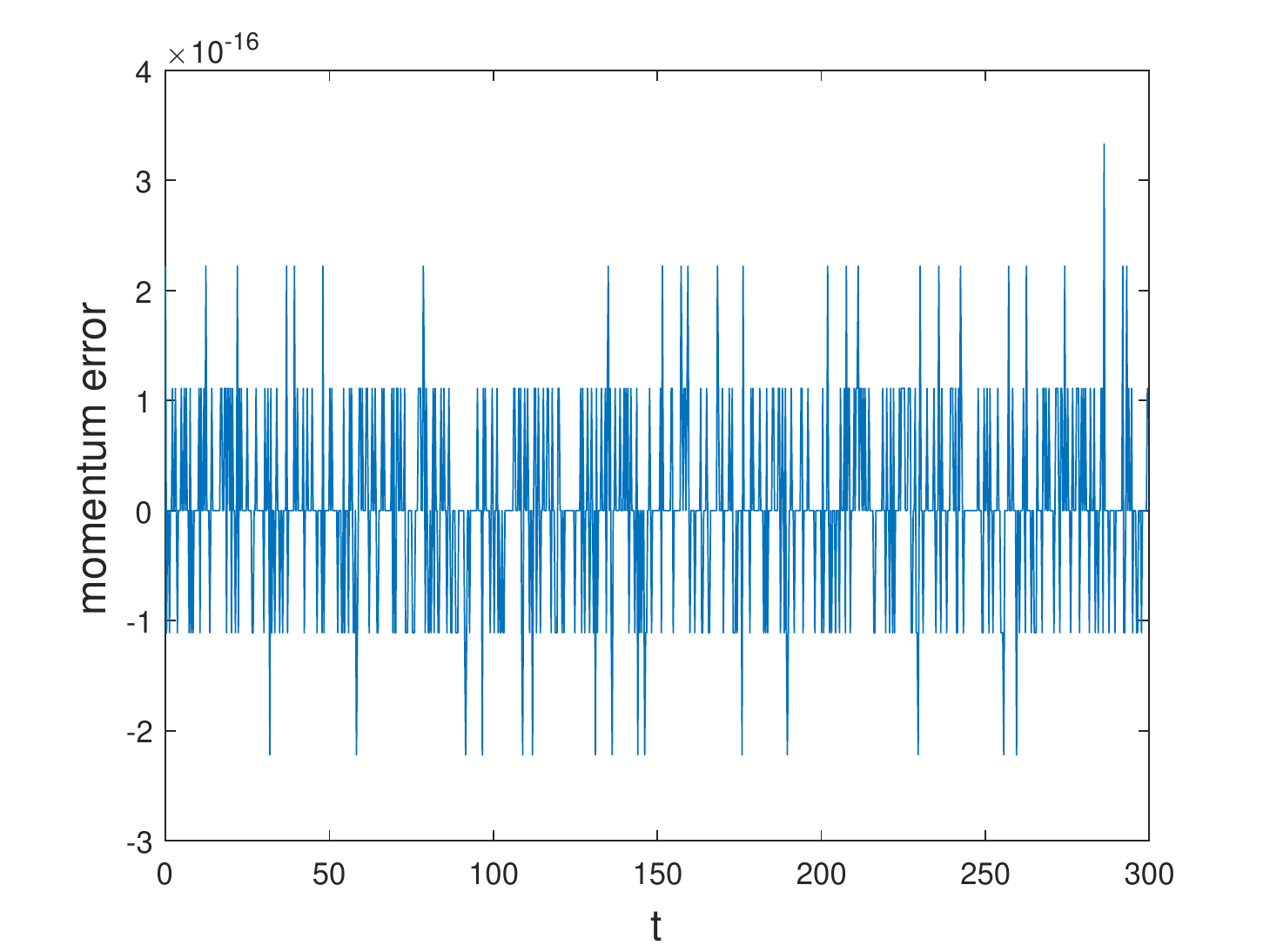}\\
\includegraphics[width=0.3\linewidth]{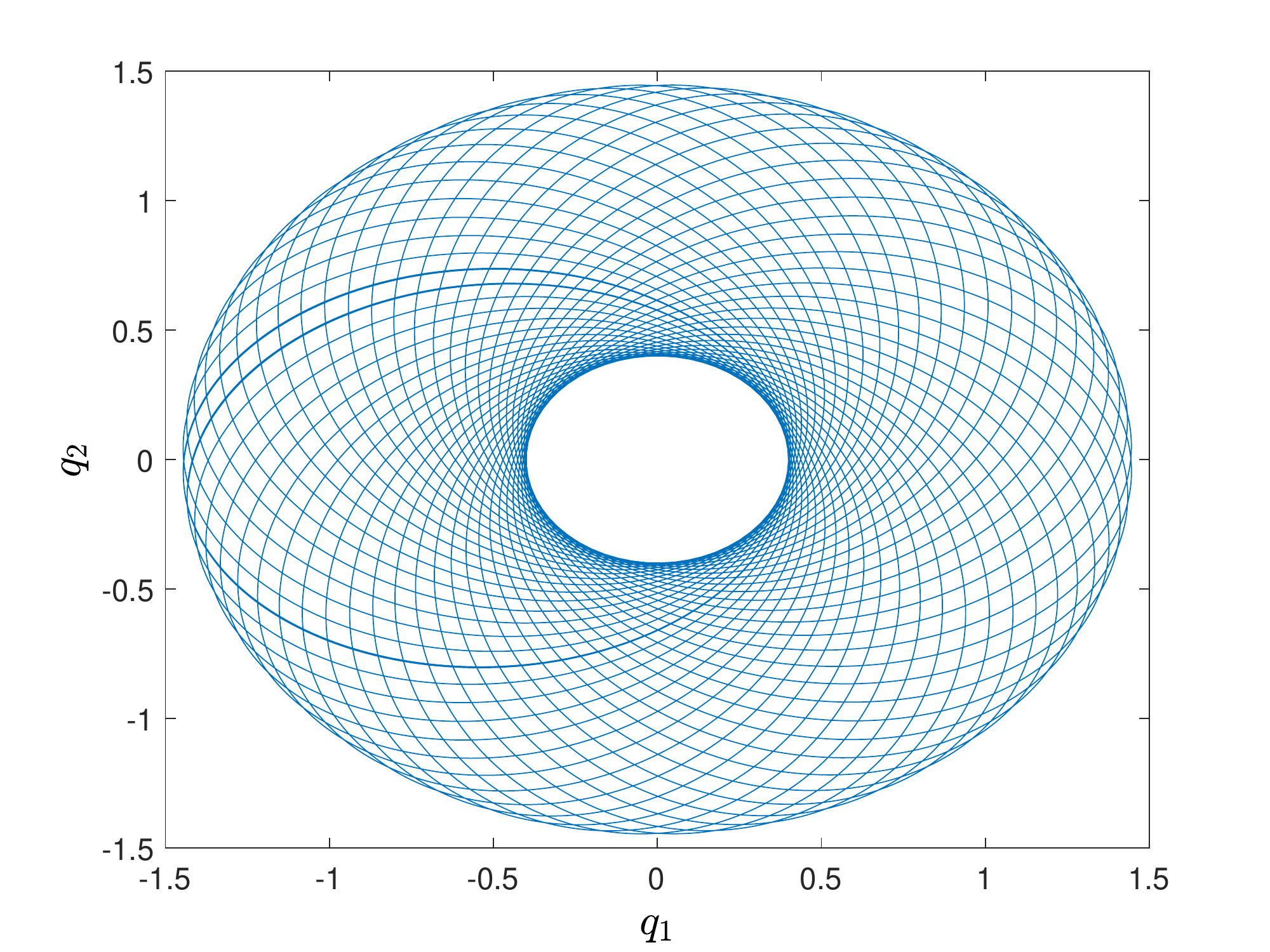}
\includegraphics[width=0.3\linewidth]{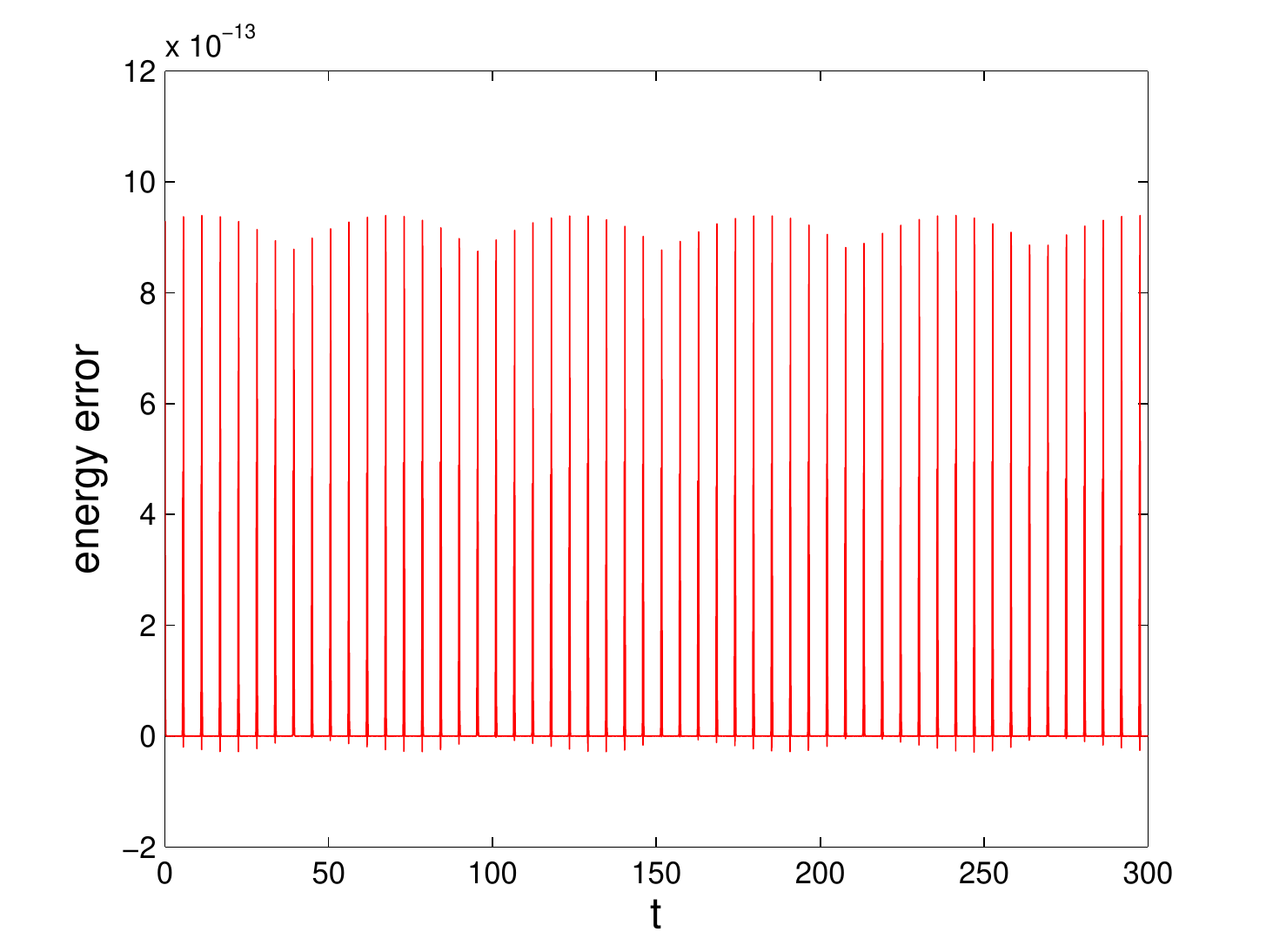}
\includegraphics[width=0.3\linewidth]{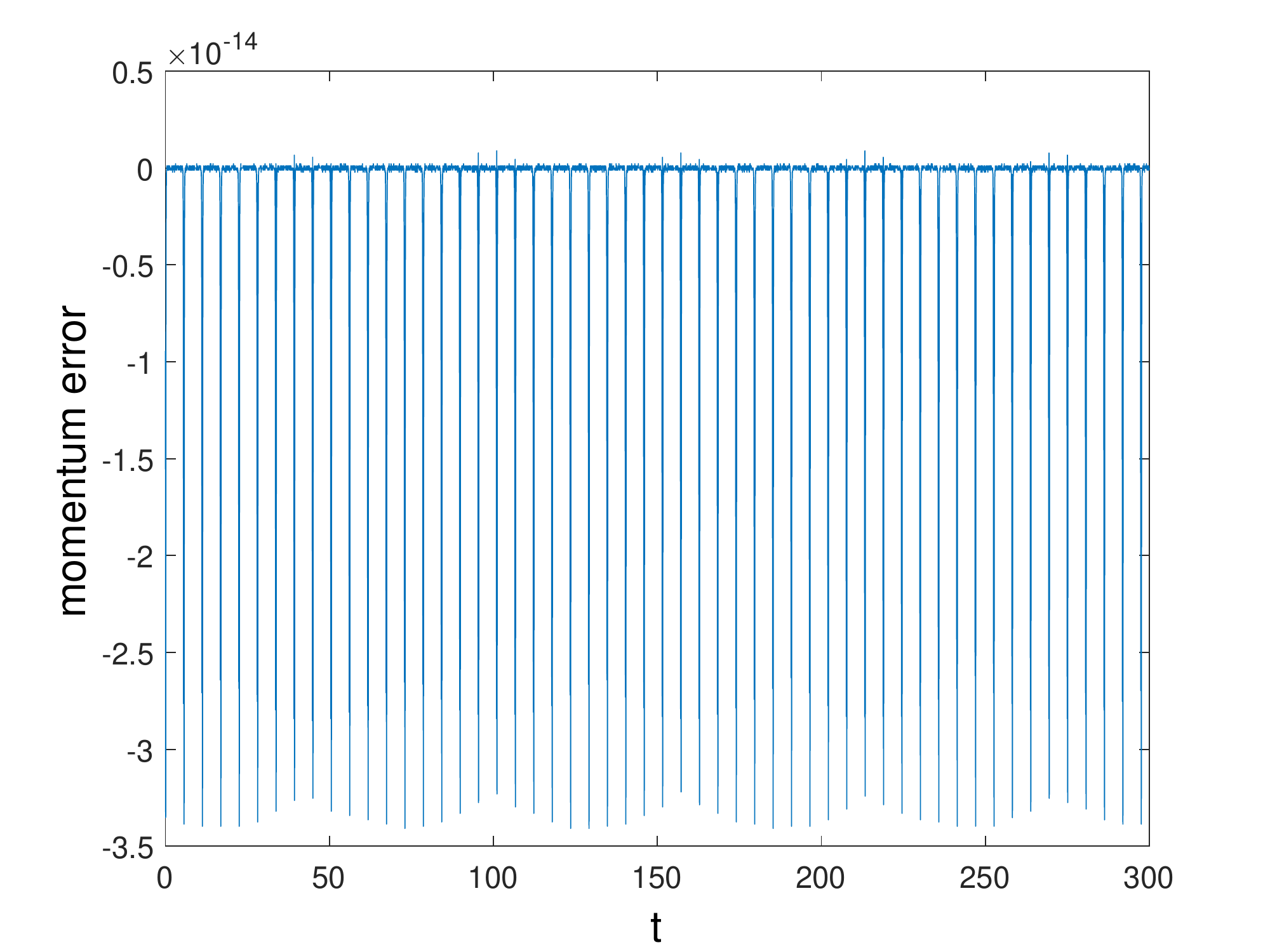}
\caption{The results in numerical orbits (left), energy error (middle) and momentum error (right) by methods \textbf{EIP-H} (upper), \textbf{EIP-L} (middle) and \textbf{EIP-HL} (bottom), respectively.}\label{ex1-fig1}
\end{figure}

Moreover, the numerical orbits as well as the errors in Hamiltonian energy and angular momentum are presented in Figures.~\ref{ex1-fig1} with time step $h=0.03$ for all the simulations. It is clear that all three methods can well perform the orbits like an ellipse that rotates slowly around one of its foci. For the conservation of invariants, the \textbf{EIP-H} and \textbf{EIP-L} methods exactly preserve the related single invariant and keep the errors in another one growing linearly under a small order of magnitude. Obviously, the \textbf{EIP-HL} method gives the best results with the two invariants being preserved to the machine accuracy, although it spends a little more CPU time than the former two EIP methods. A detailed comparison on the computational efficiency with respect to single and multiple invariants-preserving EIP methods will be carried out in the following PDE case.

\begin{figure}[H]
	\centering
	\includegraphics[width=0.3\linewidth]{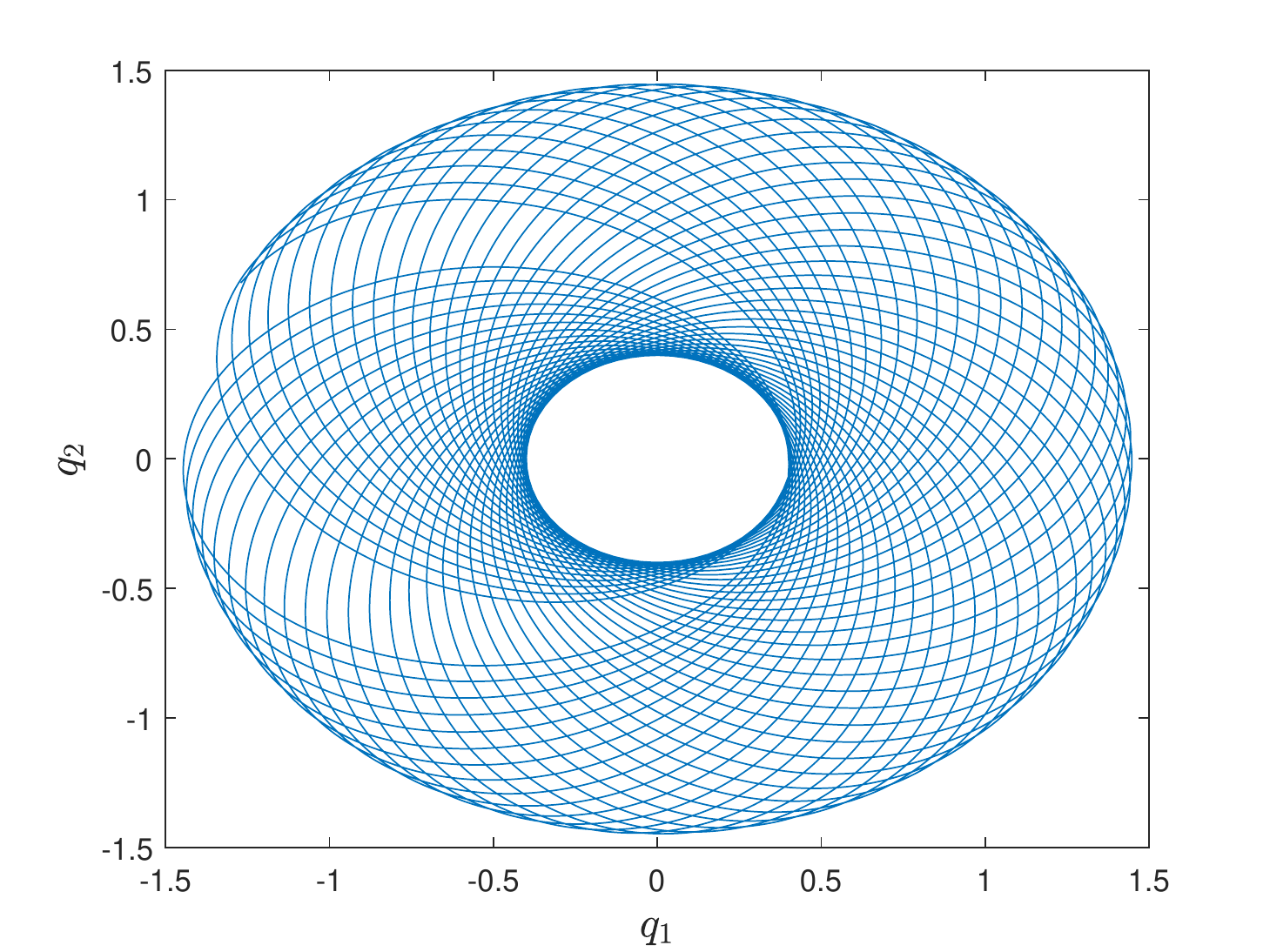}
	\includegraphics[width=0.3\linewidth]{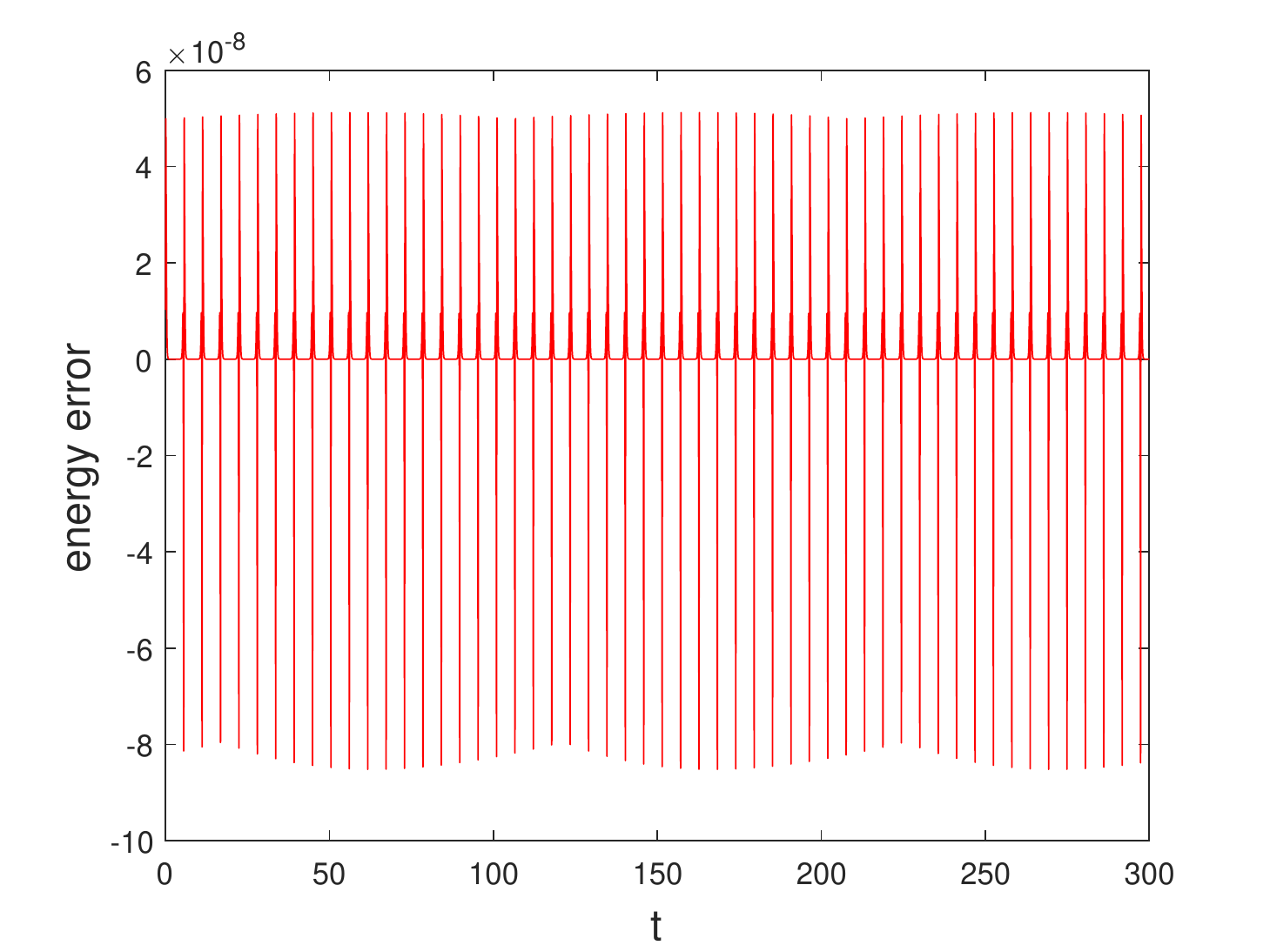}
	\includegraphics[width=0.3\linewidth]{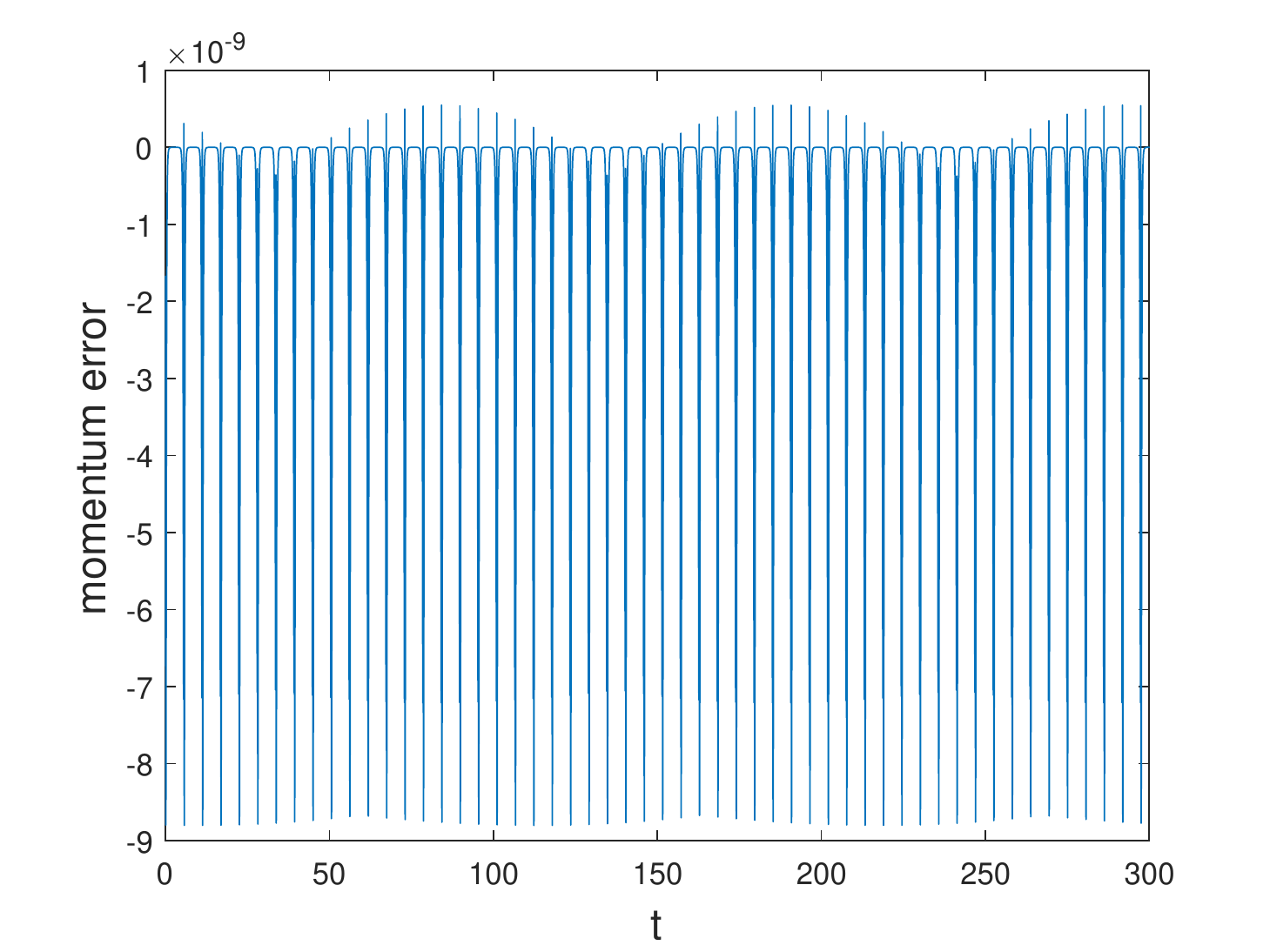}\\
	\includegraphics[width=0.3\linewidth]{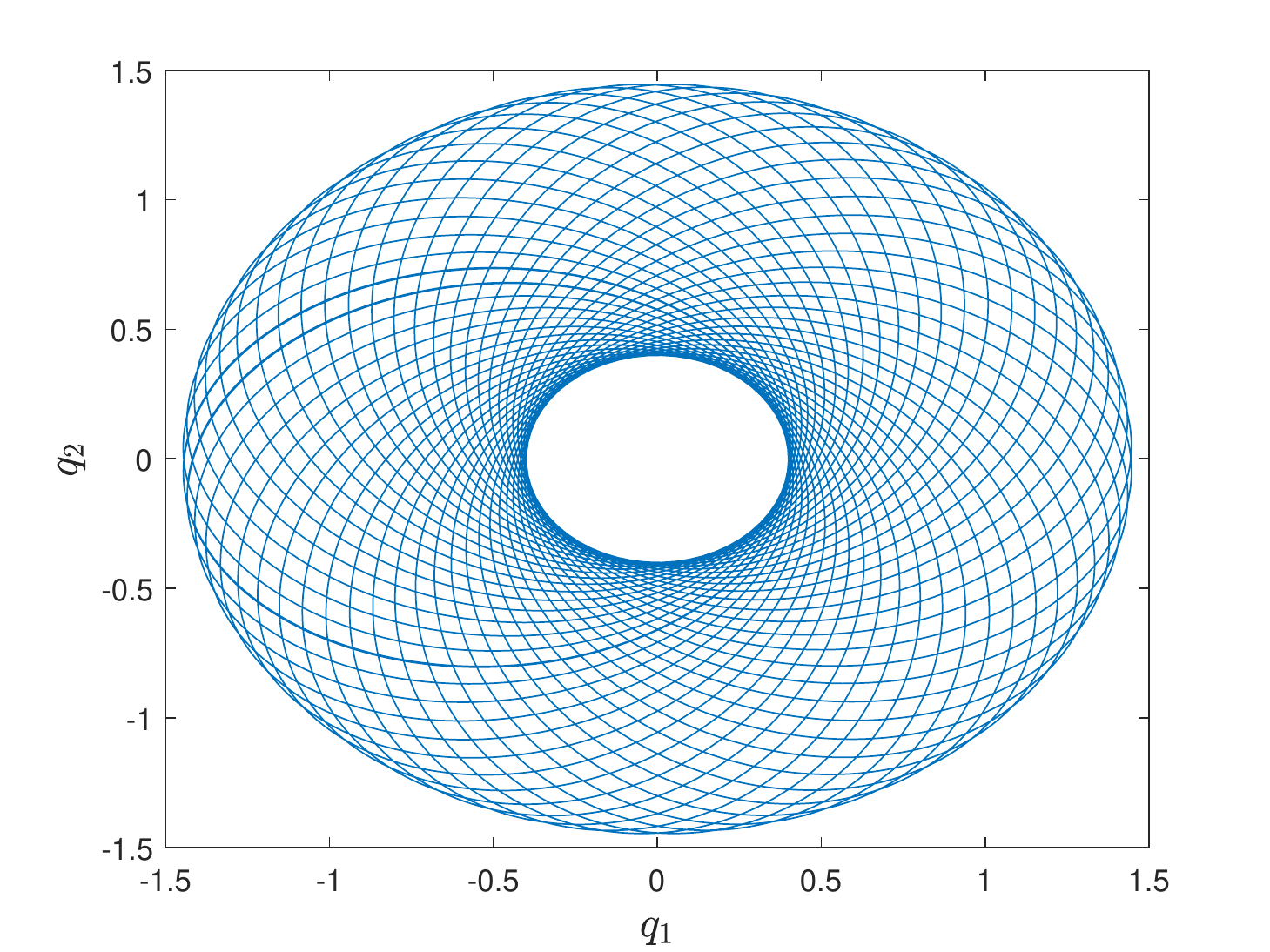}
	\includegraphics[width=0.3\linewidth]{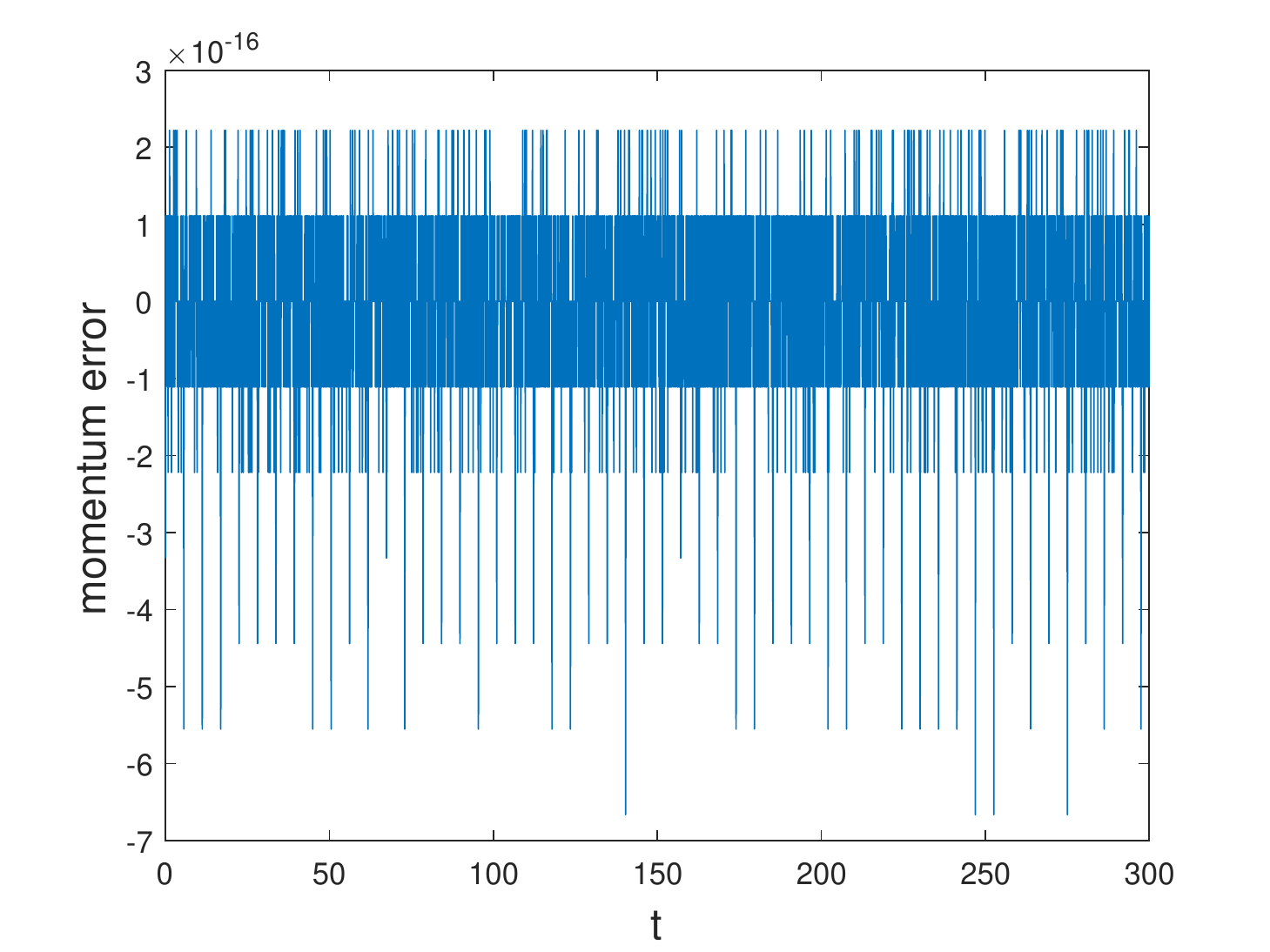}
	\includegraphics[width=0.3\linewidth]{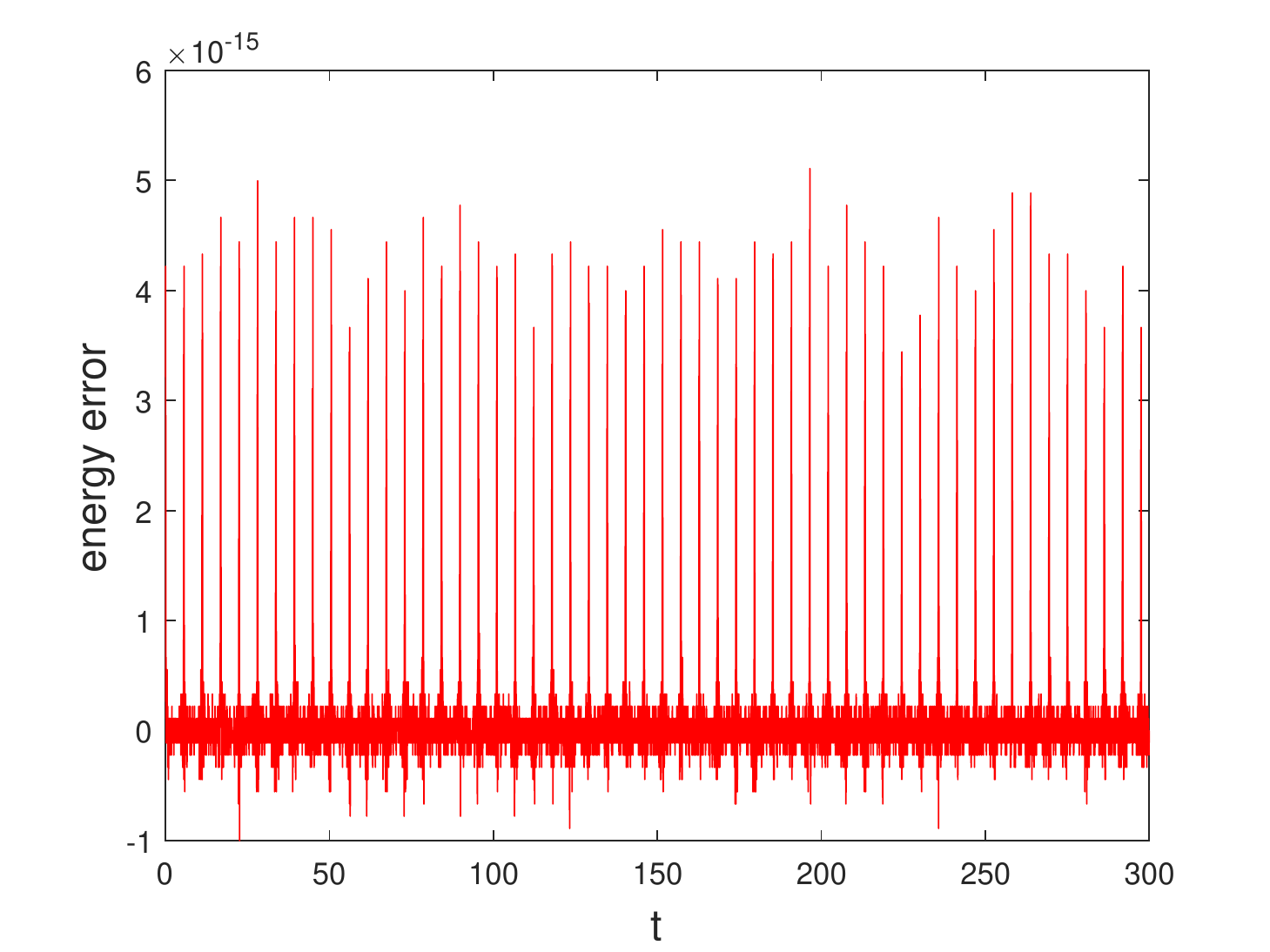}\\
	\caption{The results in numerical orbits (left), energy error (middle) and momentum error (right) by the method \textbf{EIP-HL} with RK2 (upper) and RK5 (bottom) as the underlying methods, respectively.}\label{ex1-fig2}
\end{figure}

As is mentioned that the EIP is a practical invariants-preserving method provided a high-order underlying RK method is utilized, e.g., a fourth-order one used above. For lower-order RK methods, the numerical behaviors may not be satisfied, whereas for higher-order ones the computational complexity will be significantly increased since for order $p>5$ the stages of explicit RK methods must be greater than $p$ \cite{hnw08} which makes it less cost-effective, especially for high-dimensional PDEs. To demonstrate this observation, we further present the results with the second- and fifth-order RK methods (RK2 and RK5 in short) as the underlying RK methods, where RK2 is listed above and RK5 is taken as the classic Fehlberg method with 6 stages \cite{hnw08}. From Figure.~\ref{ex1-fig2}, we can find that the orbit related to RK2 is not closed and the errors in energy and momentum are bounded but with magnitude much greater than the round-off error. While RK5 produces a correct orbit and the invariants errors exhibits even better than that in Figure.~\ref{ex1-fig1} by RK4 which coincides with the theoretical result. However, consider the requirement of practical computations, the magnitude of invariants errors by RK4 is already very satisfied and can be viewed as a conservation of invariants to round-off error. Therefore, from the perspective of cost performance we recommend RK4 as the underlying RK method of the EIP method. Following numerical tests will give a strong support of the choice.

\subsection{Experiment II:  solar system}

Next, we apply the EIP method to simulate a more realistic problem in celestial mechanics. Consider the motion of the solar system that describes eight planets and Pluto orbiting around the sun. As a generalization of the above Kepler problem, the corresponding Hamiltonian energy becomes
\begin{equation}\label{ex2-eq1}
H(p,q)=\frac{1}{2}\sum_{i=1}^{10}\frac{1}{m_i}p_i^\top p_i-\sum_{i=2}^{10}\sum_{i=1}^{i-1}\frac{Gm_im_j}{\|q_i-q_j\|},
\end{equation}
where $m_i$ represents the mass, and $p, q$ are supervectors composed by the momenta and position vectors $q_i, p_i\in\mathbb{R}^3$, respectively. Besides the energy, the solar system also admits the conservation of the angular momentum $L=\sum_{i=1}^{10}q_i\times p_i,$
which actually consists of three first integrals, namely
\begin{equation}\label{ex2-eq2}
\begin{aligned}
L_x=\sum_{i=1}^{10}\left(p_i(3)q_i(2)-p_i(2)q_i(3)\right),\\
L_y=\sum_{i=1}^{10}\left(p_i(1)q_i(3)-p_i(3)q_i(1)\right),\\
L_z=\sum_{i=1}^{10}\left(p_i(2)q_i(1)-p_i(1)q_i(2)\right),
\end{aligned}
\end{equation}
where $p_i(k), q_i(k)$ means the $k$-th component of momenta and position with respect to the $i$-th planet.

\begin{figure}[H]
	\centering
	\includegraphics[width=0.45\linewidth]{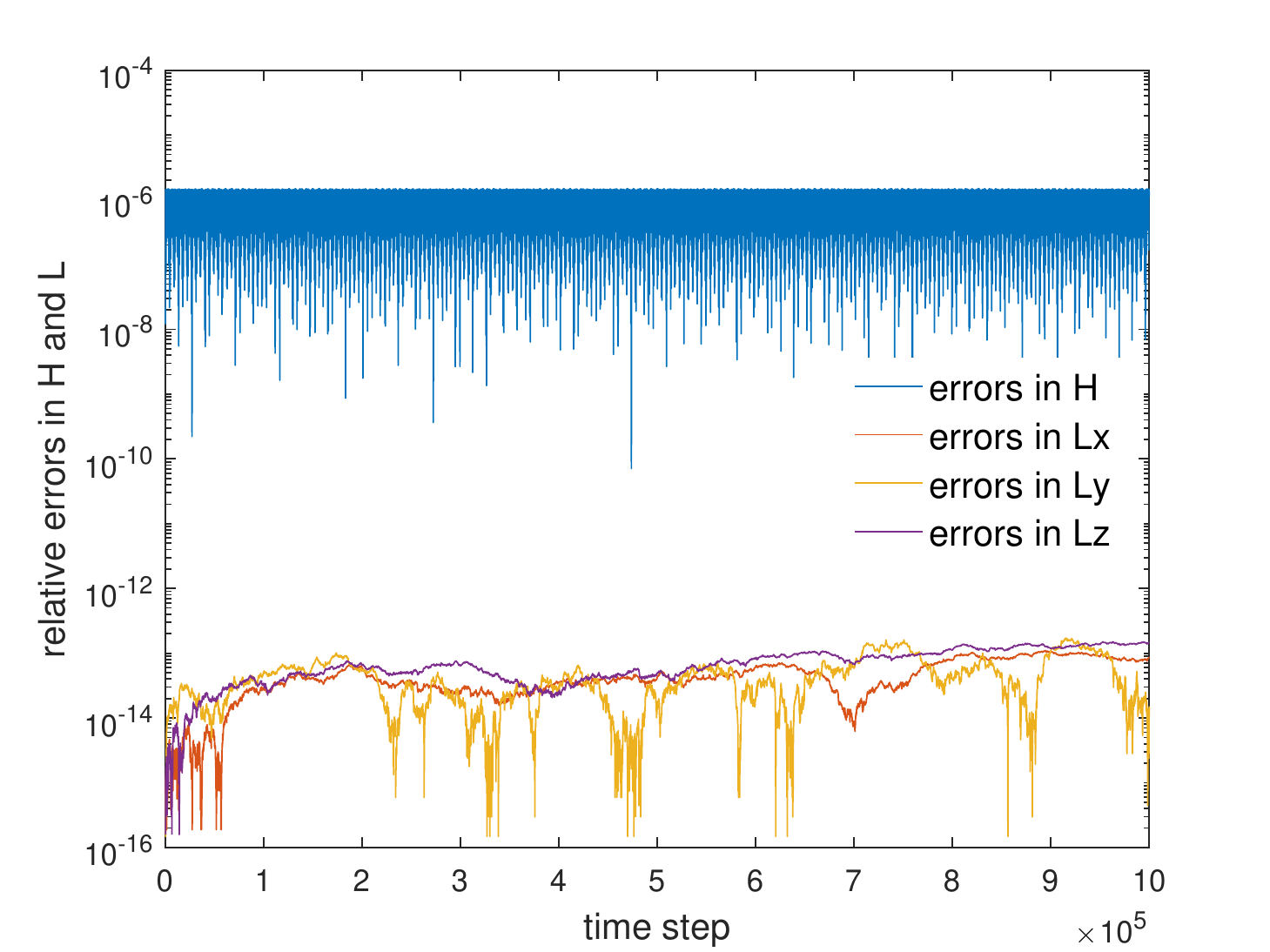}
	\includegraphics[width=0.45\linewidth]{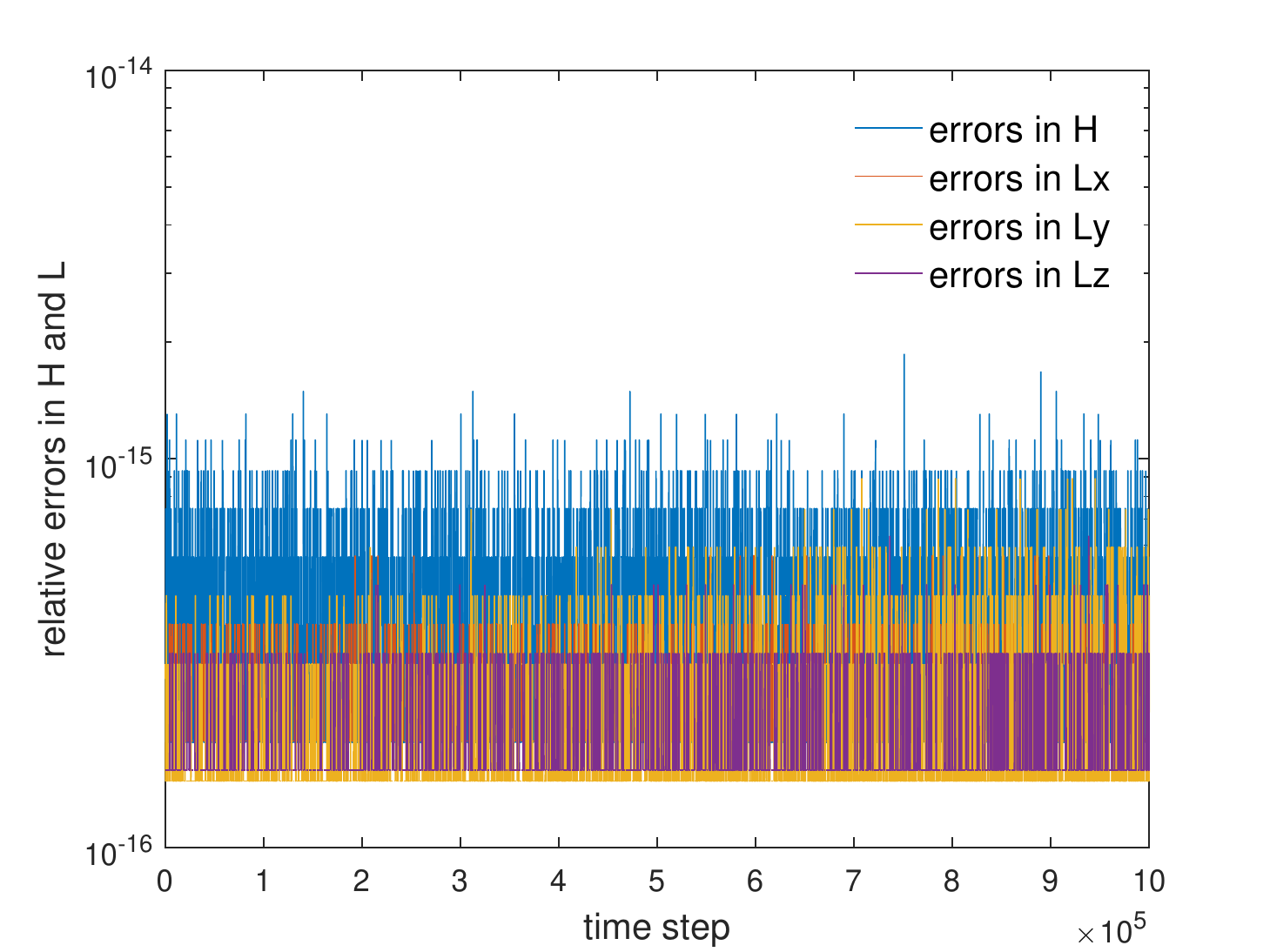}
	\caption{Relative invariant errors by \textbf{SV} (left) and the \textbf{EIP-HL} method (right).}\label{ex2-fig1}
\end{figure}

In the following experiments, to demonstrate the advantage of the proposed method in simultaneously preserving multiple invariants, we only focus on the numerical behaviors of the \textbf{EIP-HL} method. The initial datum of the solar system are taken from NASA JPL Ephemeris in the Appendix, and the simulation is carried out over 2000 years with $h=0.002$yr\footnote{Notice that the velocity in the Appendix is recorded in second and therefore the time step should also be transformed in second for practical computation.}. This long-term simulation does not only require the superior stability, but also the high precision of numerical algorithms. For comparison, we also present the numerical results by the popular second-order St\"ormer-Verlet method (\textbf{SV}) \cite{hlw06}. Since \textbf{SV} is a symplectic partitioned RK method, it can automatically preserve the three angular momentum to round-off errors and keep the energy errors bounded in a small order of amplitude in Figure.~\ref{ex2-fig1}. While the \textbf{EIP-HL} method can achieve the preservation of all the four invariants exactly and the errors behaves like random walk due to the machine accuracy.

\begin{figure}[H]
	\centering
			\includegraphics[width=0.3\linewidth]{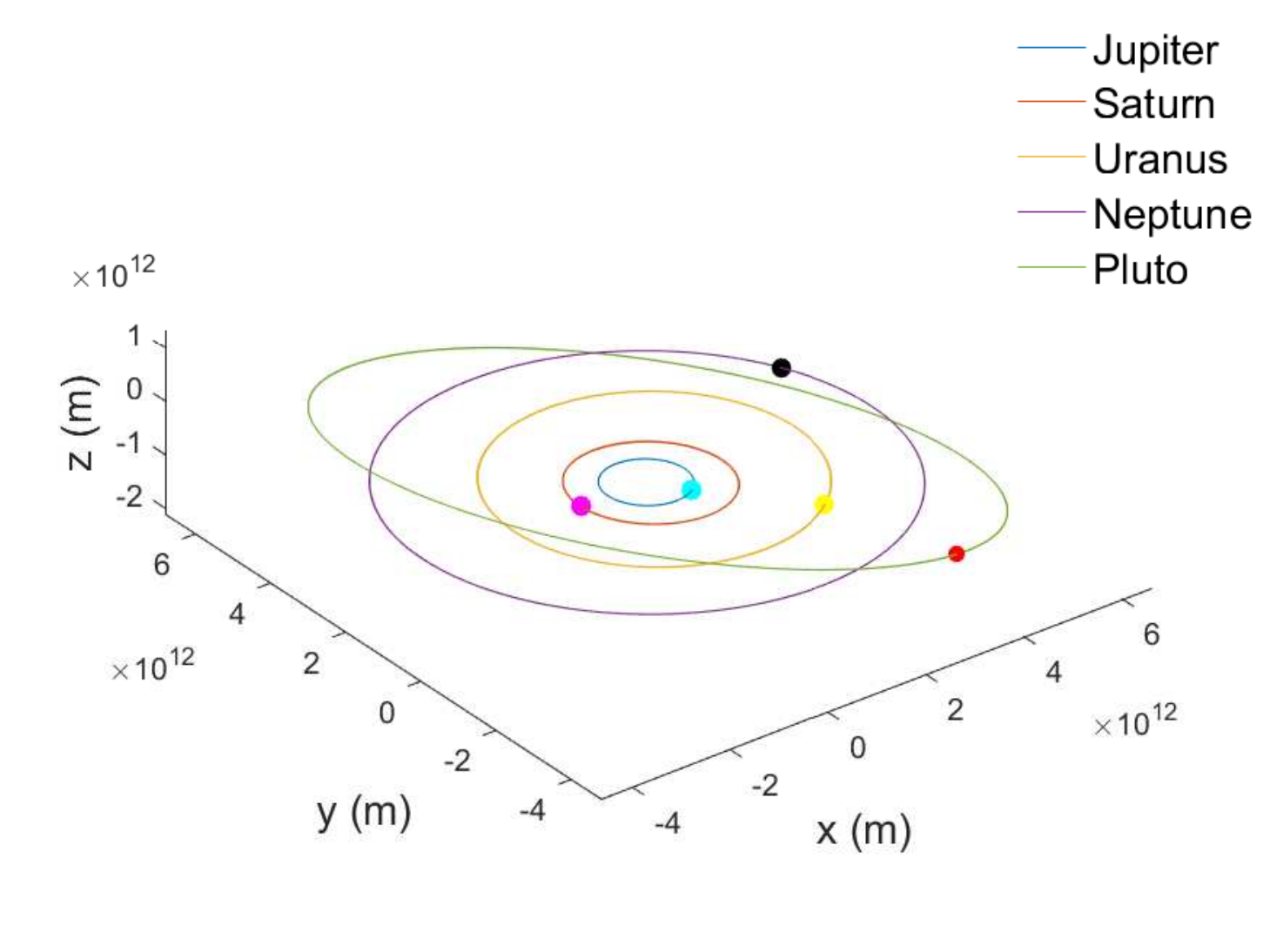}
			\includegraphics[width=0.3\linewidth]{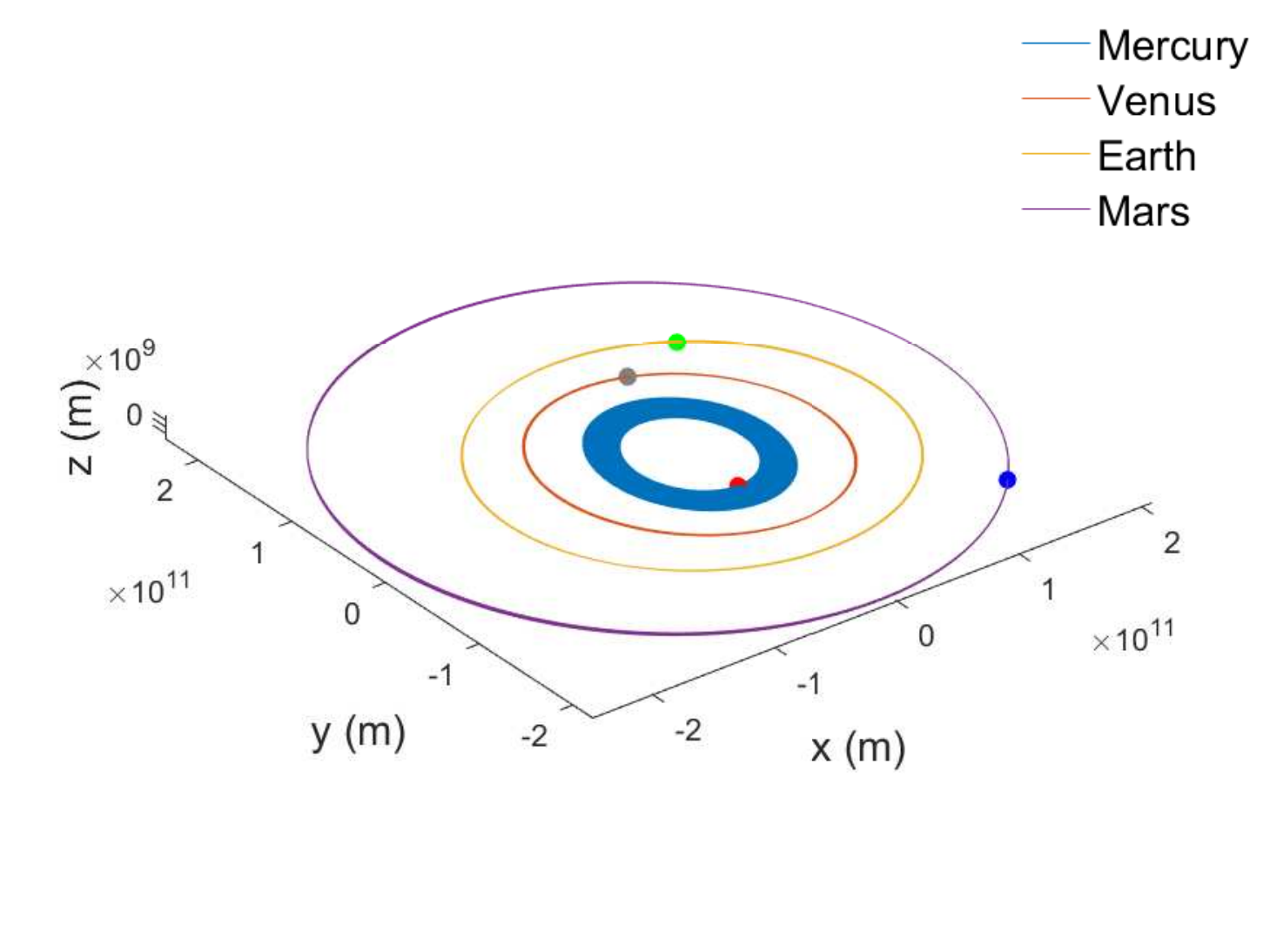}
			\includegraphics[width=0.3\linewidth]{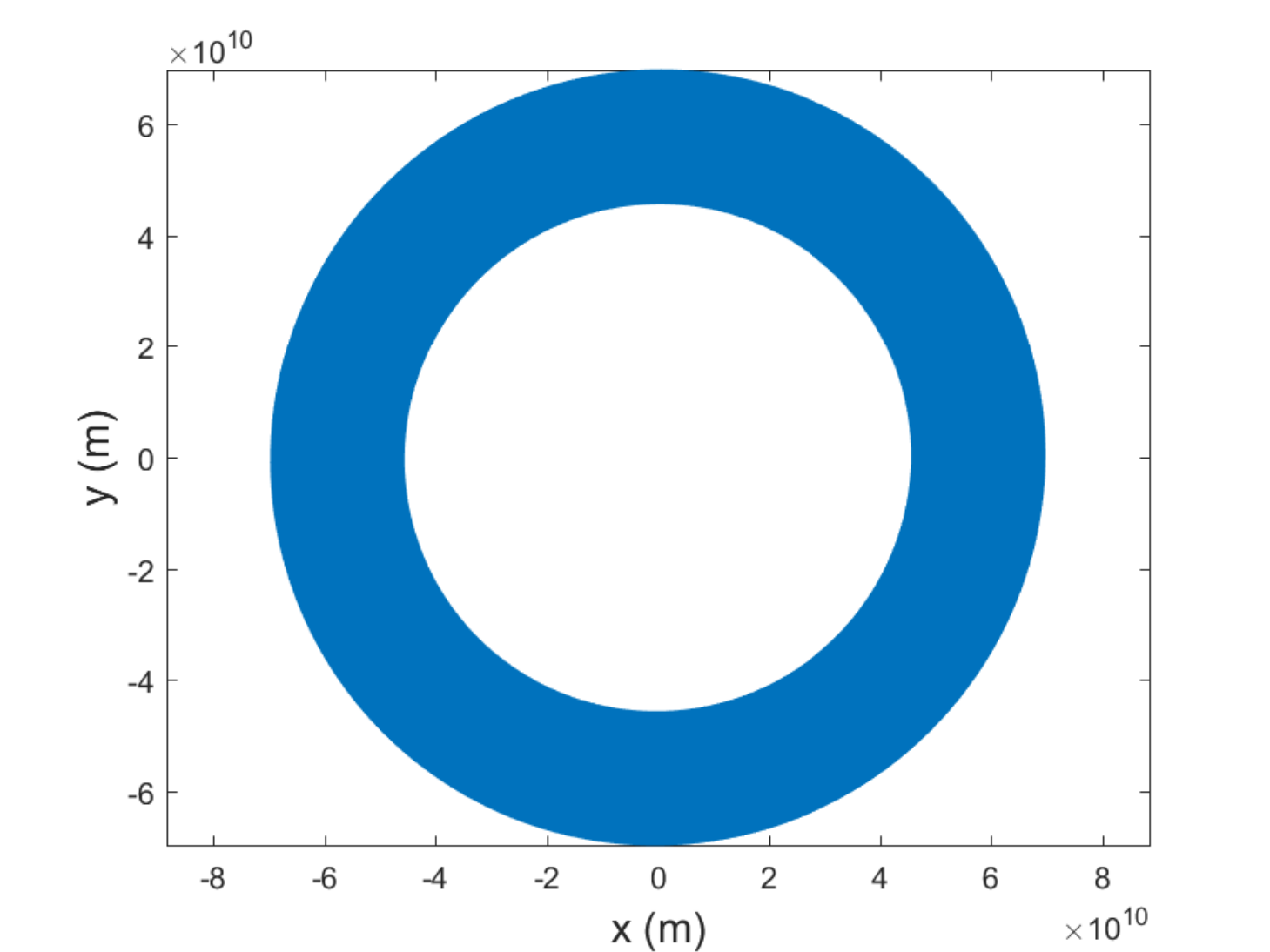}\\
			\includegraphics[width=0.3\linewidth]{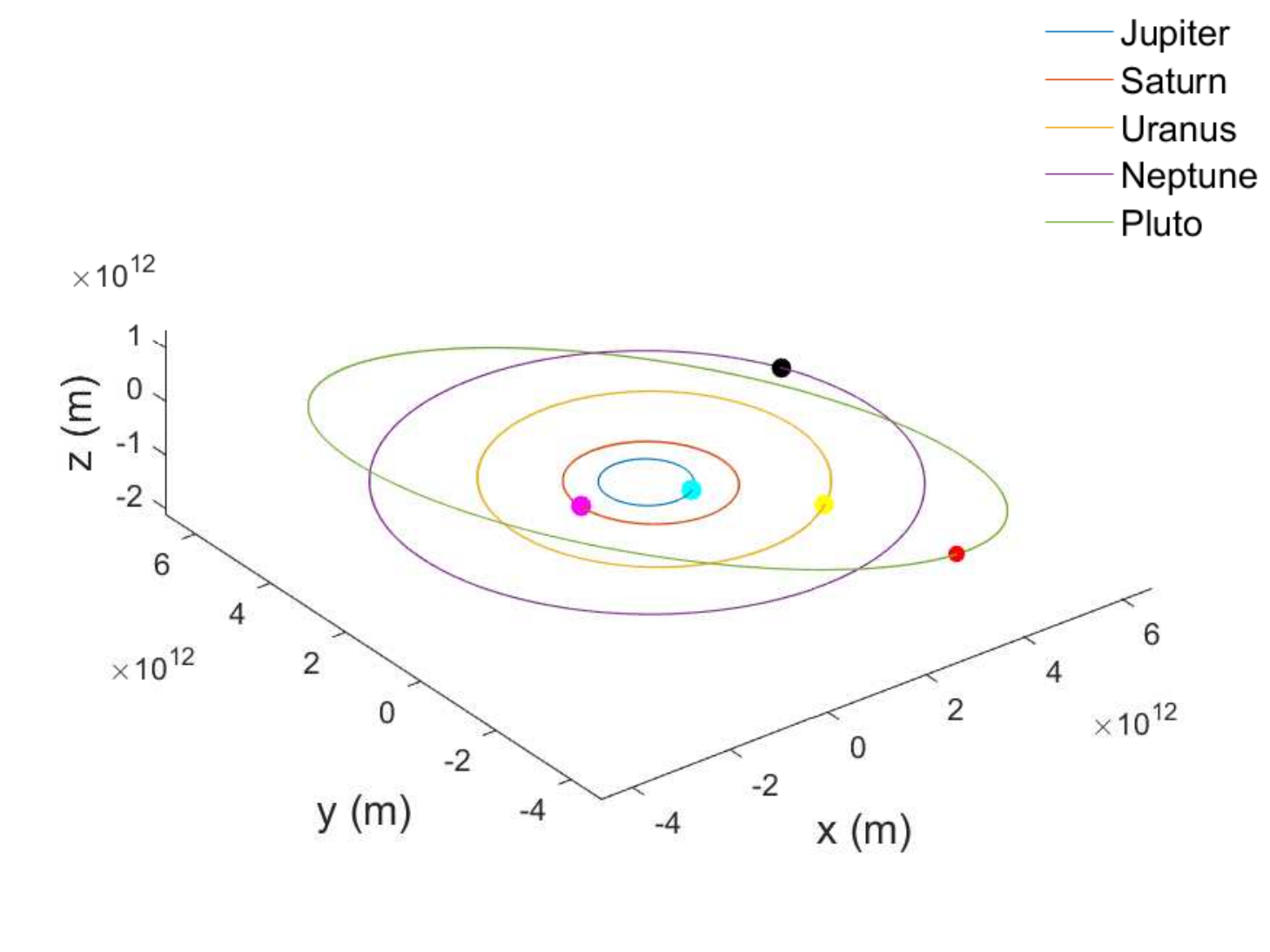}
			\includegraphics[width=0.3\linewidth]{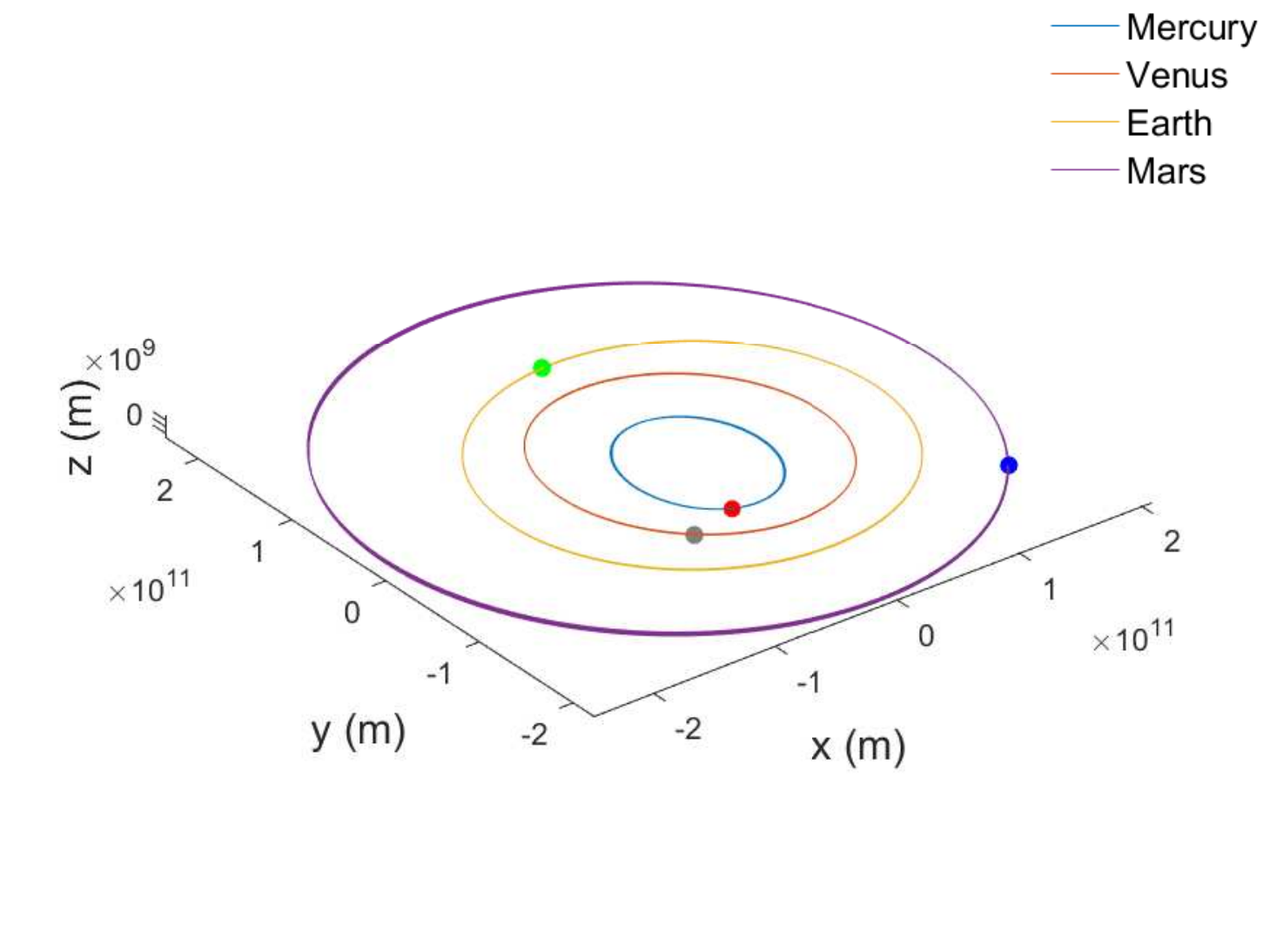}
			\includegraphics[width=0.3\linewidth]{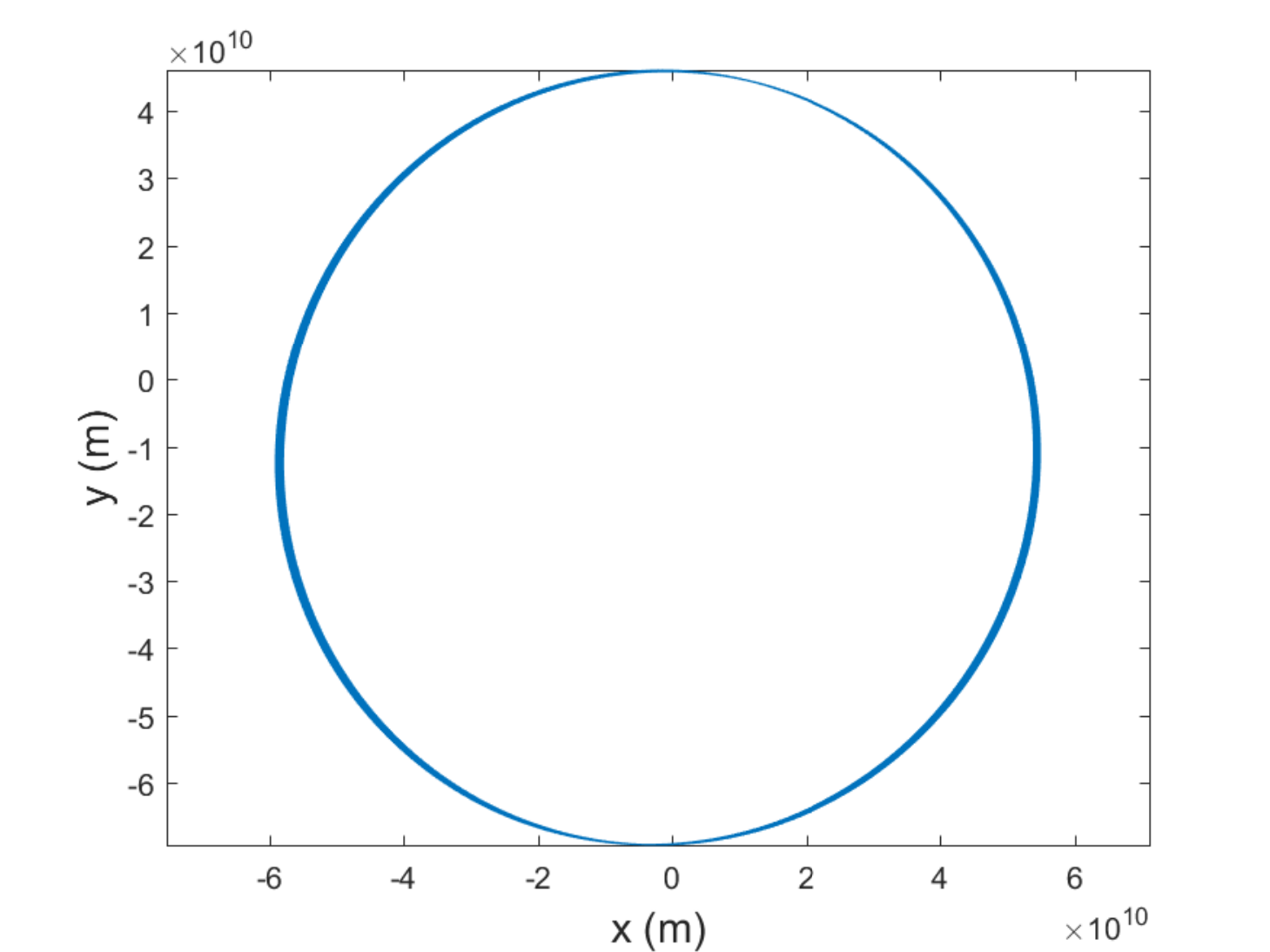}
	\caption{Orbits of outer planets and Pluto (left), orbits of inner four planets (middle) and the innermost orbit of Mercury  (right) by \textbf{SV} (upper) and the \textbf{EIP-HL} method (bottom), respectively.}\label{ex2-fig2}
\end{figure}

The high-accuracy advantage of the \textbf{EIP-HL} method for the solar system is revealed by the numerical orbits of planets in Figure.~\ref{ex2-fig2}, including the orbits of outer planets (and Pluto), the orbits of inner planets and the enlarged Mercury's orbit projected onto the $xy$-plane. We can hardly tell the difference from the outer orbits generated by \textbf{SV} and the \textbf{EIP-HL} method, no matter in the shapes of orbits or in the final positions of planets. However, it becomes clear from the inner orbits that the Mercury's orbit by \textbf{SV} is no long an ellipse. Instead, it suffers an undesirable precession effect, which can be observed more apparently from the projected picture on the $xy$-plane. The reason is mainly because of the low-order accuracy and the subsequent accumulation of phase errors, although its symplecticity guarantees the long-term stability. While for the high-order \textbf{EIP-HL} method both inner orbits and the enlarged Mercury's orbit exhibit correct elliptical shapes and the positions of planets are further improved.

\subsection{Experiment III: 2D charged particle dynamics}

As the underlying one step method, RK4 has the same order of accuracy as the presented EIP methods. Nevertheless, the numerical behaviors can be distinct. For illustration, we consider an example of 2D charged particle dynamics.

In a given electromagnetic field $\left(\bm E(\bm x), \bm B(\bm x)\right)$, the motion of a
charged particle is governed by the Lorentz force law and can be
described by
\begin{equation}\label{ex3-1}
m\bm{\ddot{x}}=q(\bm E+\bm{\dot{x}}\times \bm B),\quad \bm x\in\mathbb{R}^3
\end{equation}
where $\bm x$ is the position of the charged particle, $m$ is the mass, and $q$ denotes the electric charge. For convenience, we assume that $\bm B, \bm E$ are static and thus $\bm B=\nabla\times\bm A$ and $\bm E=-\nabla\varphi$ with $\bm A$ and $\varphi$ the potentials. Let the conjugate momentum $\bm p=m\dot{\bm x}+q\bm A(\bm x)$, the system \eqref{ex3-1} has a canonical symplectic structure \cite{he-16-VP-CPD} with the Hamiltonian given by
\begin{equation}\label{ex3-2}
H(\bm x, \bm p)=\frac{1}{2m}\left(\bm p-q\bm A(\bm x)\right)\cdot \left(\bm p-q\bm A(\bm x)\right)+q\varphi(\bm x),
\end{equation}
which however cannot be split into the form $H(\bm x,\bm p)=T(\bm p)+V(\bm x)$.
 A separable formulation of \eqref{ex3-1} can be found by recasting it with transformation $(\bm x, \bm p)\rightarrow(\bm x, \bm v): \bm x=\bm x, \bm v=\bm p/m-q\bm A(\bm x)/m$, as
\begin{equation}\label{ex3-3}
\begin{aligned}
&\dot{\bm x}=\bm v,\\
&\dot{\bm{v}}=\frac{q}{m}\left(\bm E(\bm x)+\bm v\times\bm B(\bm x)\right),
\end{aligned}
\end{equation}
where the corresponding energy becomes $H(\bm x, \bm v)=m\bm v\cdot\bm v/2+q\varphi(\bm x)$. Though several structure-preserving methods have been proposed for numerical solving \eqref{ex3-1}, including symplectic methods \cite{zhang-14-S-CPD,tao-16-explicit-S}, volume-preserving methods \cite{qin-13-Boris, he-16-VP-CPD,gong-18-Gauss-Seidel} and energy-conserving methods \cite{li-19-EP-CPD,brugnano-19-LIM-CPD}, with respect to the structures \eqref{ex3-2} or \eqref{ex3-3}. Among those methods, explicit ones are only valid for the separable formulation \eqref{ex3-3}.
To the best of our knowledge no explicit structure-preserving method, especially no energy-conserving method, exists based on the canonical form \eqref{ex3-2}. Hence, in the following experiments, we will utilize the proposed EIP method  to construct a first explicit energy-conserving method based on this form. One can follow a similar process to derive the EIP method according to the separable system \eqref{ex3-3} and we omit it here.

\begin{figure}[H]
	\centering
	\includegraphics[width=0.3\linewidth]{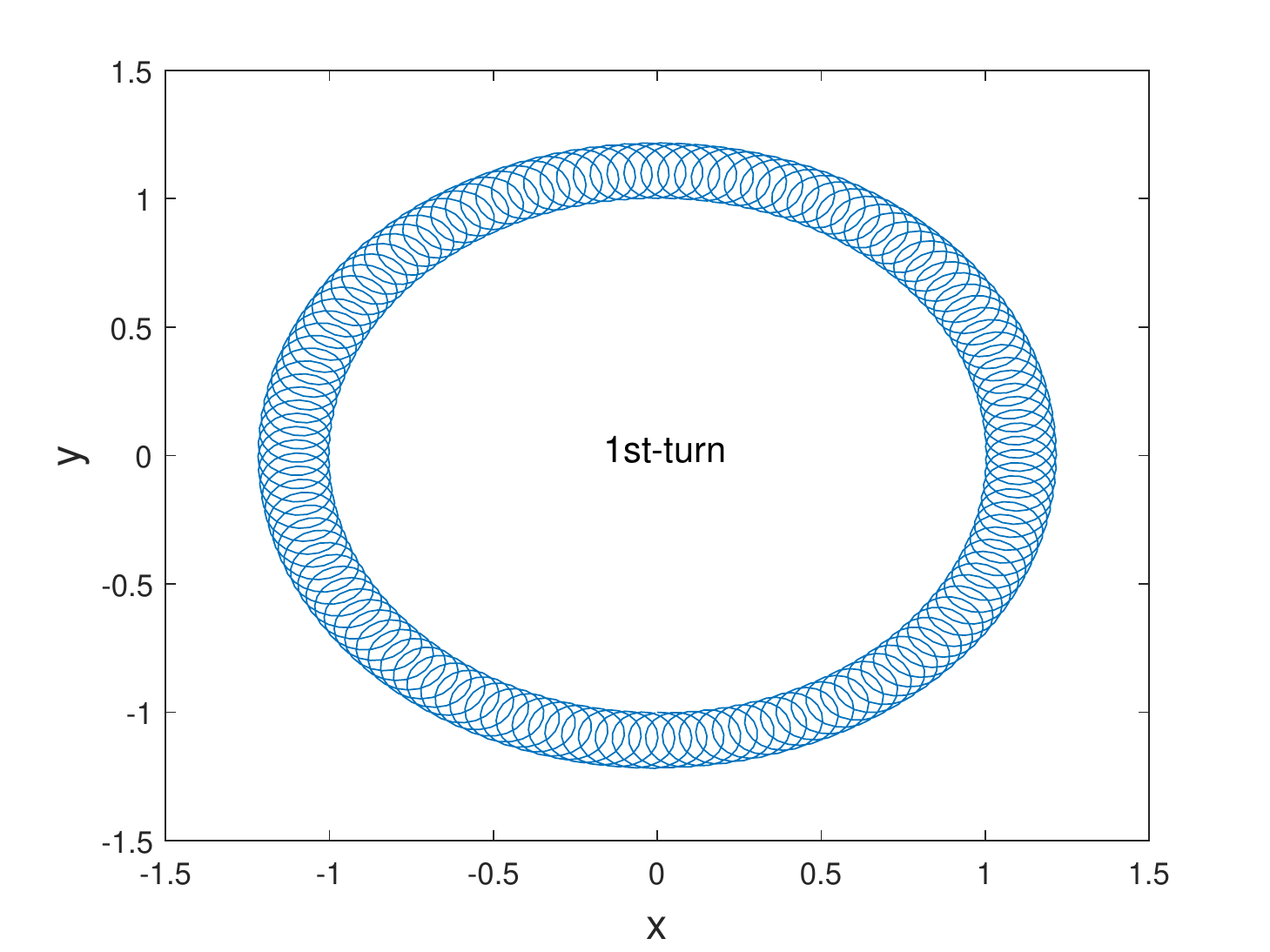}
	\includegraphics[width=0.3\linewidth]{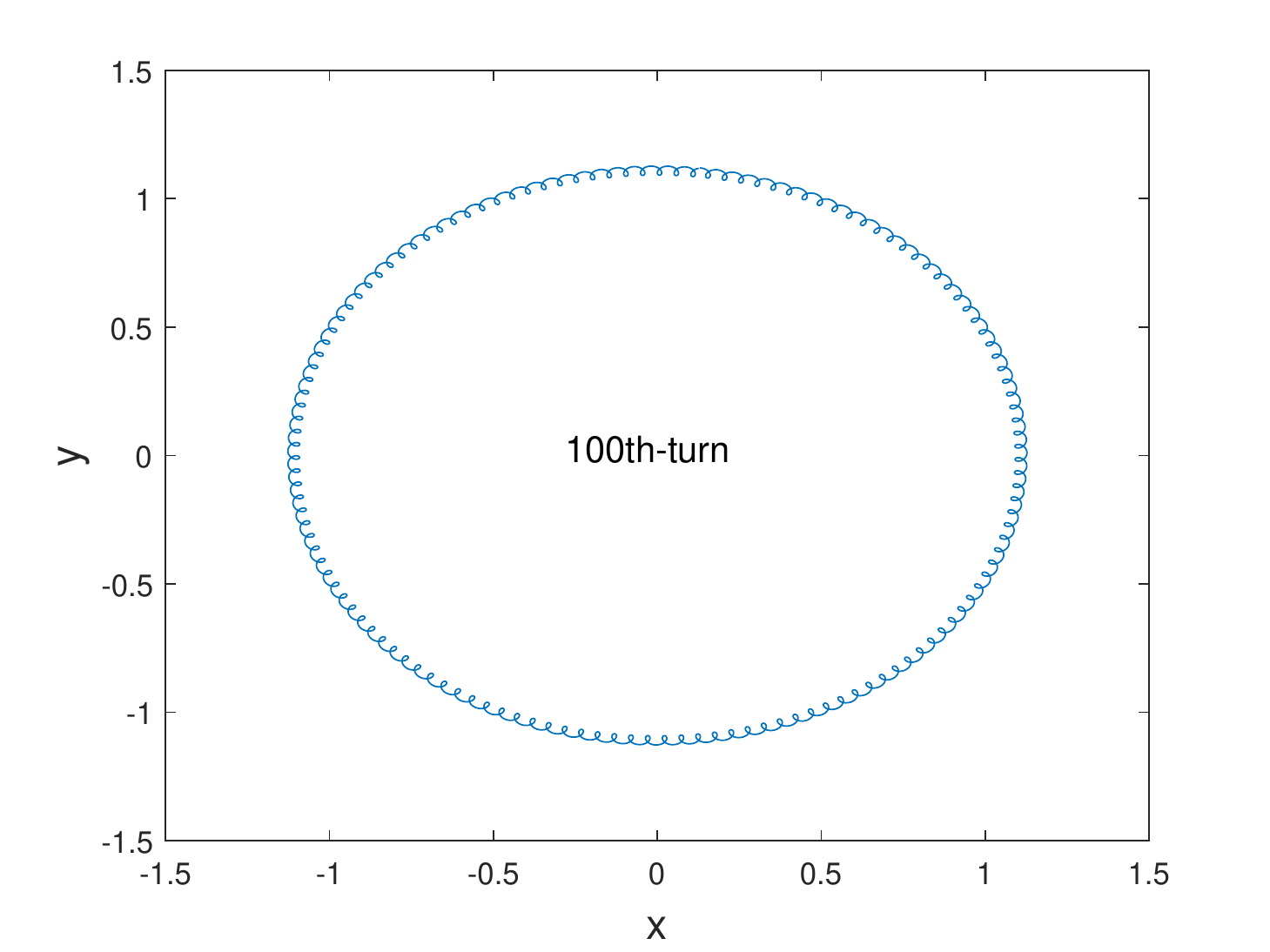}
	\includegraphics[width=0.3\linewidth]{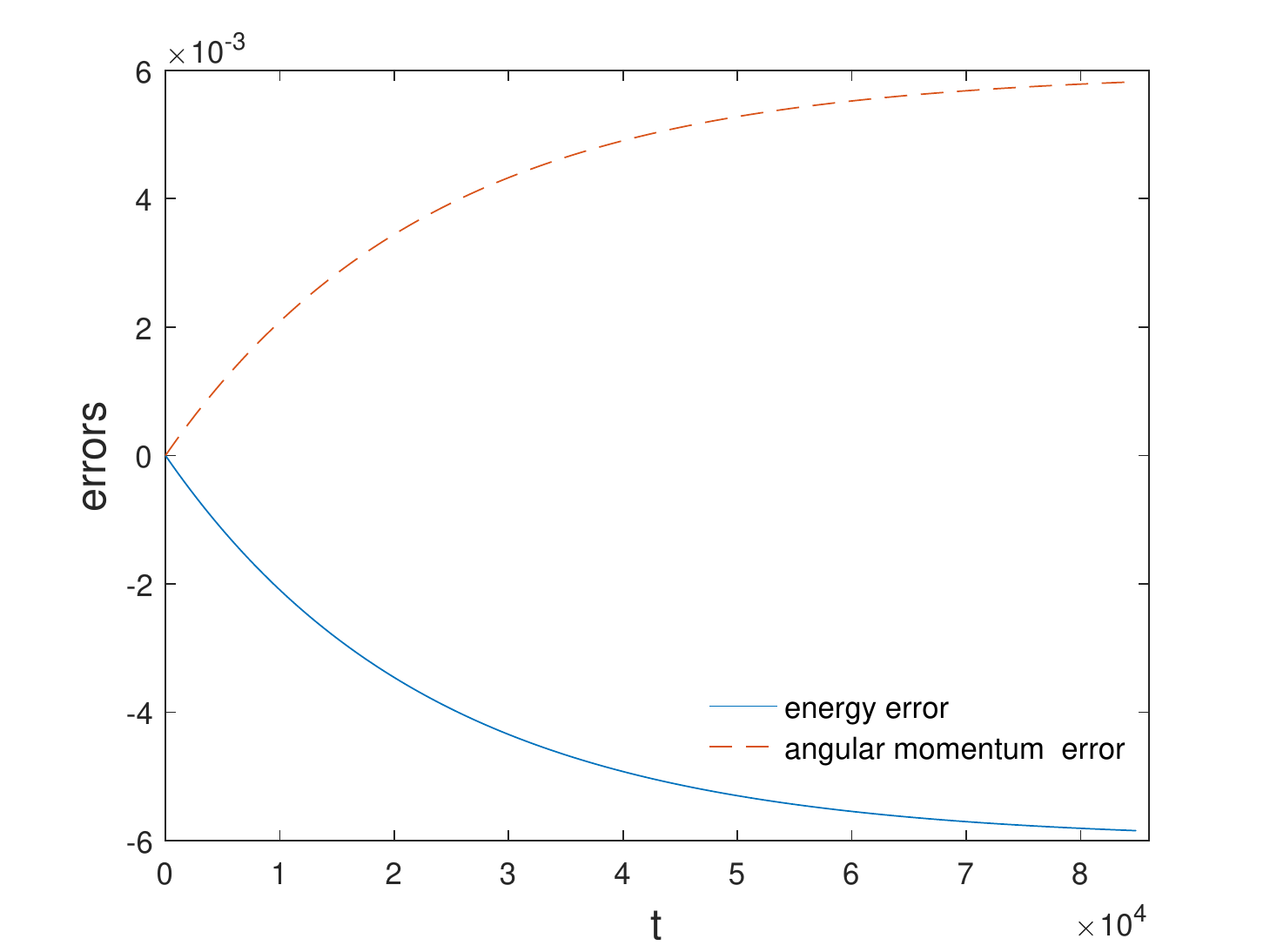}\\
	\includegraphics[width=0.3\linewidth]{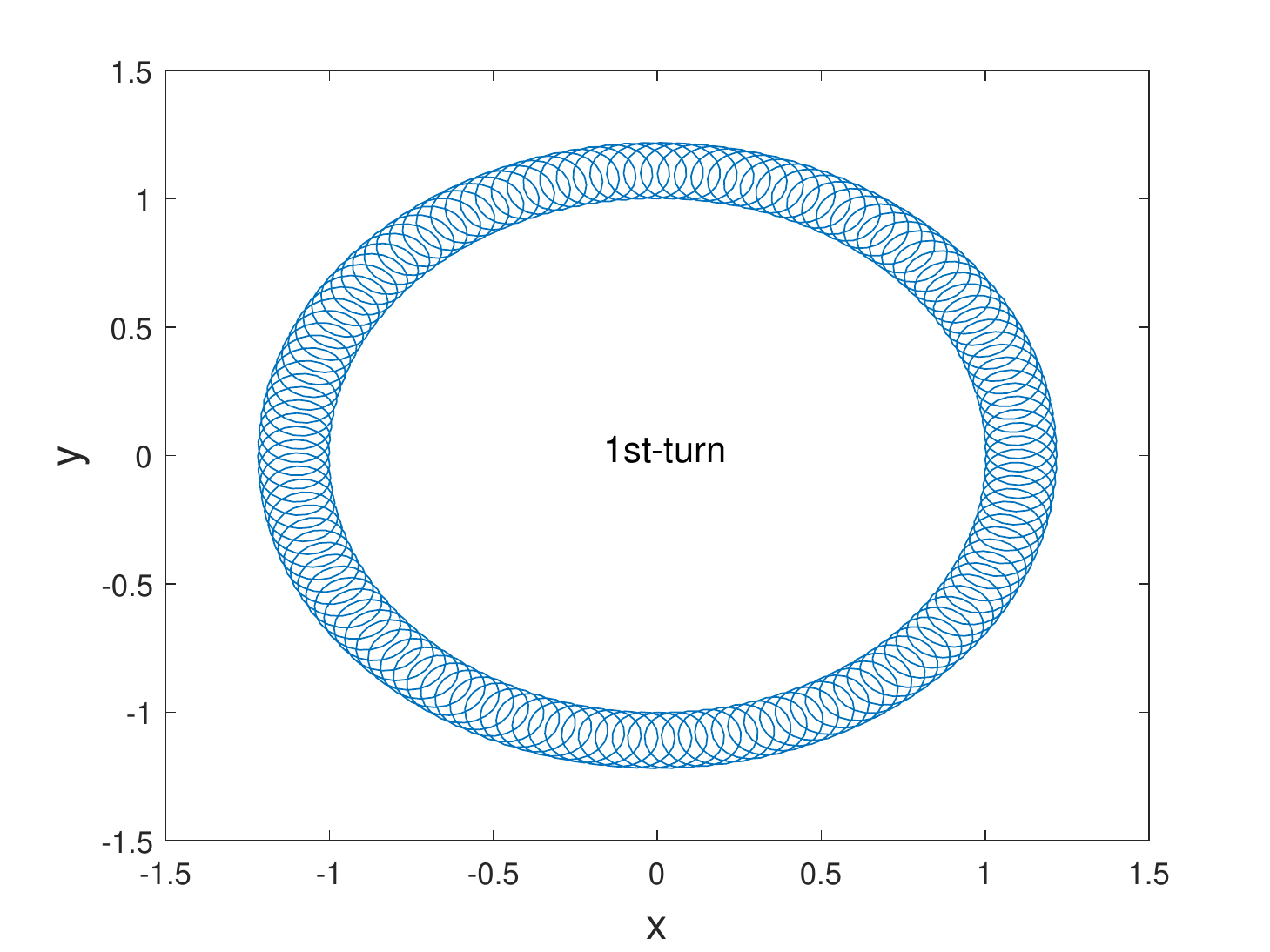}
	\includegraphics[width=0.3\linewidth]{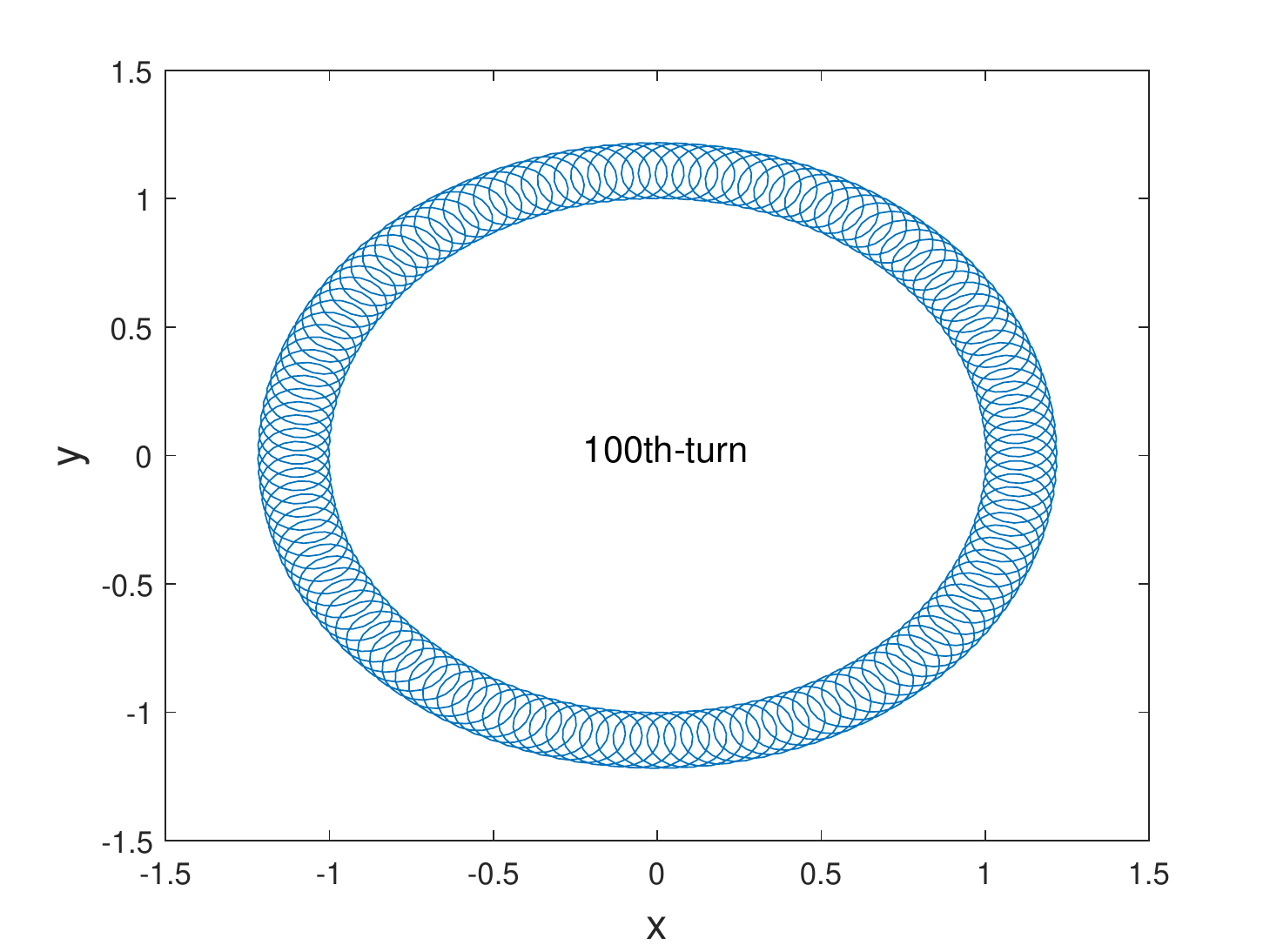}
	\includegraphics[width=0.3\linewidth]{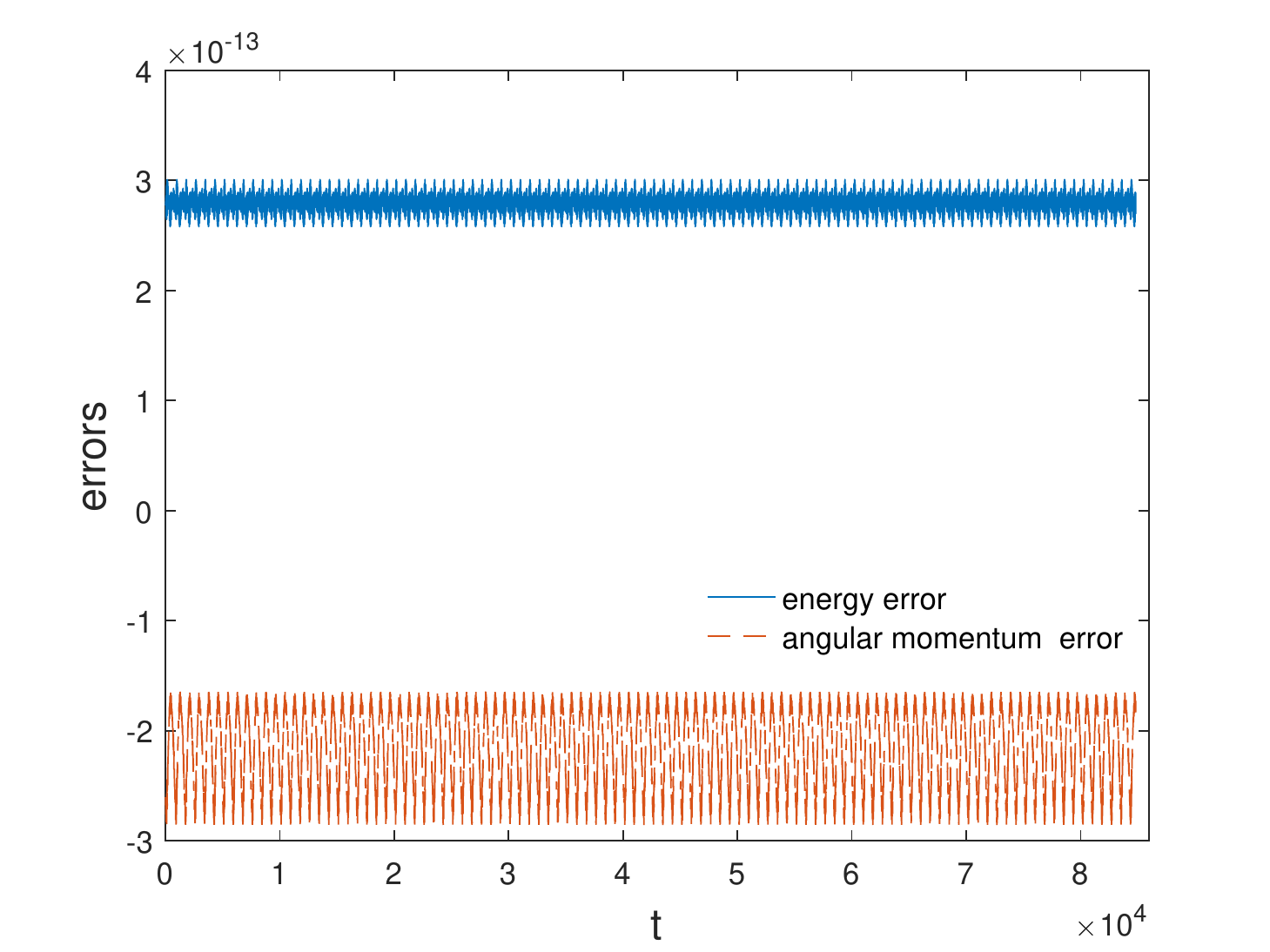}
	\caption{Numerical orbits and errors in two invariants solved on the interval $[0, 2.7\times 10^5h]$ by RK4 (upper) and the \textbf{EIP-HL} method (bottom).}\label{ex3-fig1}
\end{figure}

First, we consider the 2D dynamics of the charged particle in a static, non-uniform electromagnetic field
\[
\bm B=\nabla\times\bm A=\bm e_z,\quad \bm E=-\nabla\varphi=\frac{10^{-2}}{R^3}(x\bm e_x+y\bm e_y),
\]
where the potentials are chosen as $\bm A=[-y/2, x/2, 0]^\top, \varphi=\frac{10^{-2}}{R}, R=\sqrt{x^2+y^2}.$ In this example, the physical quantities are normalized by the system size $a$, the characteristic magnetic field $B_0$, and the gyro-frequency $\omega_0\equiv qB_0/m$ of the particle. Besides the energy \eqref{ex3-2}, another invariant of this case is given by the angular momentum
\begin{equation}\label{ex3-4}
L(\bm x, \bm p)=xp_y-xp_x.
\end{equation}
Therefore, we construct the method \textbf{EIP-HL} that preserves both energy and angular momentum. Starting from the initial conditions $\bm x_0=[0,-1,0]^\top$, $\bm v_0=[0.1,0.01,0]^\top$, and taking the step size $h=\pi/10$, which is the $1/20$ of the characteristic gyro-period, we run the \textbf{EIP-HL} method for $2.7\times 10^5$ steps. The particle's exact orbit is a spiraling circle with a constant radius. For comparison, we also apply RK4 for this problem. Figure.~\ref{ex3-fig1} presents the orbits generated by these two methods and the corresponding errors in the invariants. Though RK4 has a same fourth-order accuracy,
the numerical error accumulation after $2.7\times 10^5$ steps gives rise to a complete wrong solution orbit during the 100th-turn, where gyro-motion is numerically dissipated. This can be also confirmed from the errors in two invariants. While the \textbf{EIP-HL} method can not only provide a correct and stable gyro-motion over such a long-term simulation but also preserve the energy and angular momentum to round-off errors.

Next, we test our method for the 2D dynamics in an axisymmetric tokamak geometry without inductive electric field where
\[
\bm B=-\frac{2y+xz}{2R^2}\bm e_x+\frac{2x-yz}{2R^2}\bm e_y+\frac{R-1}{2R}\bm e_z,
\]
and
\[
\bm A=\left[\frac{xz}{2R^2}-\frac{\left((1-R)^2+z^2\right)y}{4R^2}\right]\bm e_x+\left[\frac{yz}{2R^2}+\frac{\left((1-R)^2+z^2\right)x}{4R^2}\right]\bm e_y-\frac{1}{2}\log(R)\bm e_z.
\]
The initial conditions are taken as $\bm x_0=[1.05,0,0]^\top, \bm v_0=[0,4.816\mbox{e-4}, -2.059\mbox{e-3}]^\top$. The exact orbit projected on $(R,z)$ space is a banana orbit, and it will turn into a transit orbit when the initial velocity is changed to $\bm v_0=[0,2\times4.816\mbox{e-4}, -2.059\mbox{e-3}]^\top$. The time step is also set to $h=\pi/10$. The results produced by RK4 and the \textbf{EIP-HL} method are given in Figure.~\ref{ex3-fig1}. It can be observed that the banana orbit by RK4 gradually transformed into a circulating orbit and the transit orbit deviates to the right side which is mainly due to its numerical dissipation. Again, the \textbf{EIP-HL} method can provide correct orbits as well as an exact conservation of the invariants.

From the above two tests, we can conclude that the EIP method essentially improve the numerical performance of the standard fourth-order RK method.

\begin{figure}[H]
	\centering
	\includegraphics[width=0.3\linewidth]{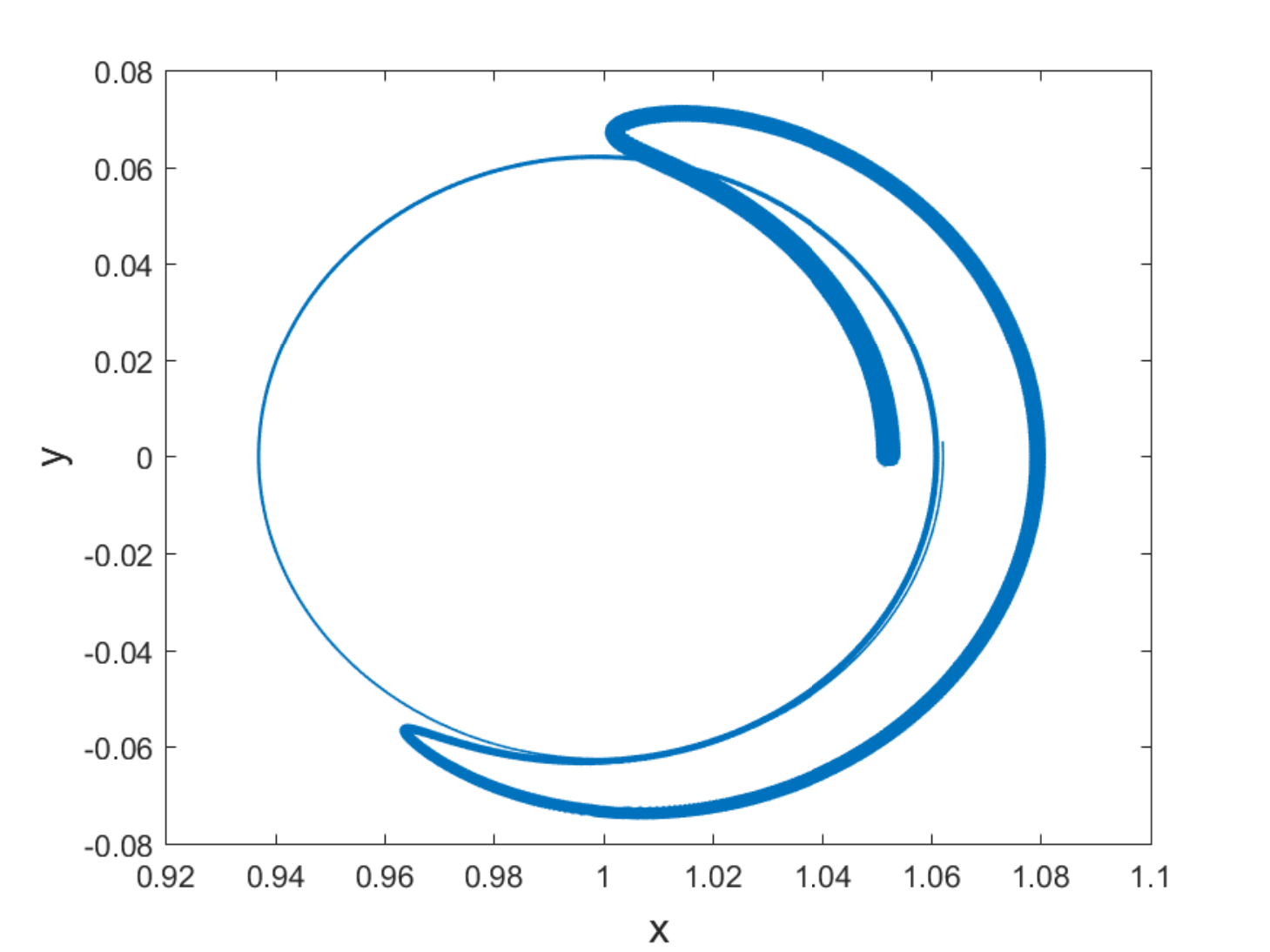}
	\includegraphics[width=0.3\linewidth]{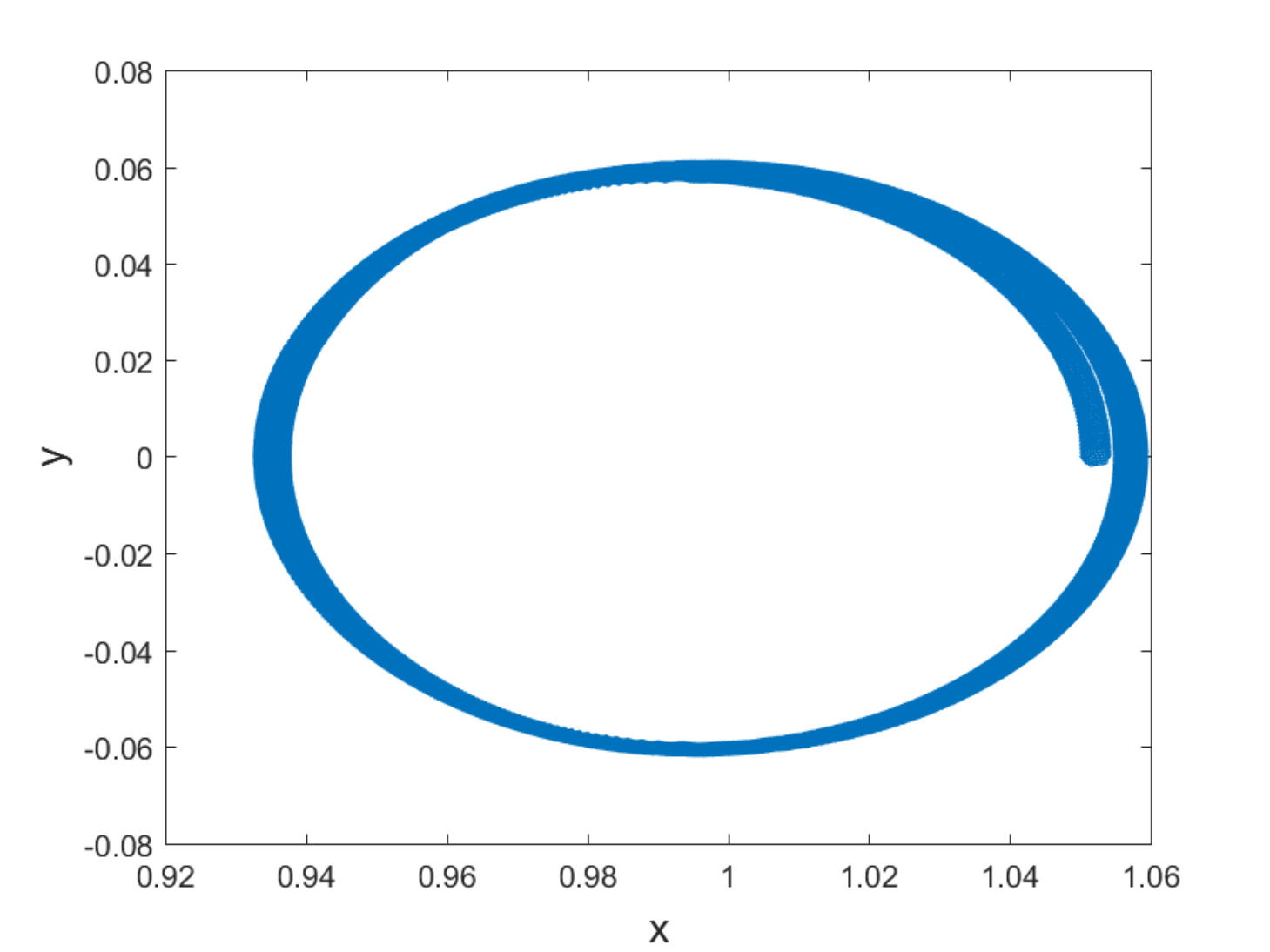}
	\includegraphics[width=0.3\linewidth]{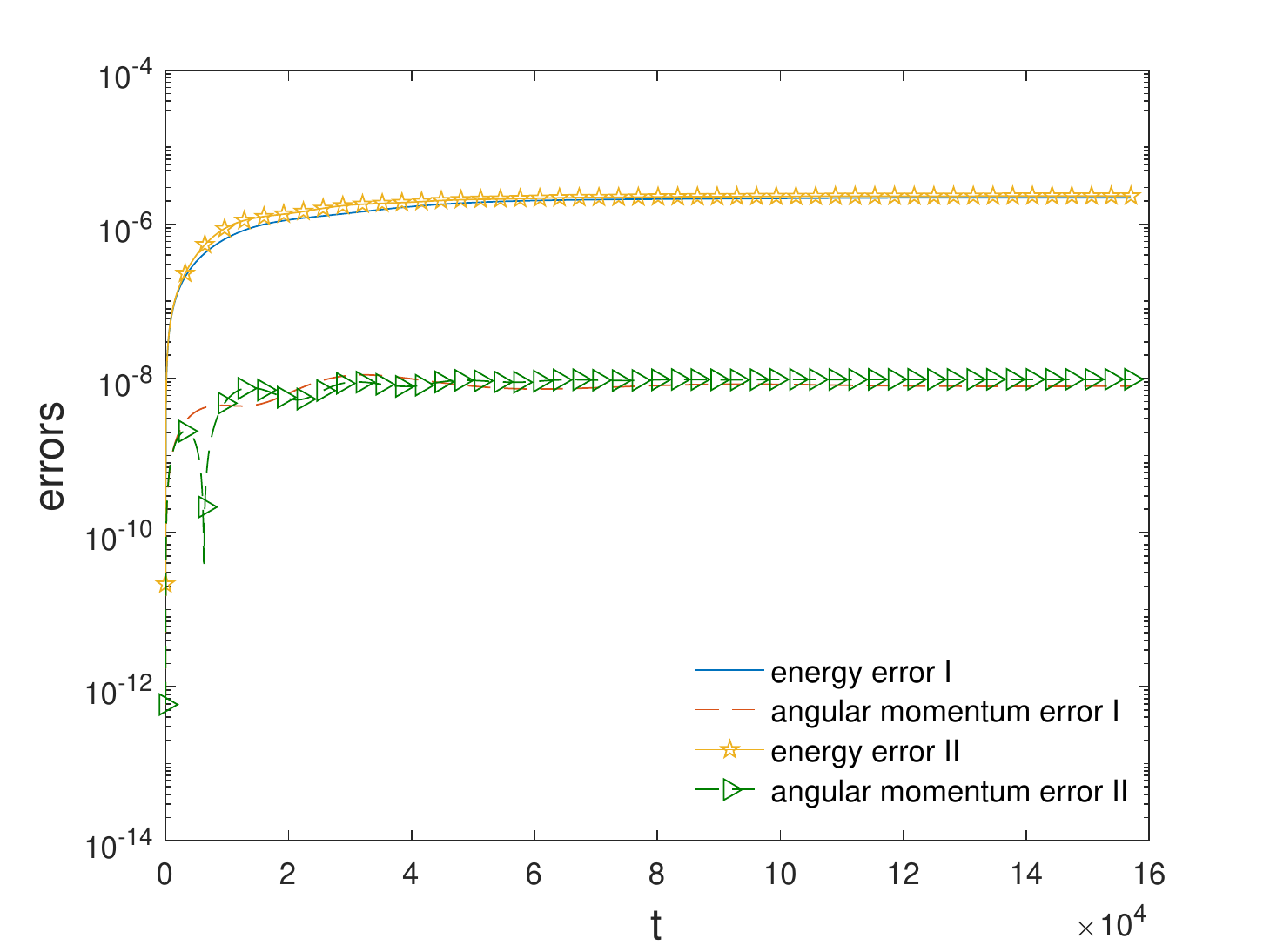}\\
	\includegraphics[width=0.3\linewidth]{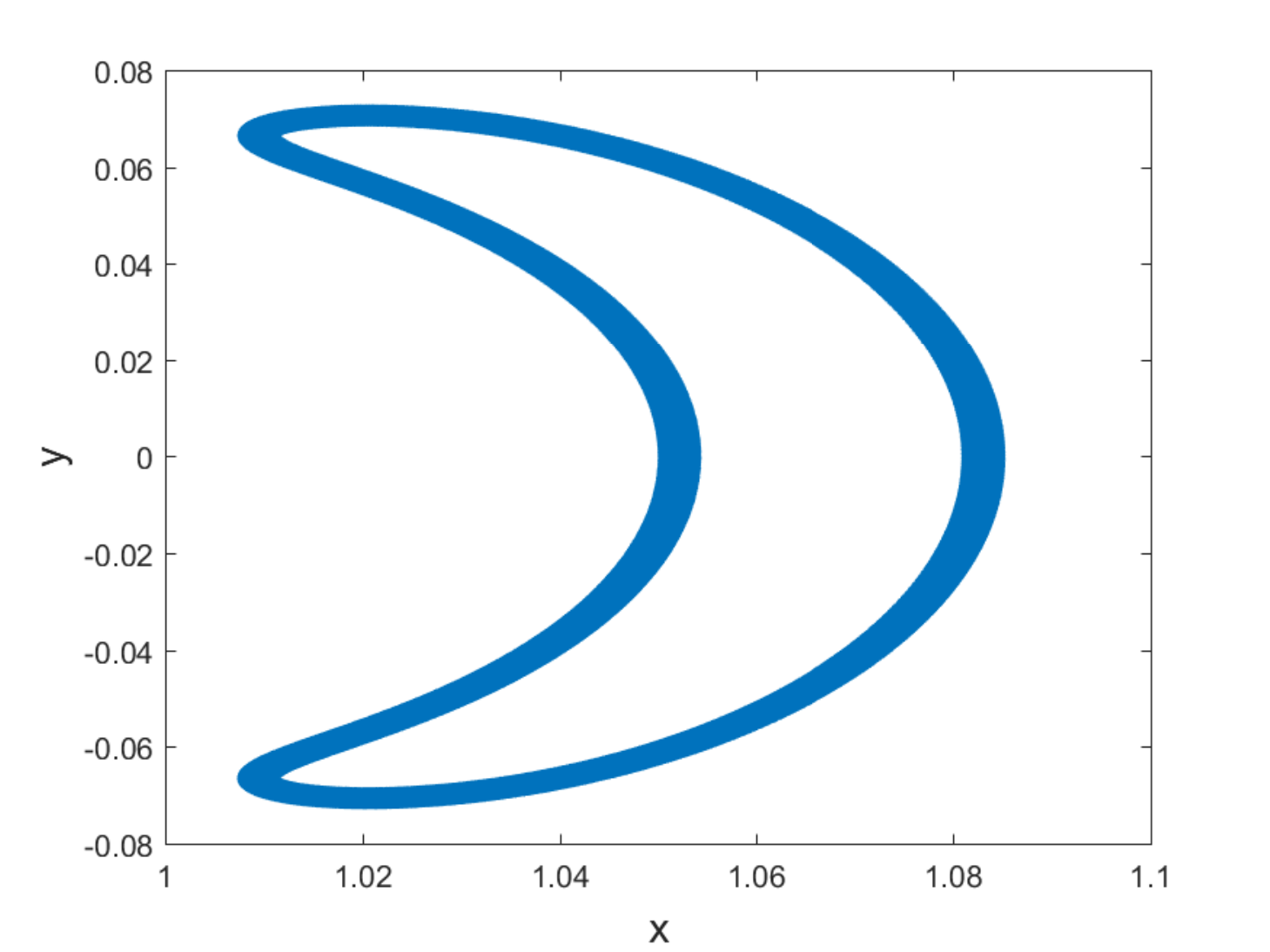}
	\includegraphics[width=0.3\linewidth]{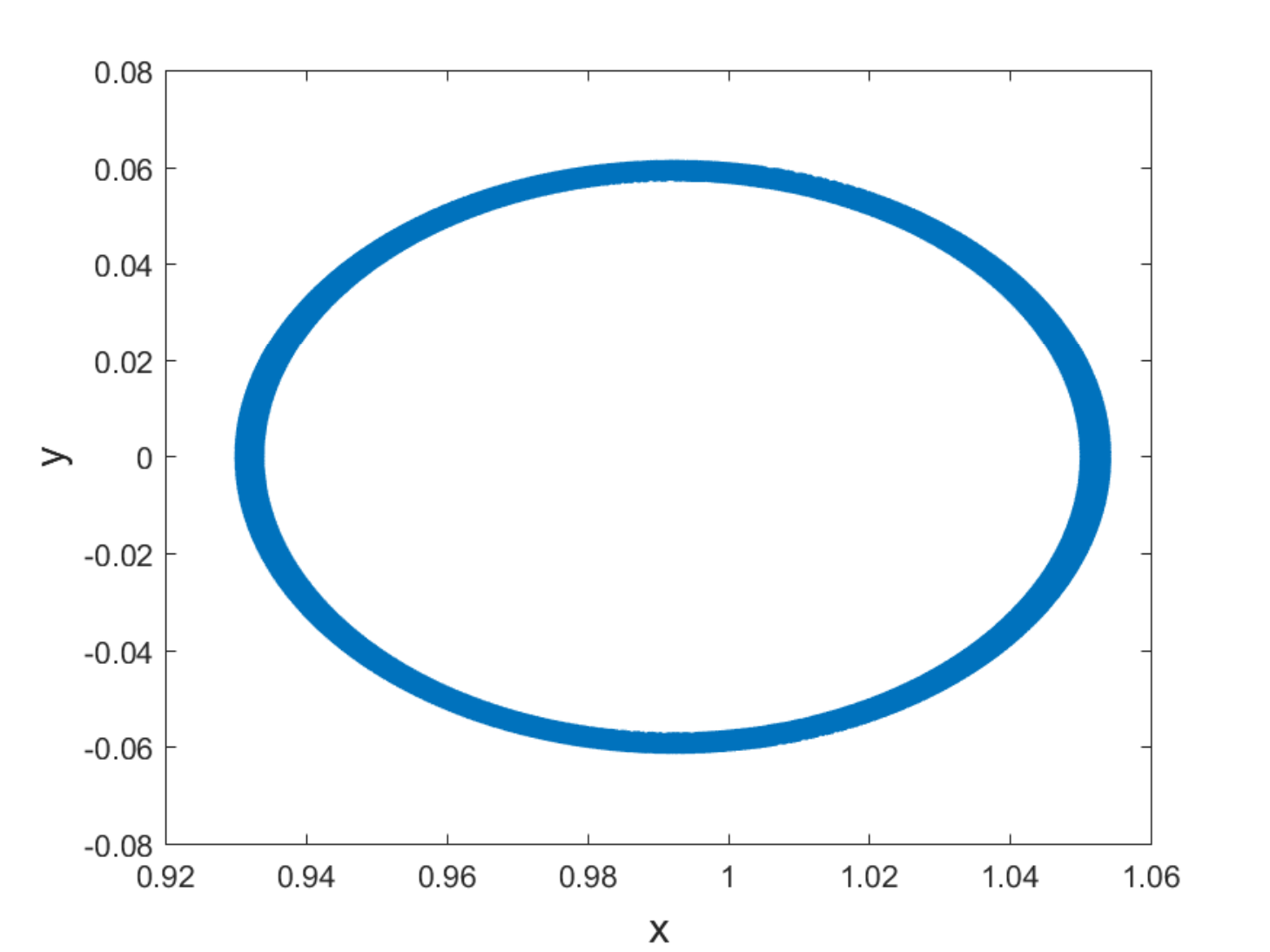}
	\includegraphics[width=0.3\linewidth]{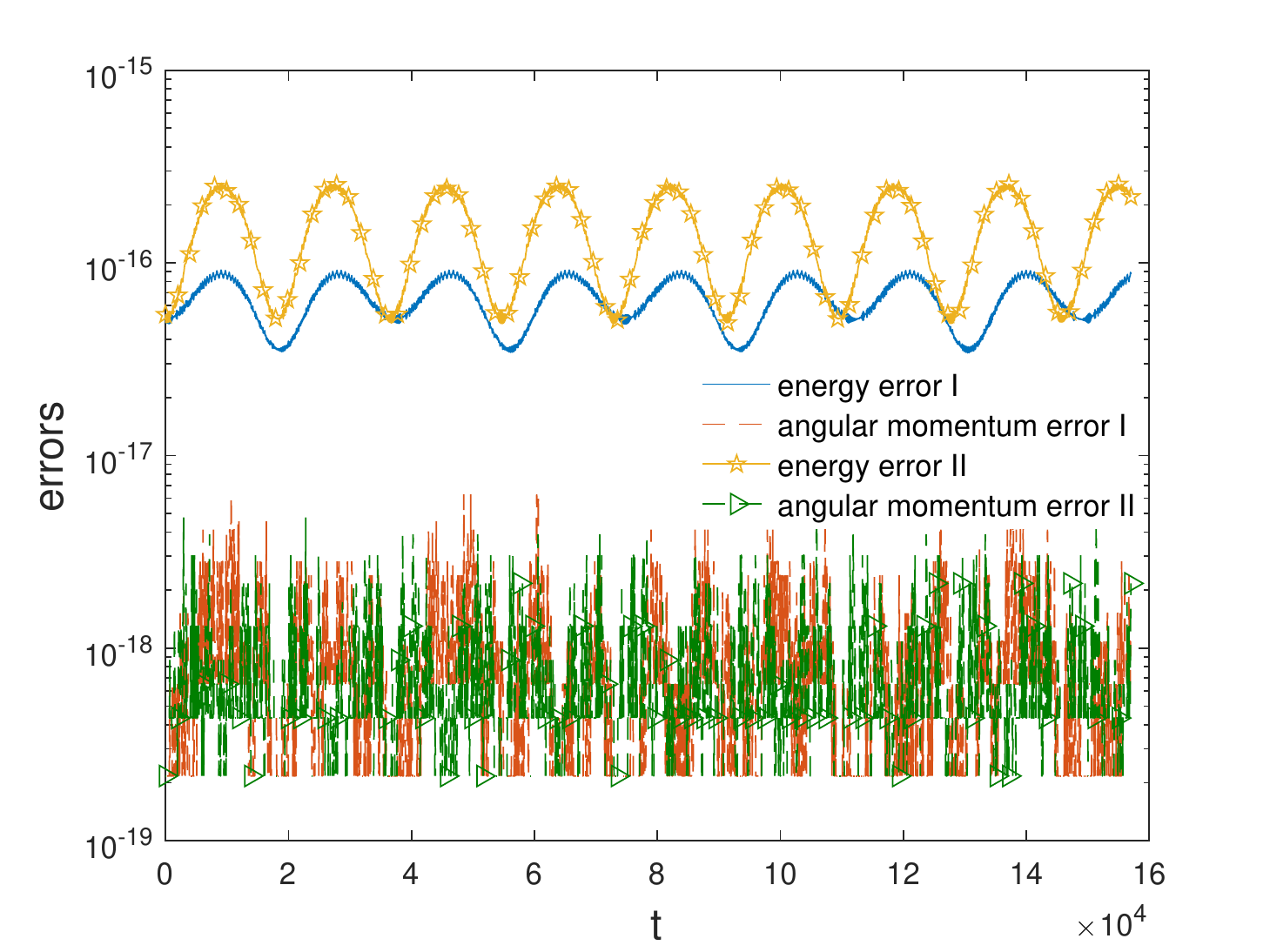}
	\caption{Numerical orbits (banana orbit and transit orbit) and errors in two invariants solved on the interval $[0, 5\times 10^5h]$ by RK4 (upper) and the \textbf{EIP-HL} method (bottom). The legends ended with I and II in the error figures correspond to the banana orbit and transit orbit, respectively.}
\end{figure}

\subsection{Example IV: the rotating Gross-Pitaevskii equation}

Finally, we apply the EIP method to the PDE case to demonstrate its computational efficiency and superior behaviors in invariants preservation. Notice that it is straightforward to extend the EIP method from ODEs to PDEs in the detailed implementation, with the only requirement that the semi-discretization of our targeted PDE is still a conservative system.

Consider the dimensionless time-dependent Gross-Pitaevskii equation (GPE) with a rotating Bose-Einstein condensate (BEC)
\begin{equation}\label{ex4-1}
i\partial_t\psi(t,\bm x)=\left(-\frac{1}{2}\Delta+V(\bm x)-\Omega L_z+\beta|\psi(t,\bm x)|^2\right)\psi(t,\bm x), \quad \bm x\in\mathcal{D}\subset\mathbb{R}^d, ~t\in(0,T],
\end{equation}
for $d=2,3$, where $\psi(t,\bm x)$ is a complex-valued condensate wave function, $V(\bm x)$ is a real-valued potential function. The parameter $\beta$ is the nonlinearity strength representing the interaction between atoms of the condensate. $\Omega$ is the angular velocity and $L_z$ is the $z$-component of the angular momentum
defined by
\[
L_z=-i(x\partial_y-y\partial_x).
\]
In the following experiments, we test the EIP method for the GPE under periodic boundary conditions in both 2D and 3D cases, and the spatial discretization is uniformly taken as the Fourier pseudospectral method (see \cite{cq01,stw11} and references therein). Consider the 2D case for illustration, the resulting semi-discrete scheme can be written as
\begin{equation}\label{ex4-2}
i\frac{d}{dt}\psi_{jk}=\left(-\frac{1}{2}\Delta_h+V_{jk}-\Omega L_z^h+\beta|\psi_{jk}|^2\right)\psi_{jk},\quad 0\leq j\leq J-1, 0\leq k\leq K-1,
\end{equation}
where $i,j$ are the grid indexes and $J, K$ are the partition numbers with respect to $x$ and $y$ directions. $\psi_{jk}$ corresponds to the approximation at the grid point. The discretization of Laplace and angular momentum are defined as
\[
\Delta_h \psi_{jk}=(D_2^x\psi+\psi D_2^y)_{jk},\quad L_z^h\psi_{jk}=-i(X\psi D_1^y-D_1^x\psi Y)_{jk},
\]
where $D_m^x, D_m^y$, $m=1,2$ are the $m$th-order spectral differentiation matrices for $x$ and $y$ directions, respectively.  $X=\mbox{diag}(x_0,x_1,\cdots, x_{J-1})$, and $Y=\mbox{diag}(y_0,y_1,\cdots, y_{K-1})$. It has be proved in \cite{cui-20-GP} that the above semi-discrete scheme \eqref{ex4-2} possesses the mass conservation law
\begin{equation}\label{ex4-3}
\frac{d}{dt}M(t)=0,\quad \mbox{with}\quad M(t)=\|\psi\|_h^2,
\end{equation}
and the energy conservation law
\begin{equation}\label{ex4-4}
\frac{d}{dt}E(t)=0,\quad \mbox{with}\quad E(t)=\frac{1}{2}\|\nabla_h \psi\|^2_h+( V,|\psi|^2)-\Omega(L_z^h\psi,\psi)+\frac{\beta}{2}\|\psi\|_h^4.
\end{equation}
Here, the inner product and the discrete norm for the 2D case are defined as
\[
(u,v)=h_xh_y\sum_{j=0}^{J-1}\sum_{k=0}^{K-1}u_{jk}\bar{v}_{jk}, \quad \|u\|_h=\sqrt{(u,u)},\quad \|\nabla_hu\|_h=\sqrt{-(\Delta_hu,u)},
\]
 where $h_x, h_y$ are the spatial grid sizes. In the following tests, we denote \textbf{EIP-M}, \textbf{EIP-E} and \textbf{EIP-ME} for the EIP methods to preserve the mass, energy and both, respectively.

\subsubsection{Computational efficiency}

When $V(\bm x)=0$ and $\Omega=0$, the GPE \eqref{ex4-1} reduces to the classic nonlinear Schr\"odinger (NLS) equation, which admits an analytical solution as
\[
\psi=A\exp\left(i(\kappa_1x+\kappa_2y-\omega t)\right),\quad \omega=\frac{1}{2}(\kappa_1^2+\kappa_2^2)+\beta A^2.
\]
Thus, we can test the computational efficiency of the EIP methods conveniently. For comparisons, we also present three kinds of conservative schemes for the NLS equation as follows:
\begin{itemize}
	\item the fully-implicit Crank-Nicolson scheme (\textbf{CN}) \cite{wang-13-EE-compact-NLS,gong-17-PS-NLS}:
	\[
	i\frac{\psi_{jk}^{n+1}-\psi_{jk}^n}{h}=\left(-\frac{1}{2}\Delta_h+V_{jk}-\Omega L_z^h\right)\psi_{jk}^{n+1/2}+\frac{\beta}{2}\left(|\psi_{jk}^{n+1}|^2+|\psi_{jk}^{n}|^2\right)\psi_{jk}^{n+1/2};
	\]
	\item the linearly-implicit central  difference scheme  (\textbf{LIC}) \cite{zhang-95-NLS}:
	\[
	i\frac{\psi_{jk}^{n+1}-\psi_{jk}^{n-1}}{2h}=\left(-\frac{1}{2}\Delta_h+V_{jk}-\Omega L_z^h\right)\frac{\psi_{jk}^{n+1}+\psi_{jk}^{n-1}}{2}+\frac{\beta}{2}|\psi_{jk}^{n}|^2\left(\psi_{jk}^{n+1}+\psi_{jk}^{n-1}\right);
	\]
	\item the linearly-implicit Crank-Nicolson scheme based on the SAV approach (\textbf{SAV/CN}) \cite{shen-18-SAV,fu-19-SAV-FNLS}.
\end{itemize}
Though \textbf{LIC} and \textbf{SAV/CN} are both linearly implicit, the corresponding algebraic systems are distinct. \textbf{LIC}  has a variable coefficient matrix while \textbf{SAV/CN} has a constant one, which represent two typical forms of linearly implicit methods.

\begin{figure}[H]
	\centering
	\includegraphics[width=0.45\linewidth]{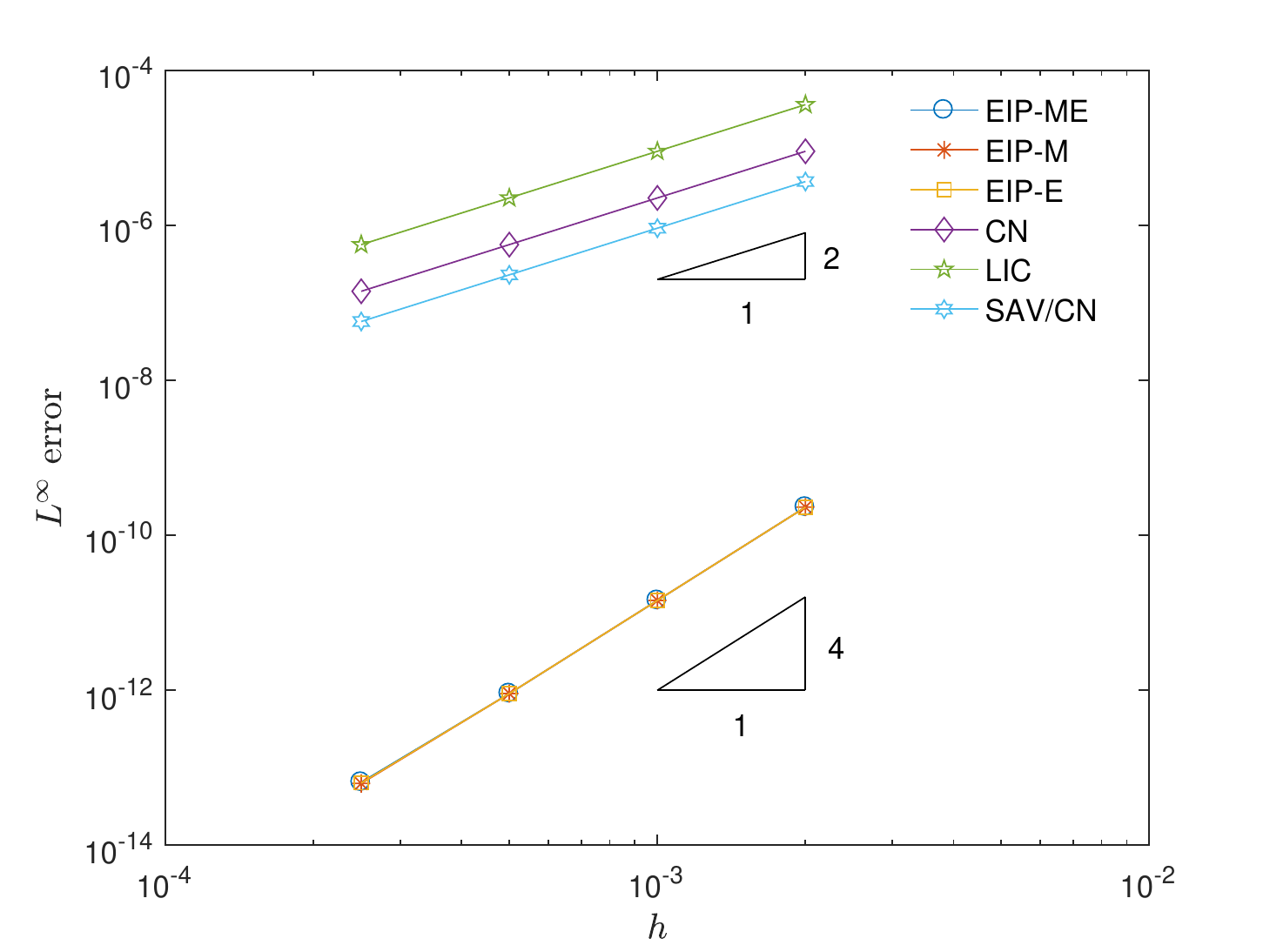}
	\includegraphics[width=0.45\linewidth]{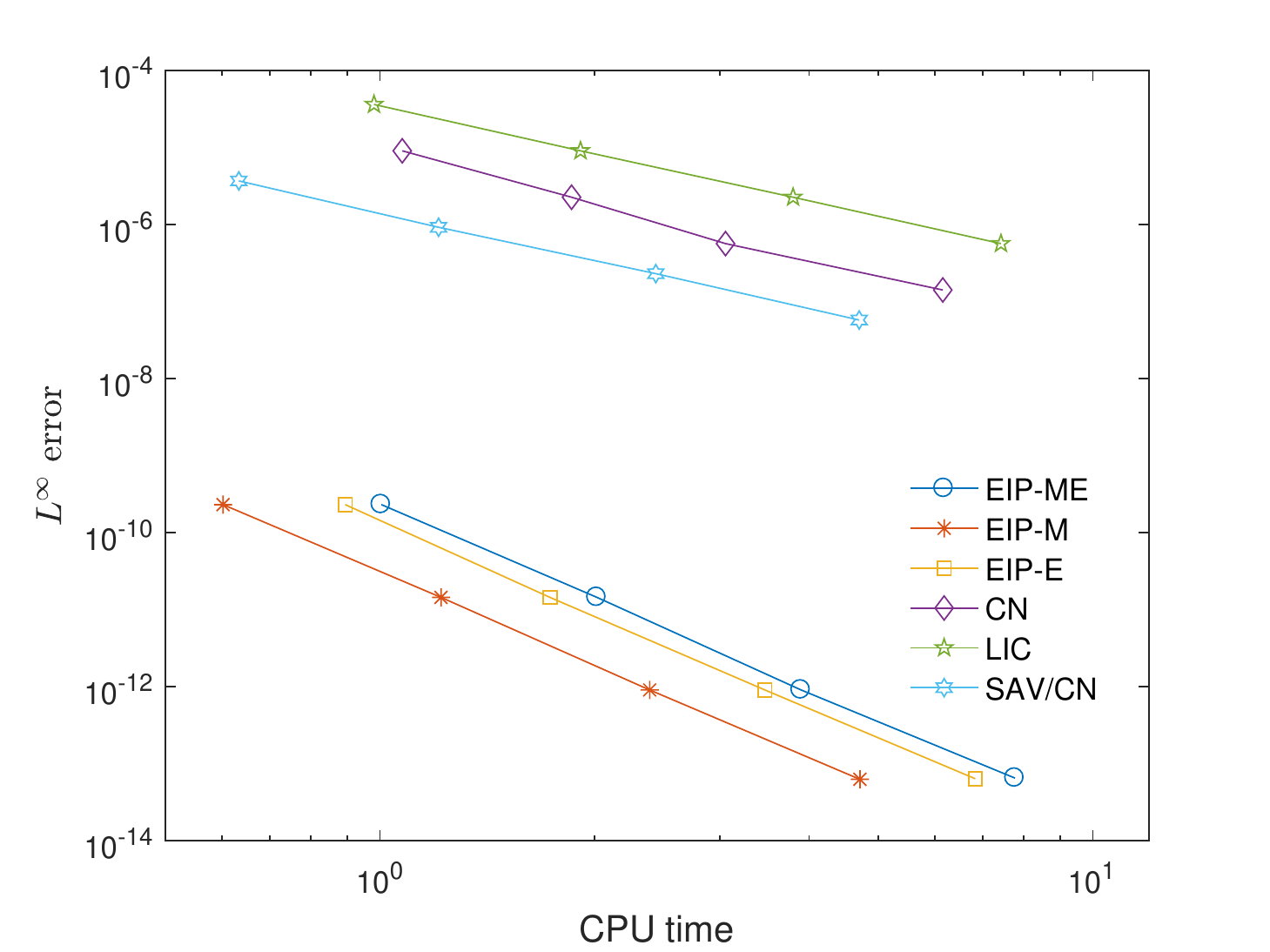}
	\caption{Accuracy tests (left) and computational efficiency (right) for various methods. }\label{ex4-fig1}
\end{figure}

First, we present the convergence test for the EIP method and the above three methods. We can see from Figure.~\ref{ex4-fig1} that all the methods exhibit correct order of accuracy (the line in the left plot overlay each other for the fourth-order EIP methods). Among the three implicit methods, \textbf{CN} is most time consuming due to the nonlinear iterations, while the computational cost of \textbf{SAV/CN} is the cheapest. Though the underlying method of the EIP methods is RK4 with four stages in the Butcher tabular which seems to have more function evaluations, due to their full explicitness, they uniformly give a better performance on the computational efficiency than the presented three methods. Furthermore, CPU times of the single invariant-preserving methods, i.e., \textbf{EIP-M} and \textbf{EIP-E}, are comparable while the \textbf{EIP-ME} method spends a little more time than the formers under the same solution error.

\subsubsection{Dynamics of a rotating BEC}

In the following experiments, we apply the EIP methods on the simulations of dynamics of a rotating BEC for the GPE \eqref{ex4-1}. Firstly, we test the convergence of the EIP methods with different angular velocity $\Omega$ based on the approach \eqref{accuracy}.  We set $\mathcal{D}=[-2,2]^2$, $V(\bm x)=\frac{1}{2}(x^2+y^2)$ and $\beta=1$. The initial condition is taken as $\psi_0=\frac{2}{\sqrt{\pi}}(x+iy)\exp(-8(x^2+y^2))$. Table.~\ref{ex4-tab1} lists the accuracy results for the three methods where $h_0=0.0002$ is the initial time step. It is clearly that all the methods can also achieve an excepted convergence order, independent of the angular velocity $\Omega$.

\begin{table}[H]
	\centering
	\begin{tabular*}{0.9\textwidth}[h]{@{\extracolsep{\fill}}llcccccc} \hline
		& \multirow{2}{*}{$h$}	& \multicolumn{2}{c}{\textbf{EIP-M}} &\multicolumn{2}{c}{\textbf{EIP-E}} & \multicolumn{2}{c}{\textbf{EIP-ME}} \\ \cline{3-4}\cline{5-6}\cline{7-8}
		& & $L^\infty$-error & order & $L^\infty$-error & order & $L^\infty$-error & order \\ \hline
		\multirow{5}{*}{$\Omega=0$}		& $h_0$  &  &  &  & & & \\[1ex]
		& $h_0/2$ & 9.0794e-11  &    &   9.0712e-11 &   &    9.0696e-11 &  \\[1ex]
		& $h_0/4$ &  5.6754e-12  & 3.9998 &    5.6727e-12 &    3.9992 &  5.6722e-12 & 3.9990 \\[1ex]
		& $h_0/8$ & 3.5474e-13    &  3.9999  &  3.5454e-13 &   4.0000 &   3.5443e-13 &  4.0004 \\ [1ex] 	\hline
		\multirow{5}{*}{$\Omega=0.5$} & $h_0$  &  &  &  & & & \\[1ex]
		& $h_0/2$ & 9.3456e-11  &    &   9.3455e-11  &   &    9.3448e-11 &  \\[1ex]
		& $h_0/4$ &  5.8442e-12  & 3.9992 &    5.8441e-12 &   3.9992 &  5.8439e-12 & 3.9992 \\[1ex]
		& $h_0/8$ & 3.6541e-13    &   3.9994  &  3.6528e-13 &   3.9999 &   3.6539e-13 &  3.9994\\ [1ex] 	\hline
		\multirow{5}{*}{$\Omega=0.9$} & $h_0$  &  &  &  & & & \\[1ex]
		& $h_0/2$ & 1.1316e-10  &    &   1.1315e-10 &   &    1.1309e-10 &  \\[1ex]
		& $h_0/4$ &  7.0722e-12  & 4.0001 &    7.0715e-12 &    4.0000 &  7.0698e-12 & 3.9996 \\[1ex]
		& $h_0/8$ & 4.4213e-13    &  3.9996 &  4.4224e-13 &   3.9991 &   4.4220e-13 &  3.9989 \\ [1ex] \hline	
	\end{tabular*}
	\caption{Errors between two adjacent time steps and the corresponding orders with three different $\Omega$ for the EIP methods. The partition numbers in $x$ and $y$ directions are set to $J=K=128$ and the computational time is $t=0.5$.}\label{ex4-tab1}
\end{table}

\begin{figure}[H]
	\centering
	\subfigure[$t=0$]{
		\begin{minipage}[b]{0.3\linewidth}
			\includegraphics[width=1\linewidth]{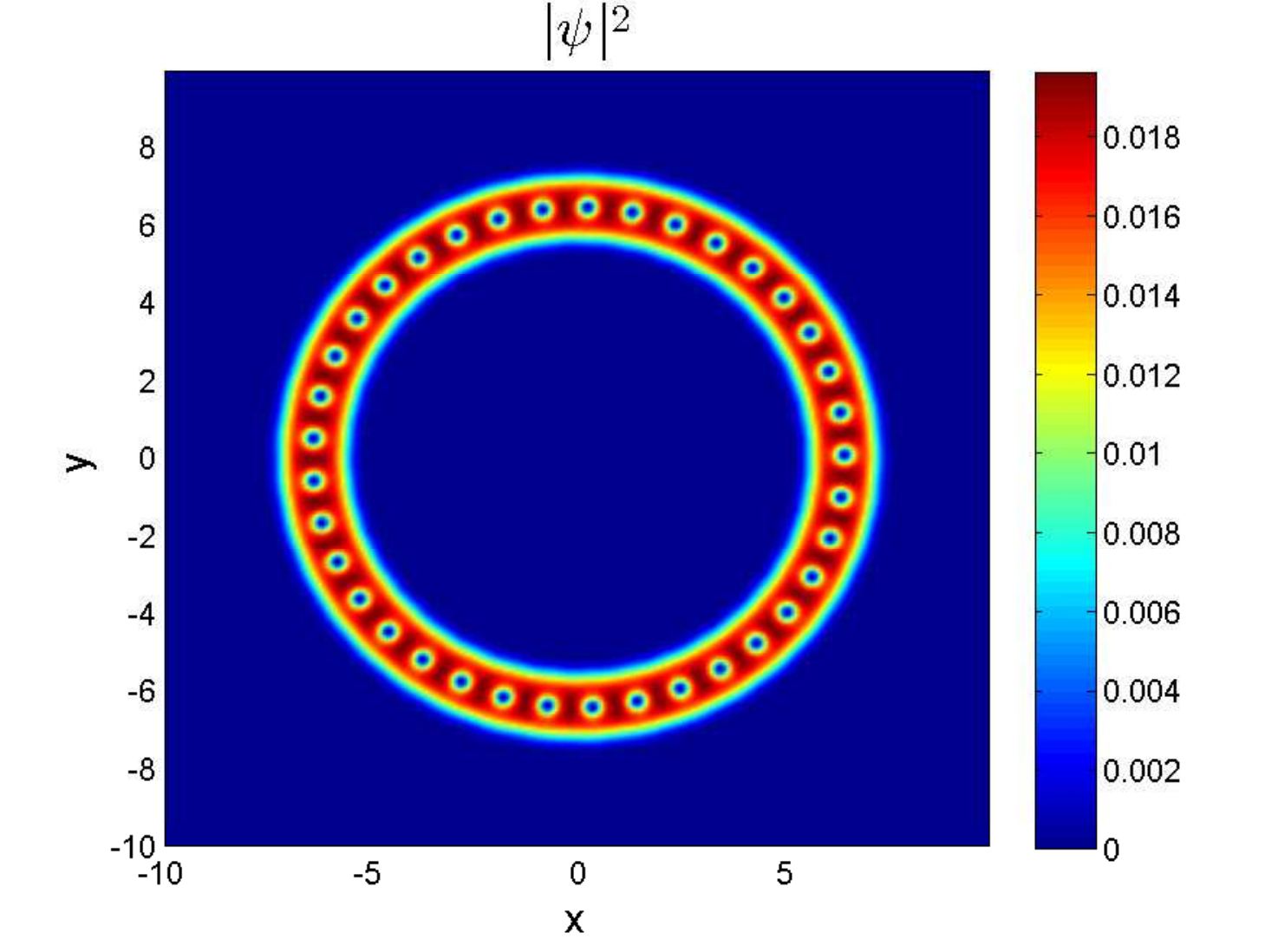}
	\end{minipage}}
	\subfigure[$t=0.24$]{
		\begin{minipage}[b]{0.3\linewidth}
			\includegraphics[width=1\linewidth]{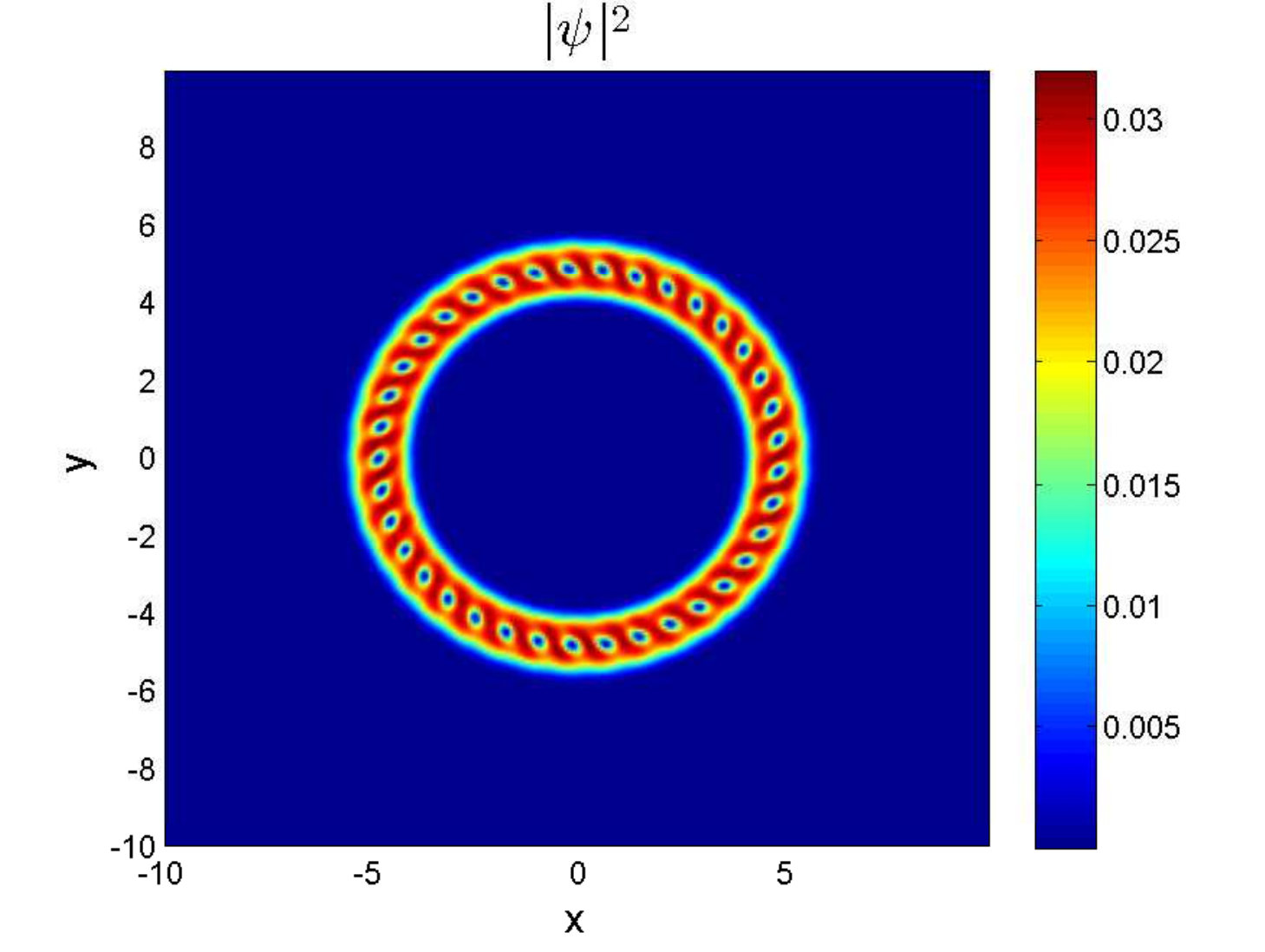}
	\end{minipage}}
	\subfigure[$t=0.44$]{
	\begin{minipage}[b]{0.3\linewidth}
		\includegraphics[width=1\linewidth]{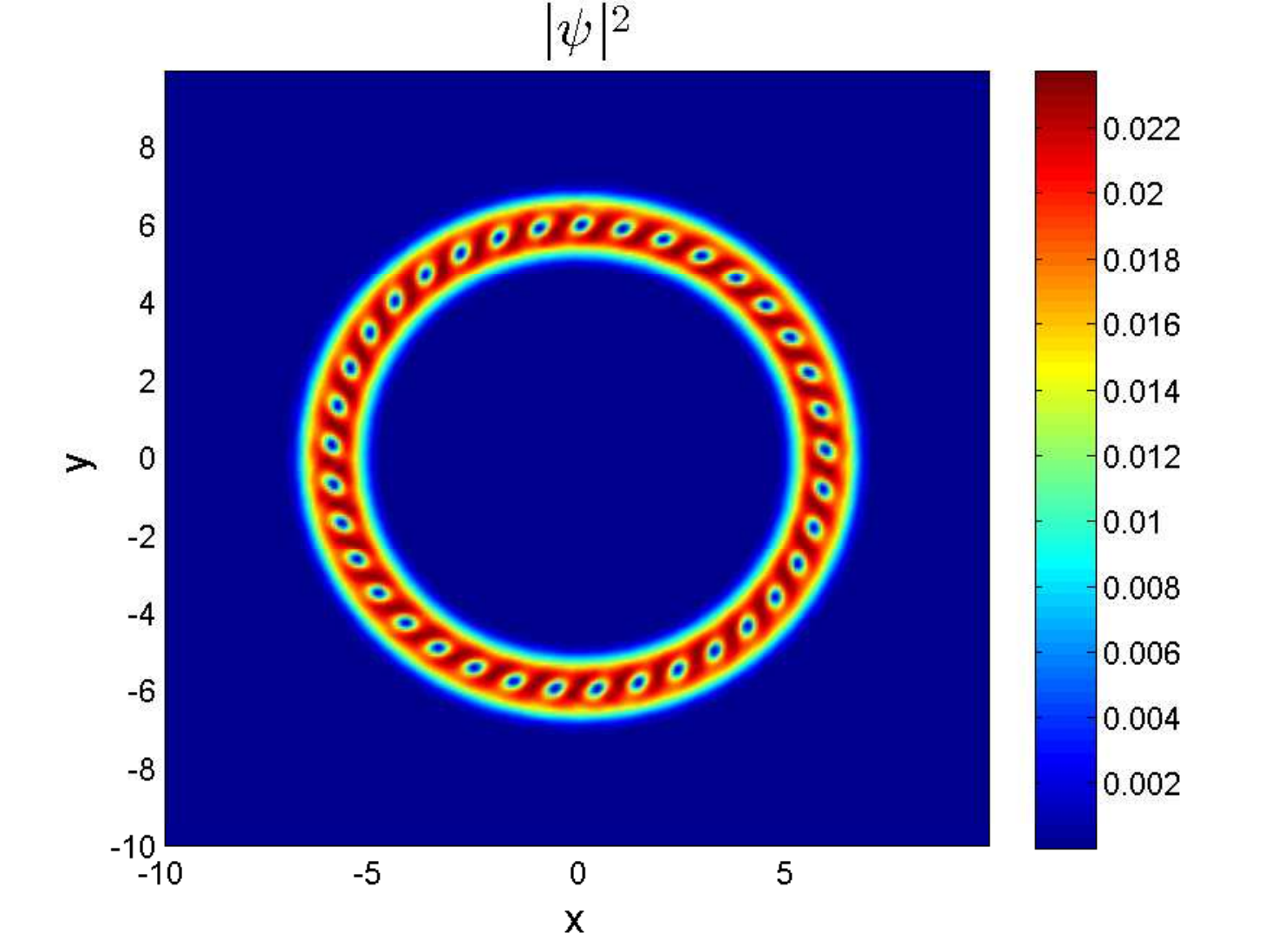}
	\end{minipage}}
	\subfigure[$t=0.72$]{
		\begin{minipage}[b]{0.3\linewidth}
			\includegraphics[width=1\linewidth]{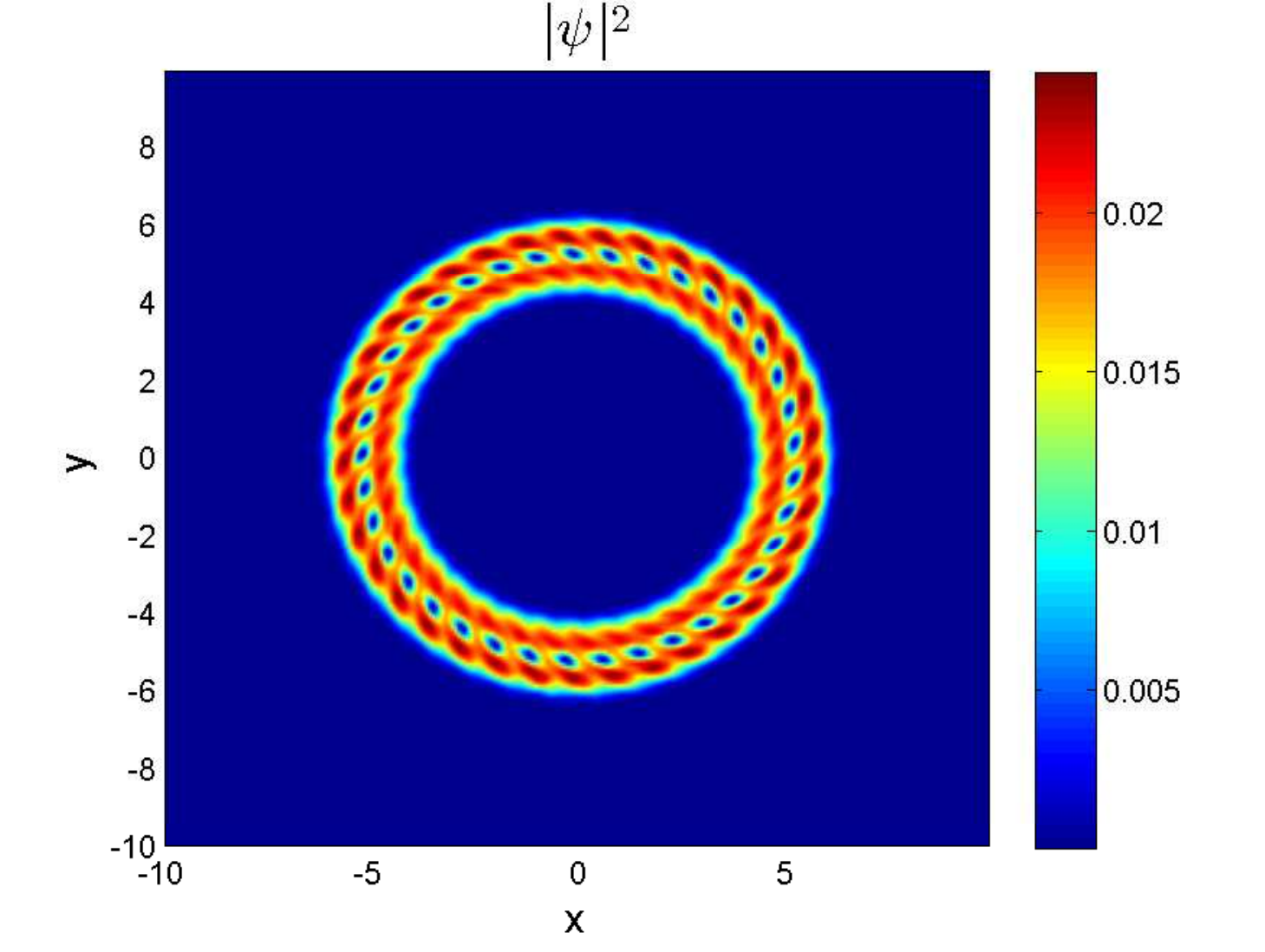}
	\end{minipage}}
		\subfigure[$t=0.8$]{
		\begin{minipage}[b]{0.3\linewidth}
			\includegraphics[width=1\linewidth]{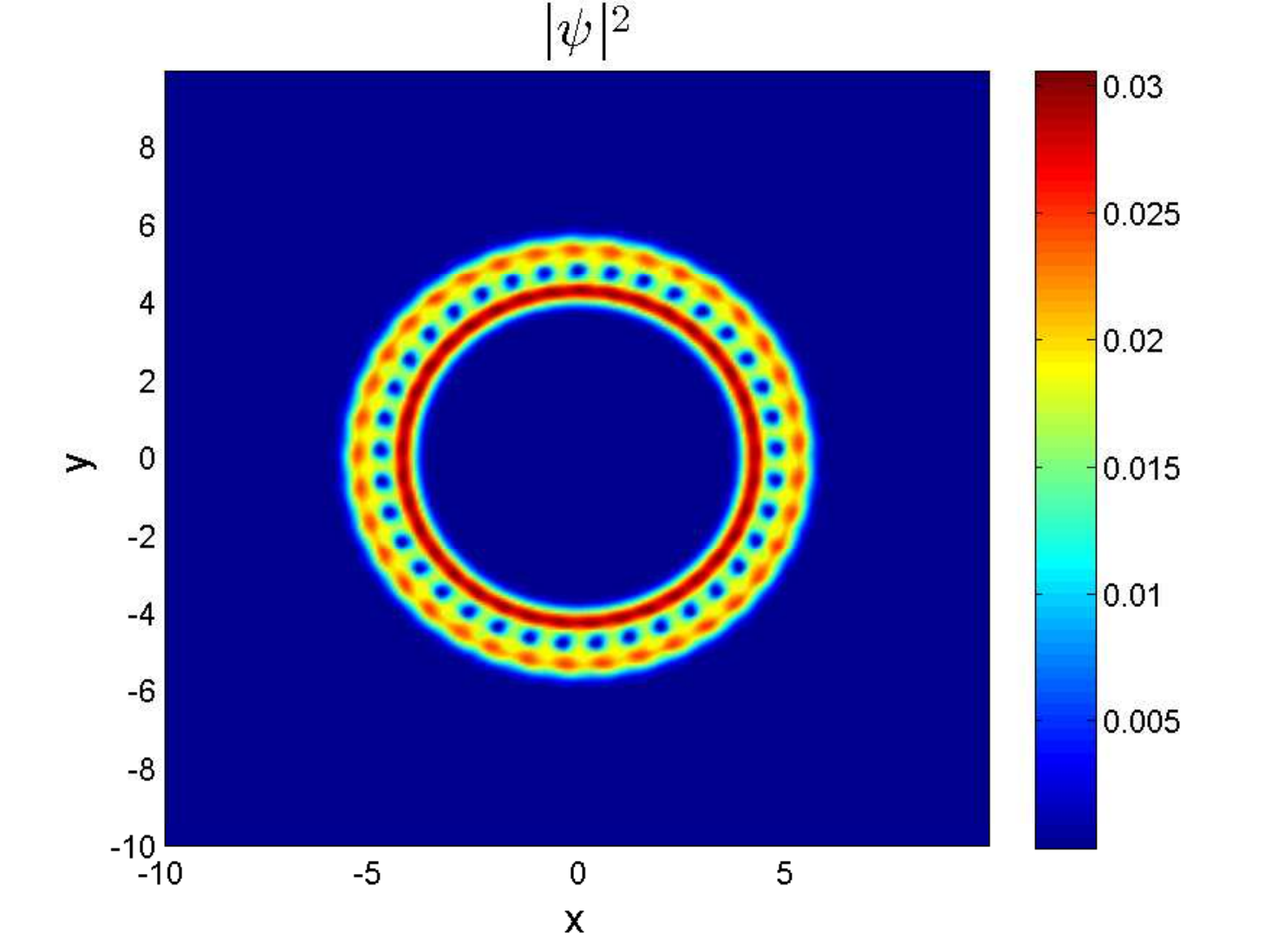}
	\end{minipage}}
		\subfigure[$t=1$]{
		\begin{minipage}[b]{0.3\linewidth}
			\includegraphics[width=1\linewidth]{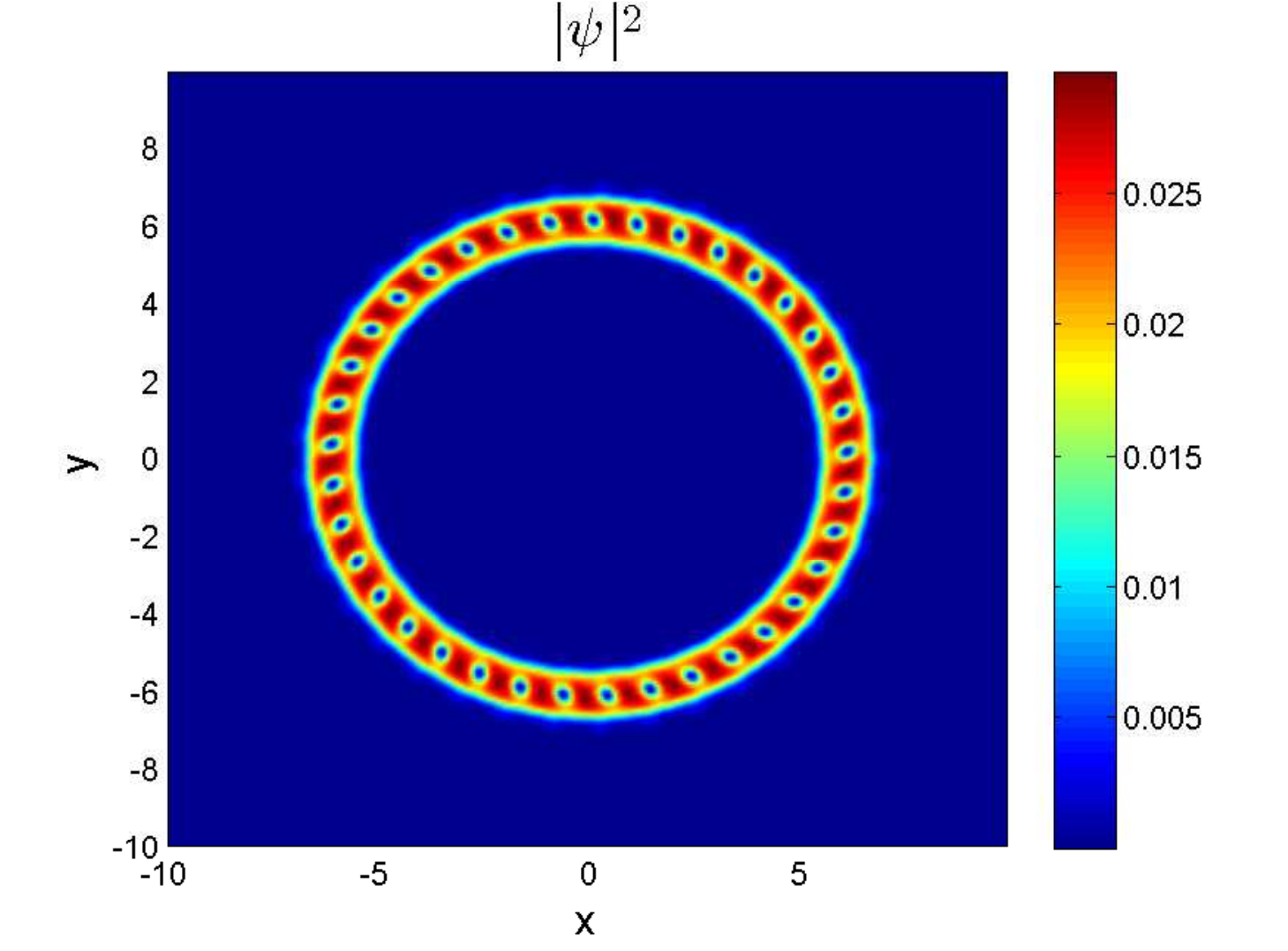}
	\end{minipage}}
	\caption{Snapshots of the density $|\psi|^2$ computed by the \textbf{EIP-M} method for the ring BEC.}\label{ex4-fig2}
\end{figure}

\begin{figure}[H]
	\centering
	\subfigure[\textbf{EIP-M}]{
		\begin{minipage}[b]{0.3\linewidth}
			\includegraphics[width=1\linewidth]{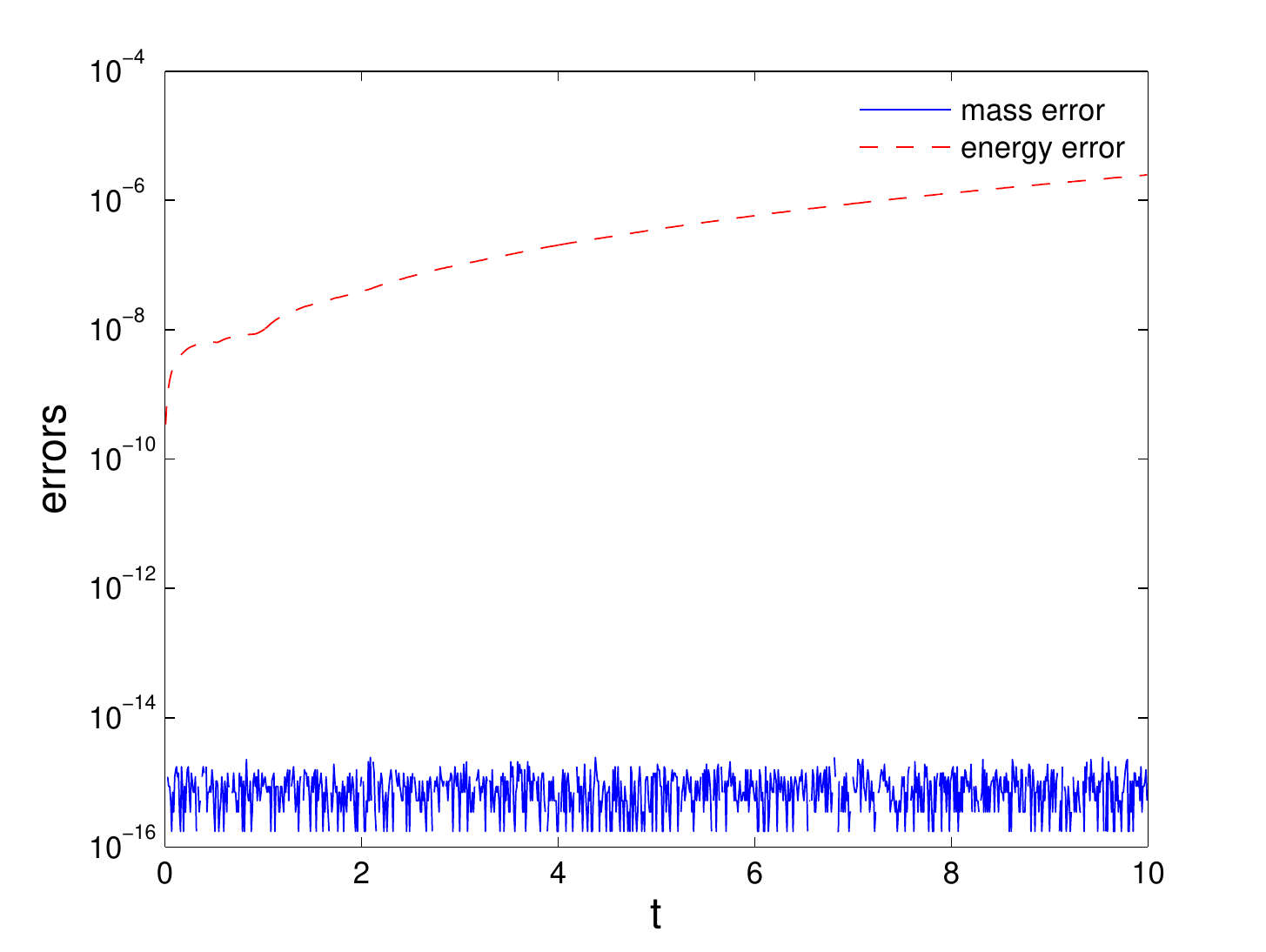}
	\end{minipage}}
	\subfigure[\textbf{EIP-E}]{
		\begin{minipage}[b]{0.3\linewidth}
			\includegraphics[width=1\linewidth]{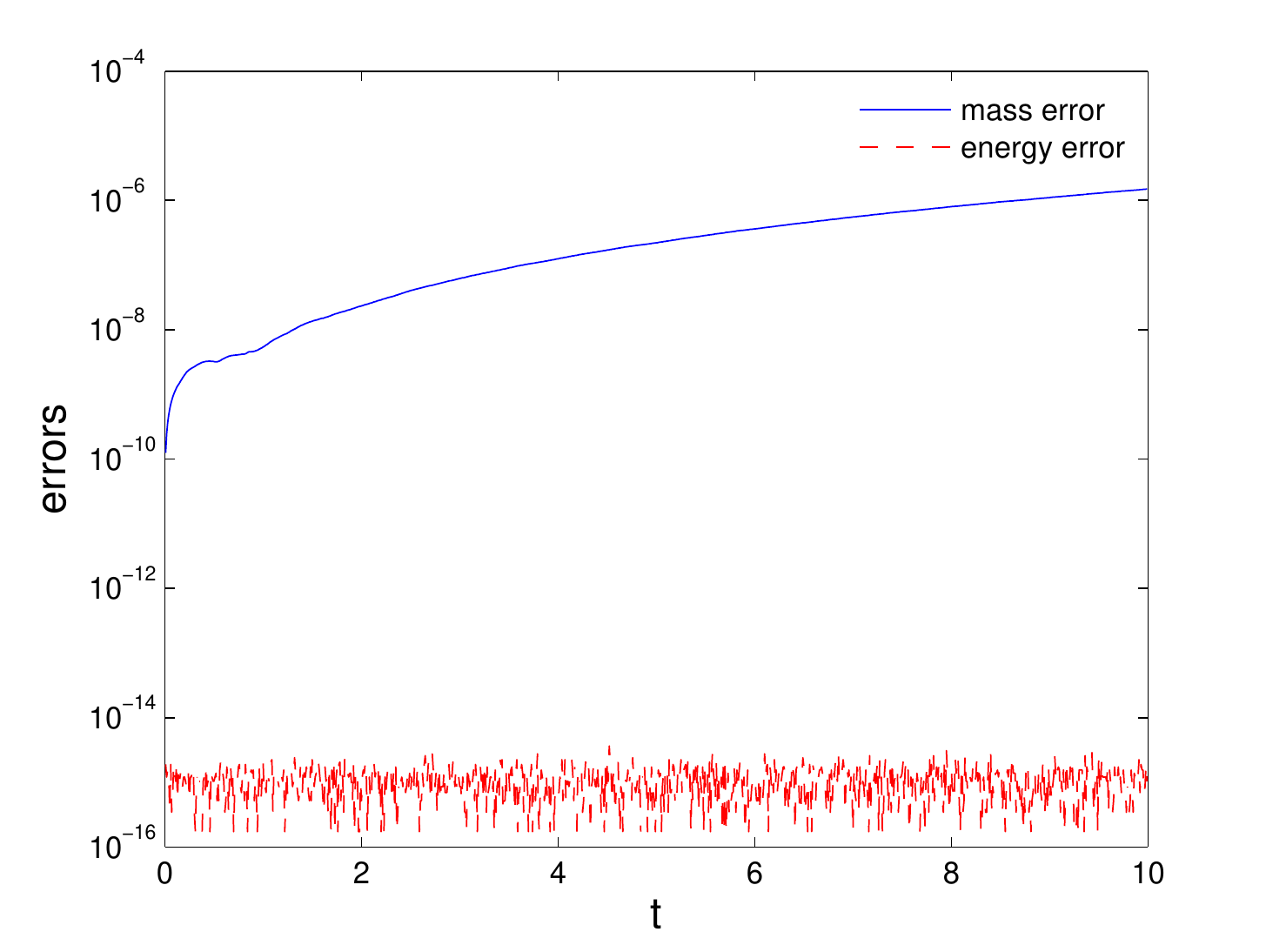}
	\end{minipage}}
	\subfigure[\textbf{EIP-ME}]{
		\begin{minipage}[b]{0.3\linewidth}
			\includegraphics[width=1\linewidth]{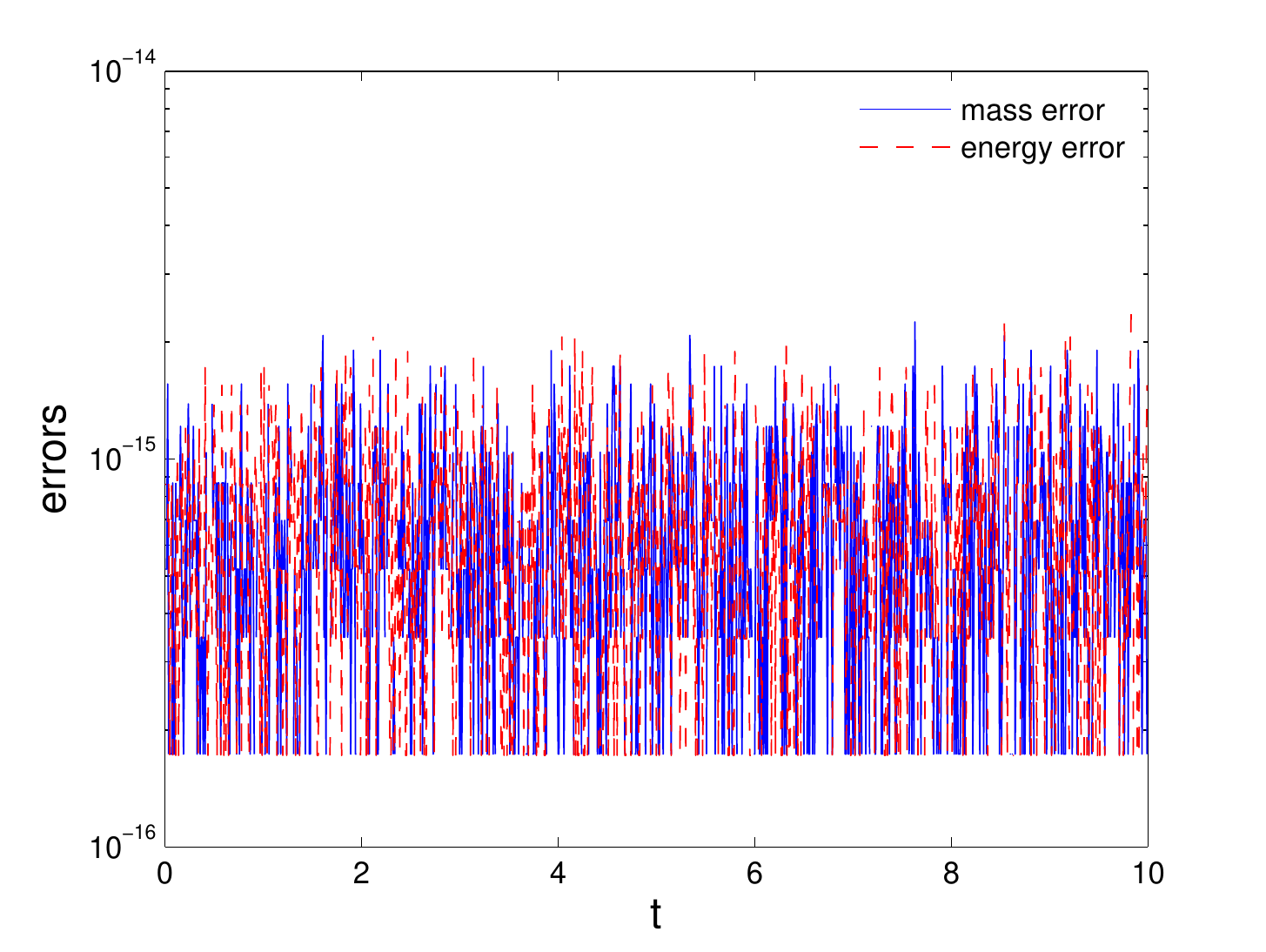}
	\end{minipage}}
	\caption{Errors in two invariants computed by the EIP methods for the ring BEC.}\label{ex4-fig3}
\end{figure}

Next, we consider the dynamics of vortex lattices in rotating BECs with $\Omega=3.5$, $\beta=1000$. The domain $\mathcal{D}=[-10,10]^2$, the partition numbers $J=K=256$ and $h=0.0002$.  The initial datum is chosen as the $L^2$-normalized ground state eigenvector of the Gross-Pitaevskii operator $G_0(v):=\left(-\frac{1}{2}\Delta +V_0(\bm x)-\Omega L_z+\beta|v|^2\right)v$ where the potential function $V_0$ is set to a quadratic-plus-quartic potential \cite{ad15b} as
\begin{equation}
V_0(\bm x)=\frac{1-\alpha}{2}(\gamma_x^2 x^2+\gamma_y^2 y^2)+\frac{\kappa}{4}(\gamma_x^2 x^2+\gamma_y^2 y^2)^2,
\end{equation}
with $\gamma_x=\gamma_y=1$, $\alpha=1.2$ and $\kappa=0.3$. The practical computation of the ground state is done by the Matlab package named GPELab \cite{ad14,ad15a}. This stationary state is a circular ring with many uniformly distributed vortices (see the first plot in Figure.~\ref{ex4-fig2}). To simulate the dynamics in this ring BEC, we perturb the parameter $\kappa=0.3$ to $0.7$ and generate a typical example of a fast rotating BEC which demands high-precision numerical methods to capture the movement of each vortex. Figure.~\ref{ex4-fig2} displays the snapshots of solutions obtained by the \textbf{EIP-M} method, and the solutions computed by the \textbf{EIP-E} and \textbf{EIP-ME} methods look the same. We can observe a complex dynamics in the ring BEC. All the vortices exhibit a clockwise rotation and the number of which is also conserved during the simulation. We further demonstrate the long-time behavior by carrying out a larger time period $t=10$. As we can see from Figure.~\ref{ex4-fig3}, the error in the mass or energy invariant is preserved exactly by the \textbf{EIP-M} or \textbf{EIP-H} method respectively, while both errors reach machine accuracy by the \textbf{EIP-ME} method.

\begin{figure}[H]
	\centering
	\subfigure[$t=0$]{
		\begin{minipage}[b]{0.3\linewidth}
			\includegraphics[width=1\linewidth]{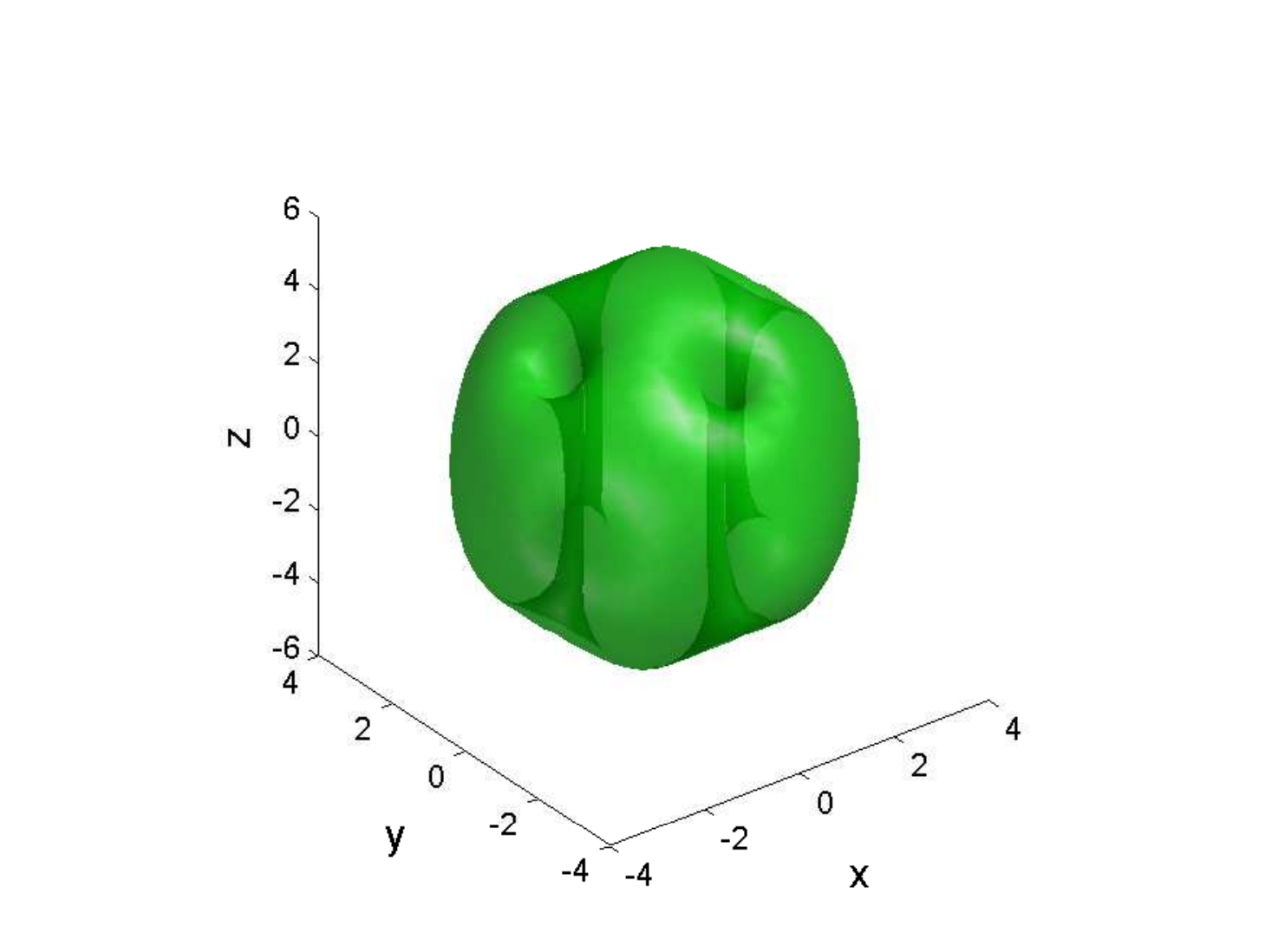}
	\end{minipage}}
	\subfigure[$t=4$]{
		\begin{minipage}[b]{0.3\linewidth}
			\includegraphics[width=1\linewidth]{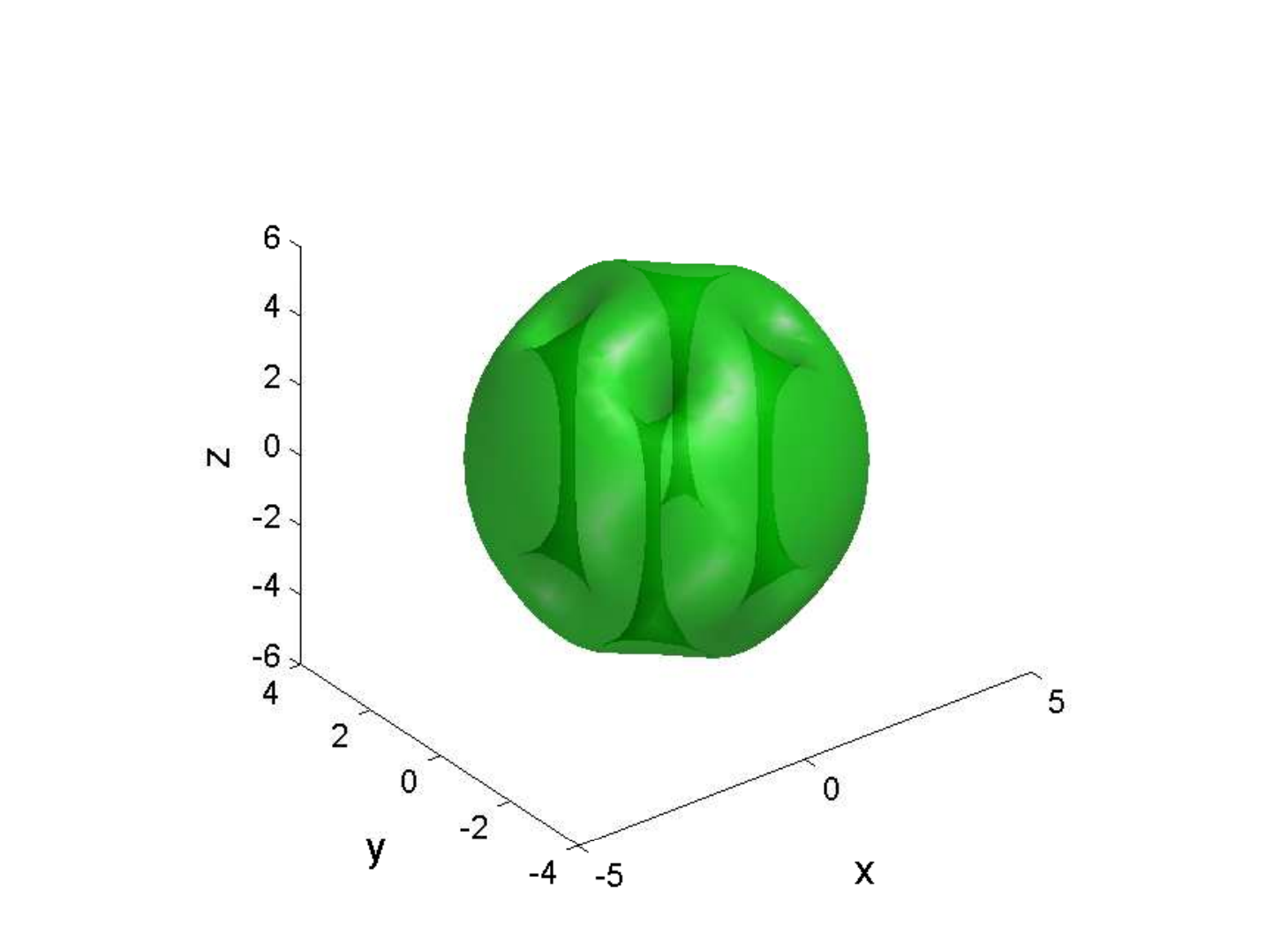}
	\end{minipage}}
	\subfigure[$t=8$]{
		\begin{minipage}[b]{0.3\linewidth}
			\includegraphics[width=1\linewidth]{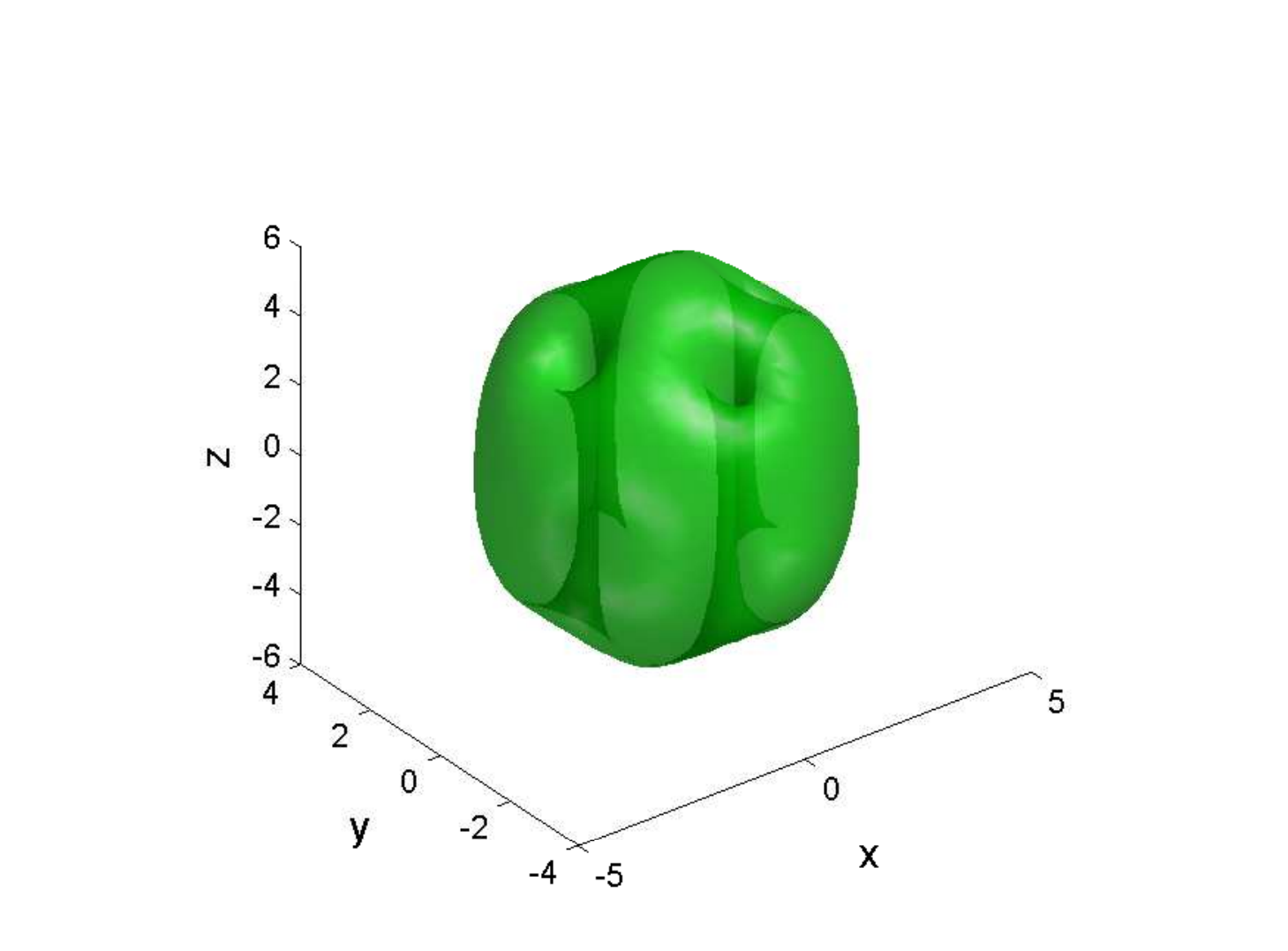}
	\end{minipage}}
	\subfigure[$t=12$]{
		\begin{minipage}[b]{0.3\linewidth}
			\includegraphics[width=1\linewidth]{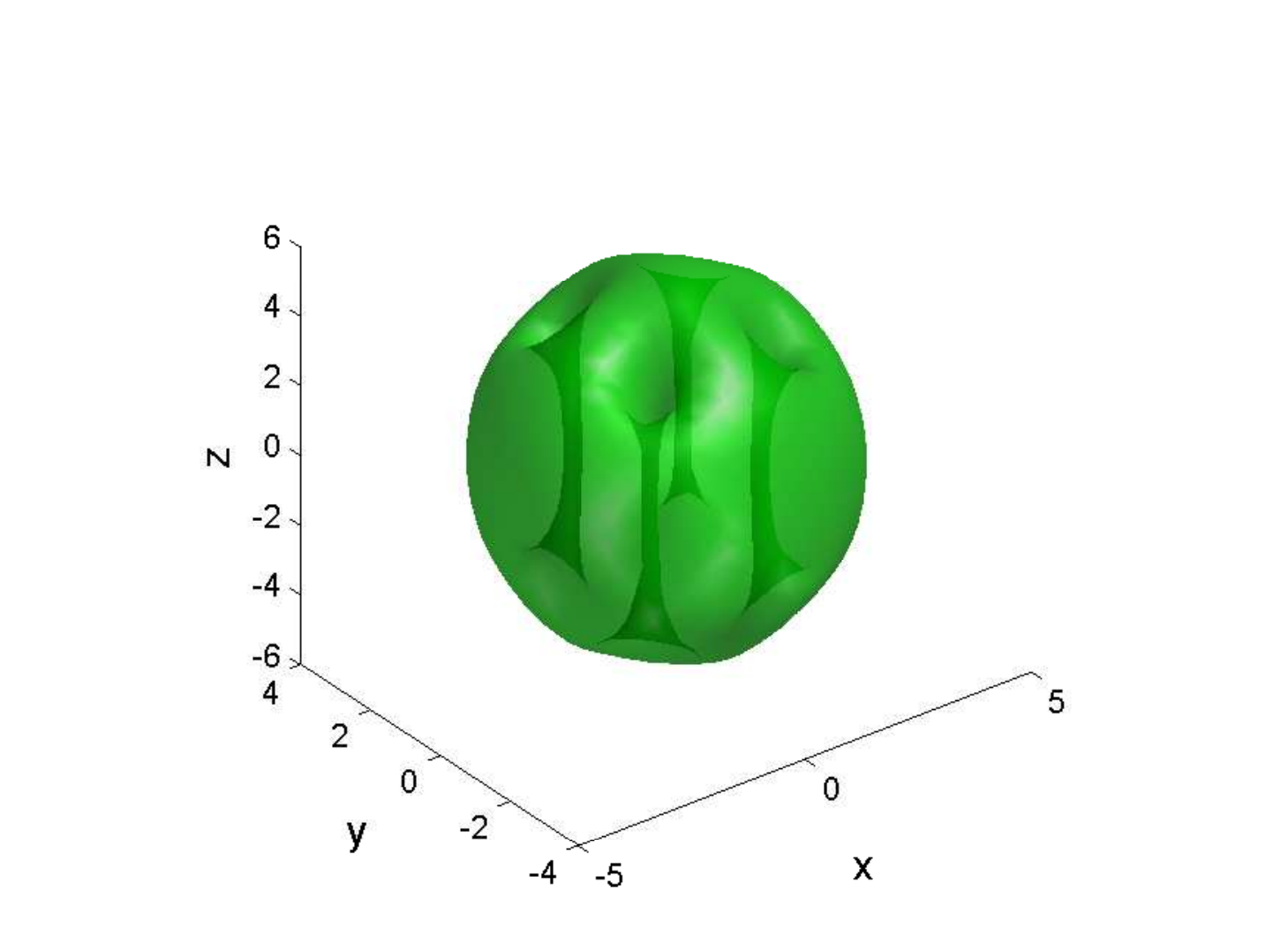}
	\end{minipage}}
	\subfigure[$t=16$]{
		\begin{minipage}[b]{0.3\linewidth}
			\includegraphics[width=1\linewidth]{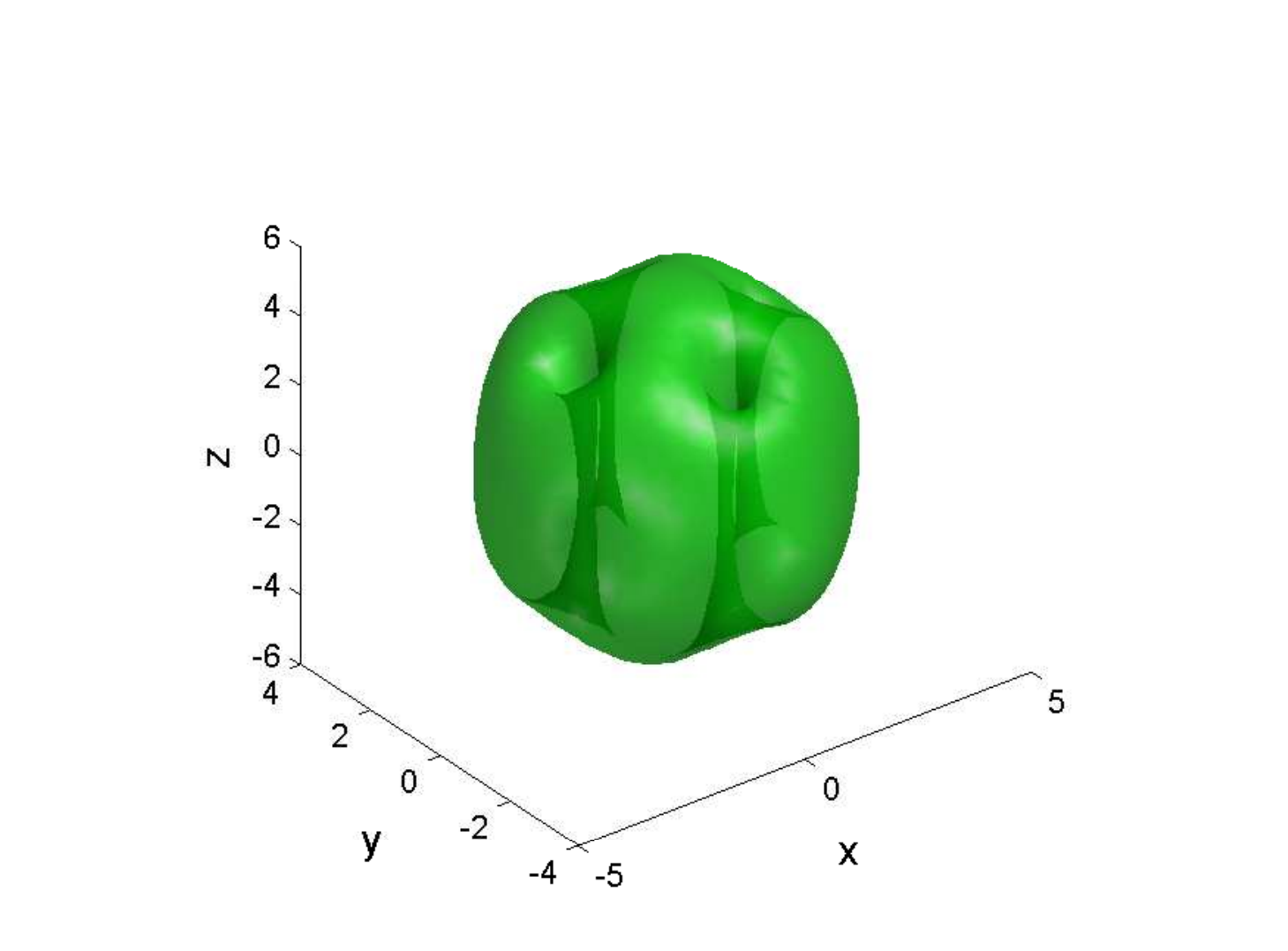}
	\end{minipage}}
	\subfigure[$t=20$]{
	\begin{minipage}[b]{0.3\linewidth}
			\includegraphics[width=1\linewidth]{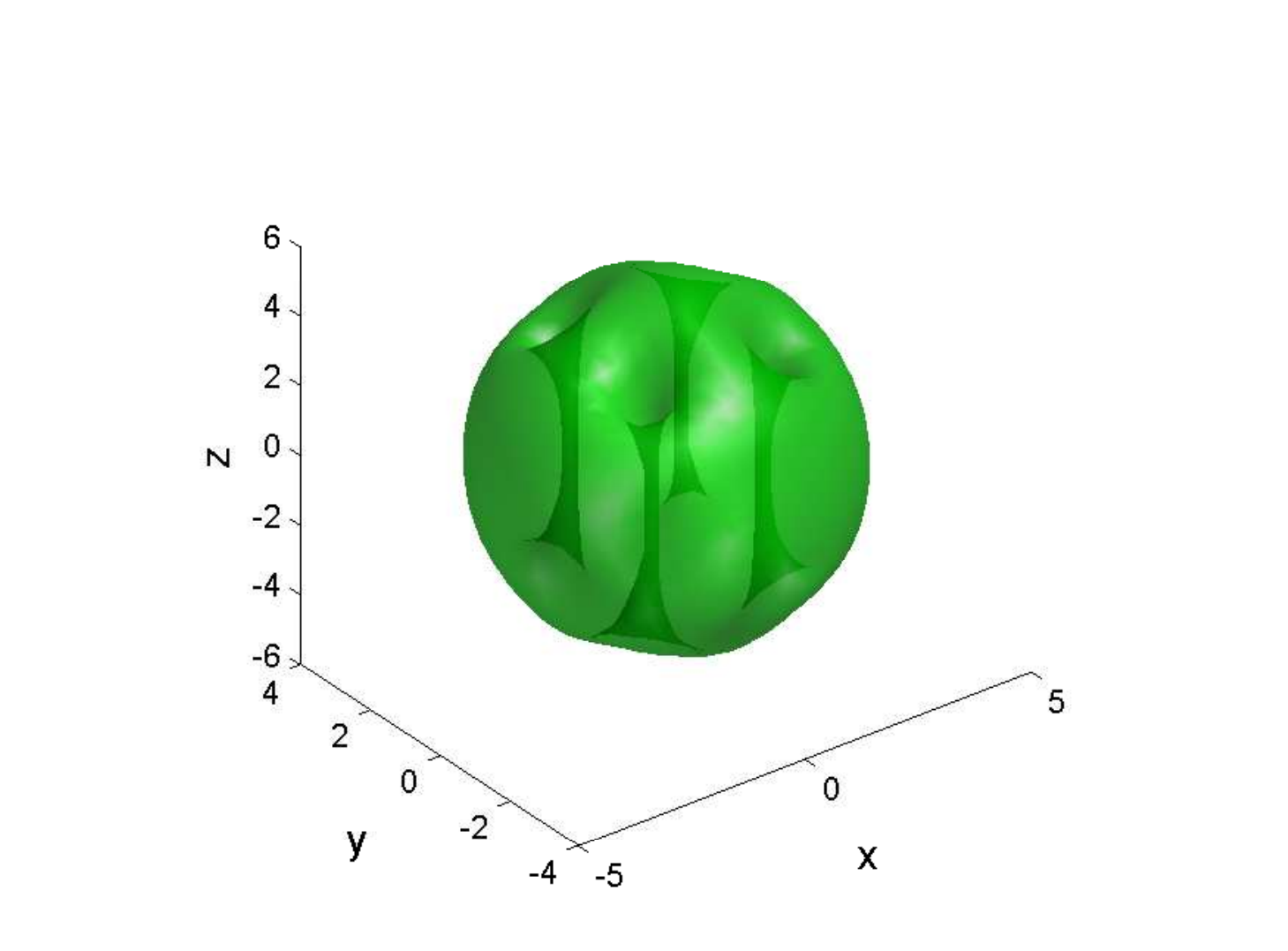}
	\end{minipage}}
	\caption{Snapshots of $10^{-3}$-isosurfaces of $|\psi|^2$ at different times computed by the \textbf{EIP-M} method for the dynamics of vortex lines.}\label{ex4-fig4}
\end{figure}

\begin{figure}[H]
 	\centering
 	\subfigure[$t=0$]{
 		\begin{minipage}[b]{0.3\linewidth}
 			\includegraphics[width=1\linewidth]{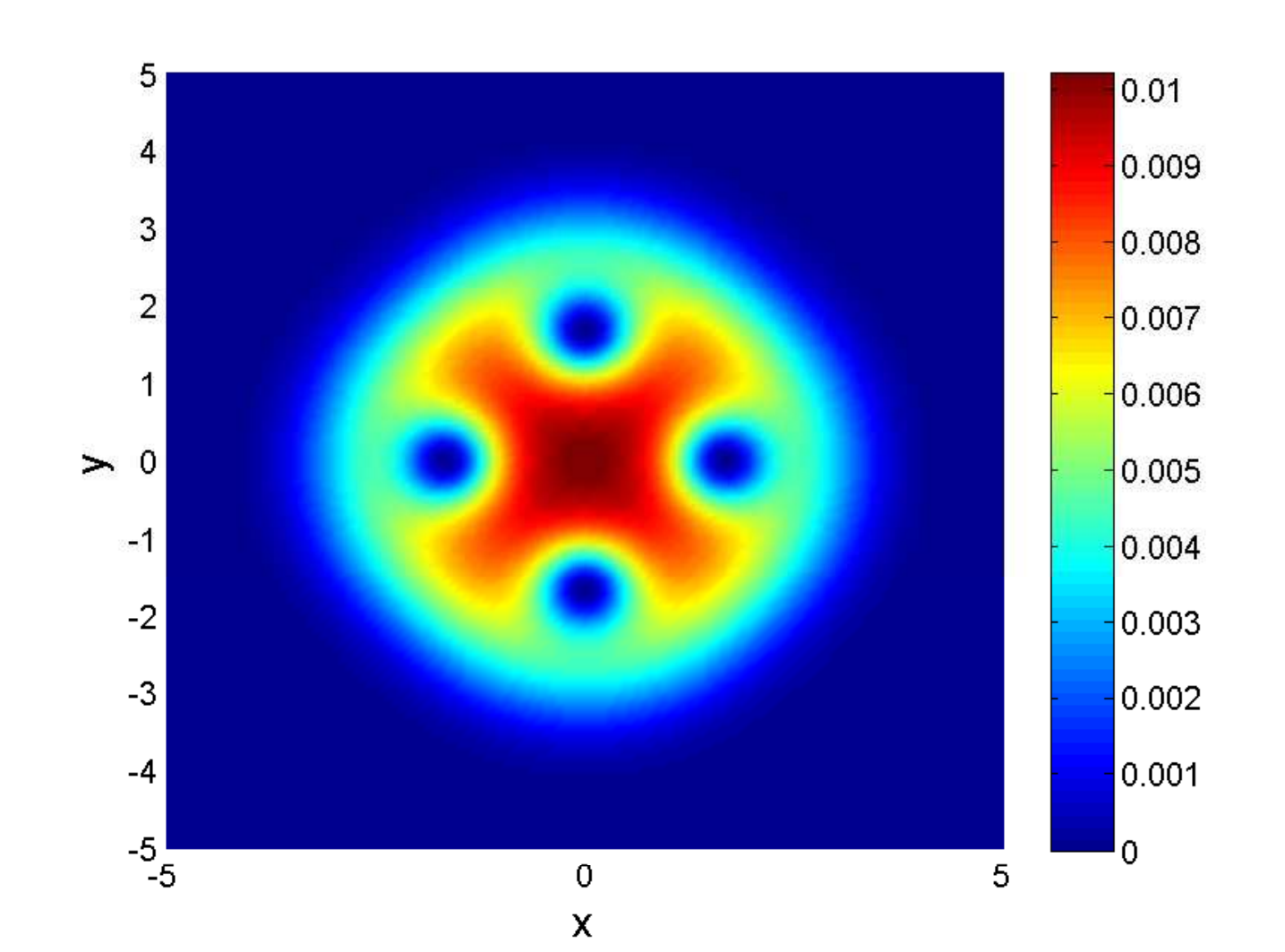}
 	\end{minipage}}
 	\subfigure[$t=4$]{
 		\begin{minipage}[b]{0.3\linewidth}
 			\includegraphics[width=1\linewidth]{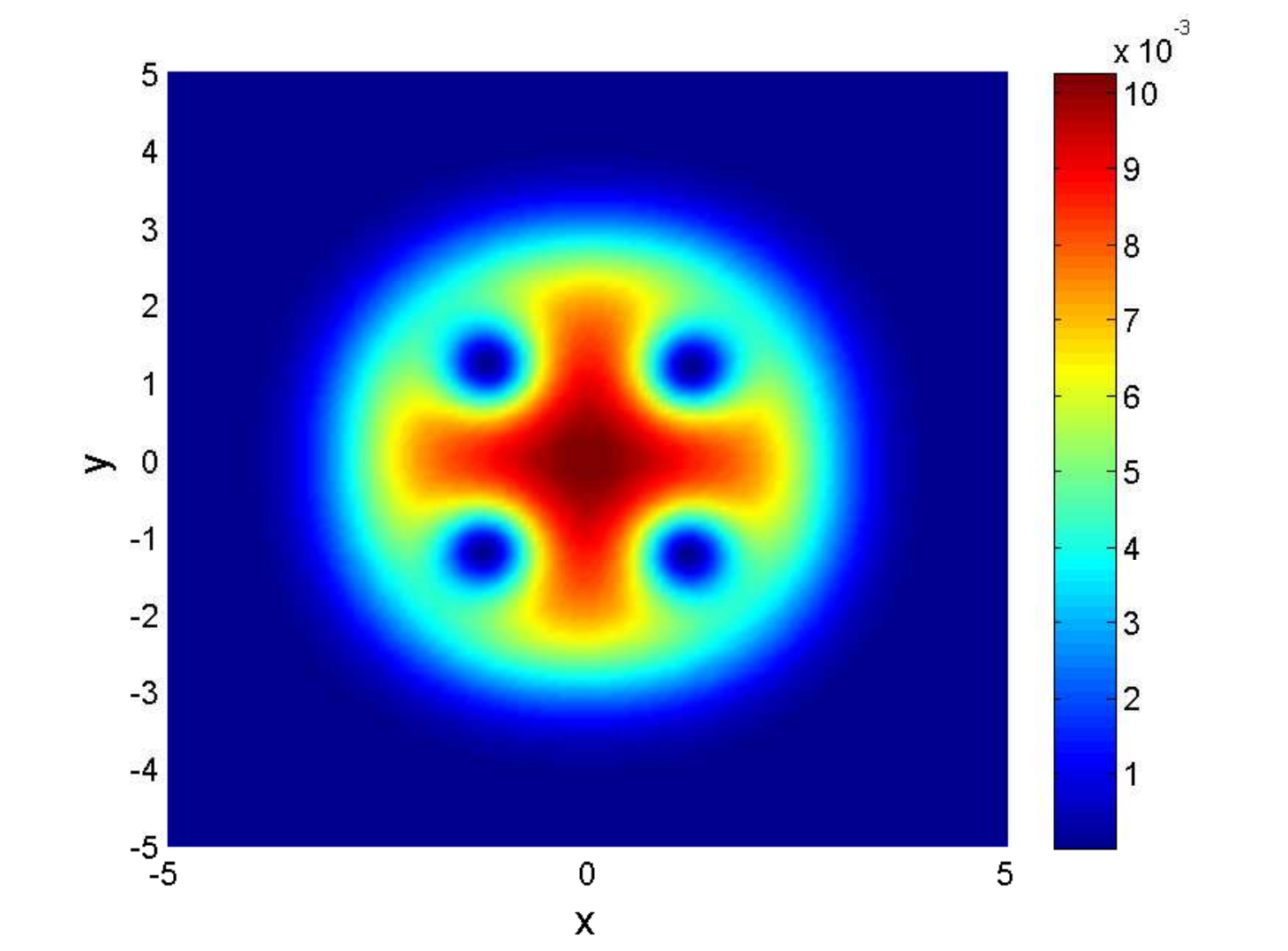}
 	\end{minipage}}
 	\subfigure[$t=8$]{
 		\begin{minipage}[b]{0.3\linewidth}
 			\includegraphics[width=1\linewidth]{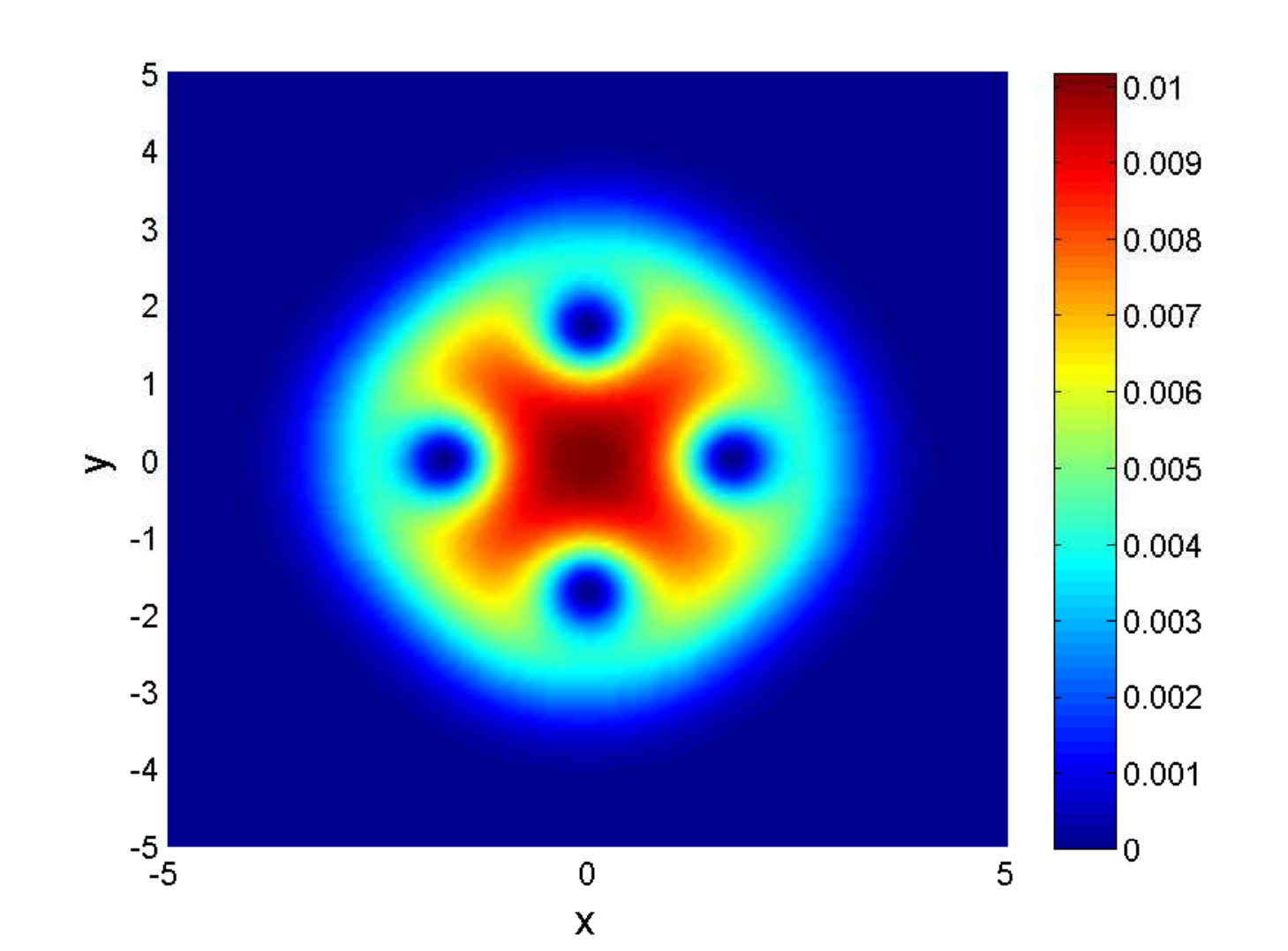}
 	\end{minipage}}
 	\subfigure[$t=12$]{
 		\begin{minipage}[b]{0.3\linewidth}
 			\includegraphics[width=1\linewidth]{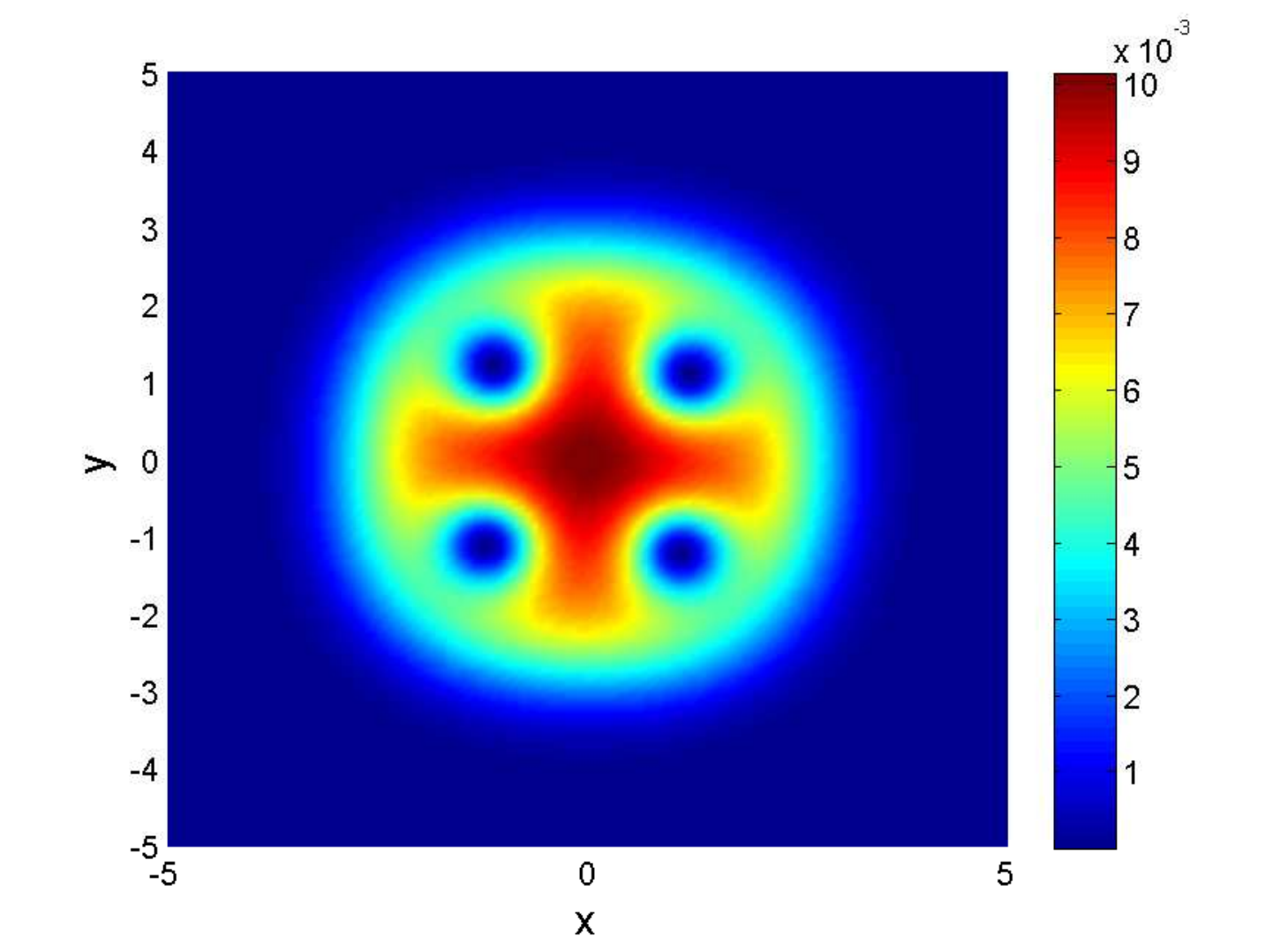}
 	\end{minipage}}
 	\subfigure[$t=16$]{
 		\begin{minipage}[b]{0.3\linewidth}
 			\includegraphics[width=1\linewidth]{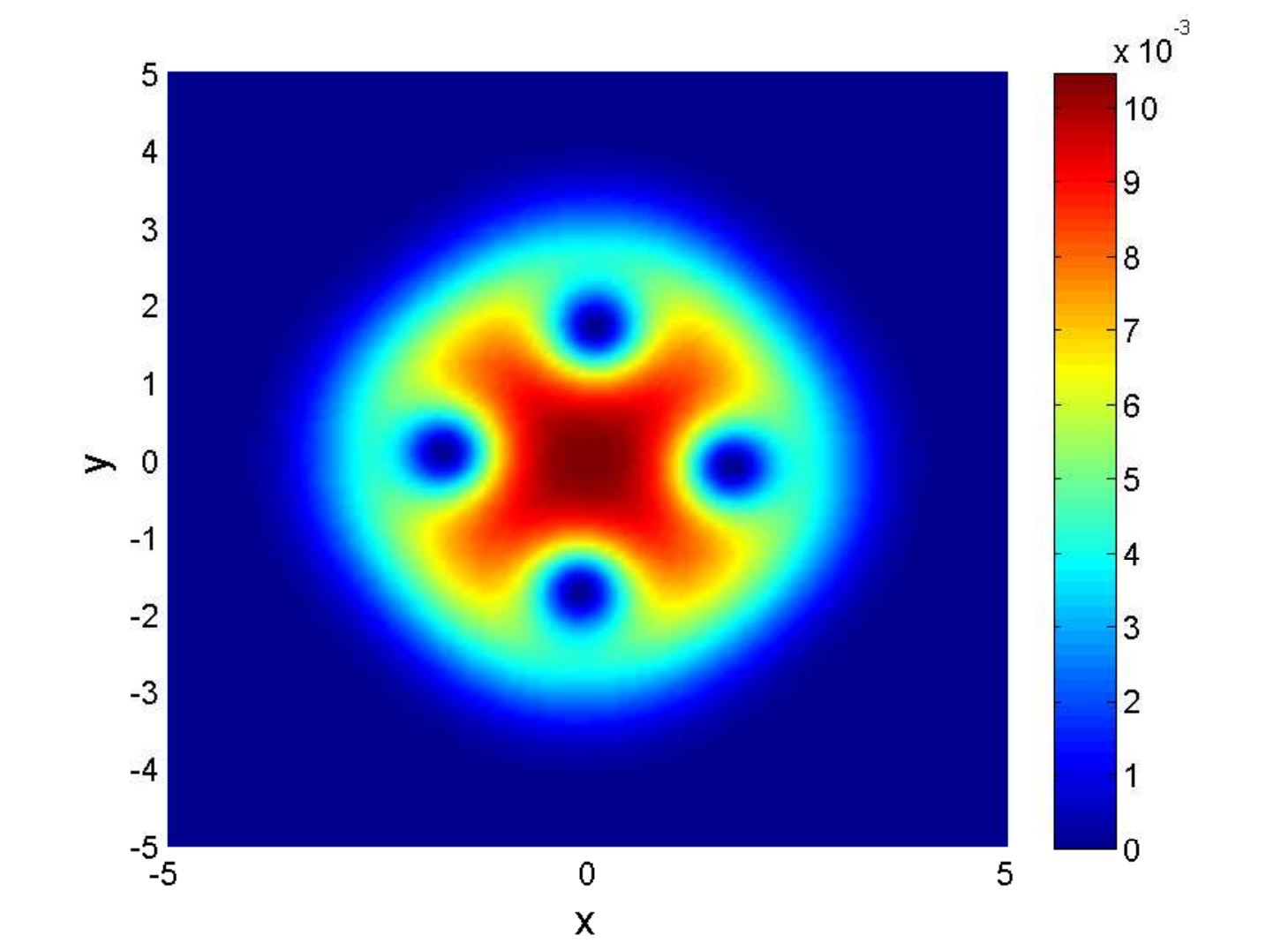}
 	\end{minipage}}
 	\subfigure[$t=20$]{
 		\begin{minipage}[b]{0.3\linewidth}
 			\includegraphics[width=1\linewidth]{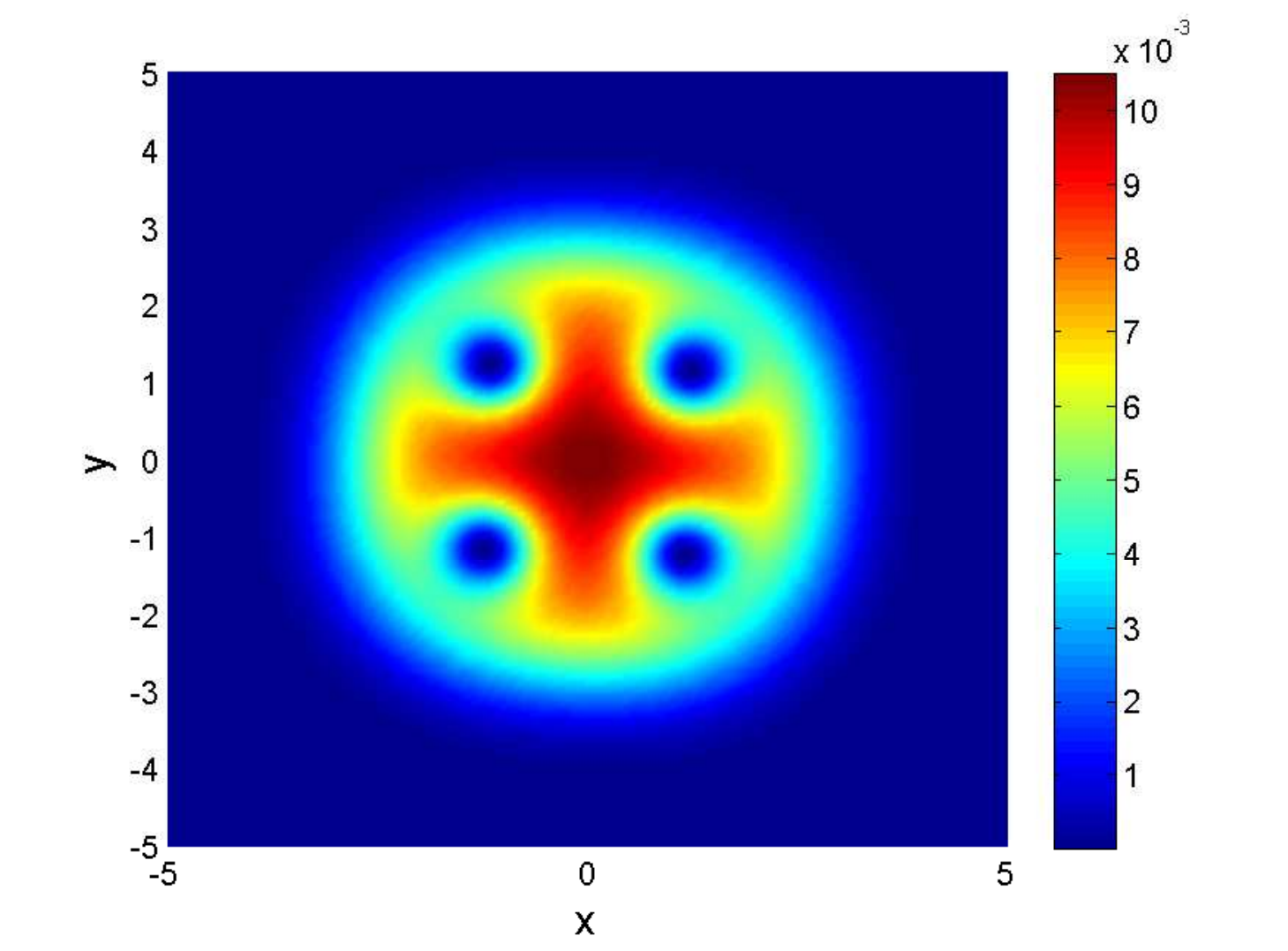}
 	\end{minipage}}
 	\caption{Slices along the surface $z=0$ of $|\psi|^2$ at different times computed by the \textbf{EIP-M} method for the dynamics of vortex lines, intercepted on the domain $[-5, 5]^2$. }\label{ex4-fig5}
 \end{figure}

\begin{figure}[H]
	\centering
	\subfigure[\textbf{EIP-M}]{
		\begin{minipage}[b]{0.3\linewidth}
			\includegraphics[width=1\linewidth]{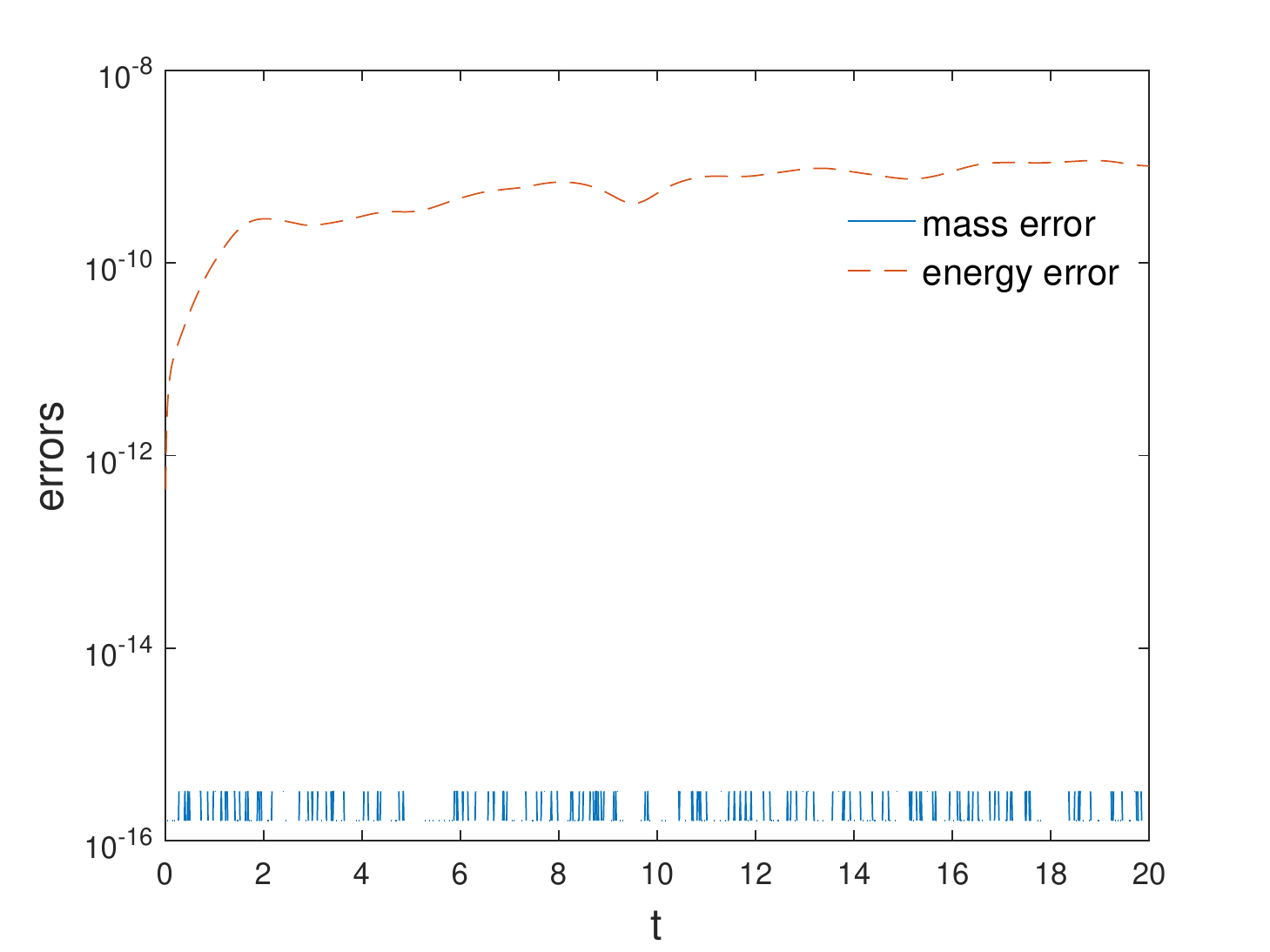}
	\end{minipage}}
	\subfigure[\textbf{EIP-E}]{
		\begin{minipage}[b]{0.3\linewidth}
			\includegraphics[width=1\linewidth]{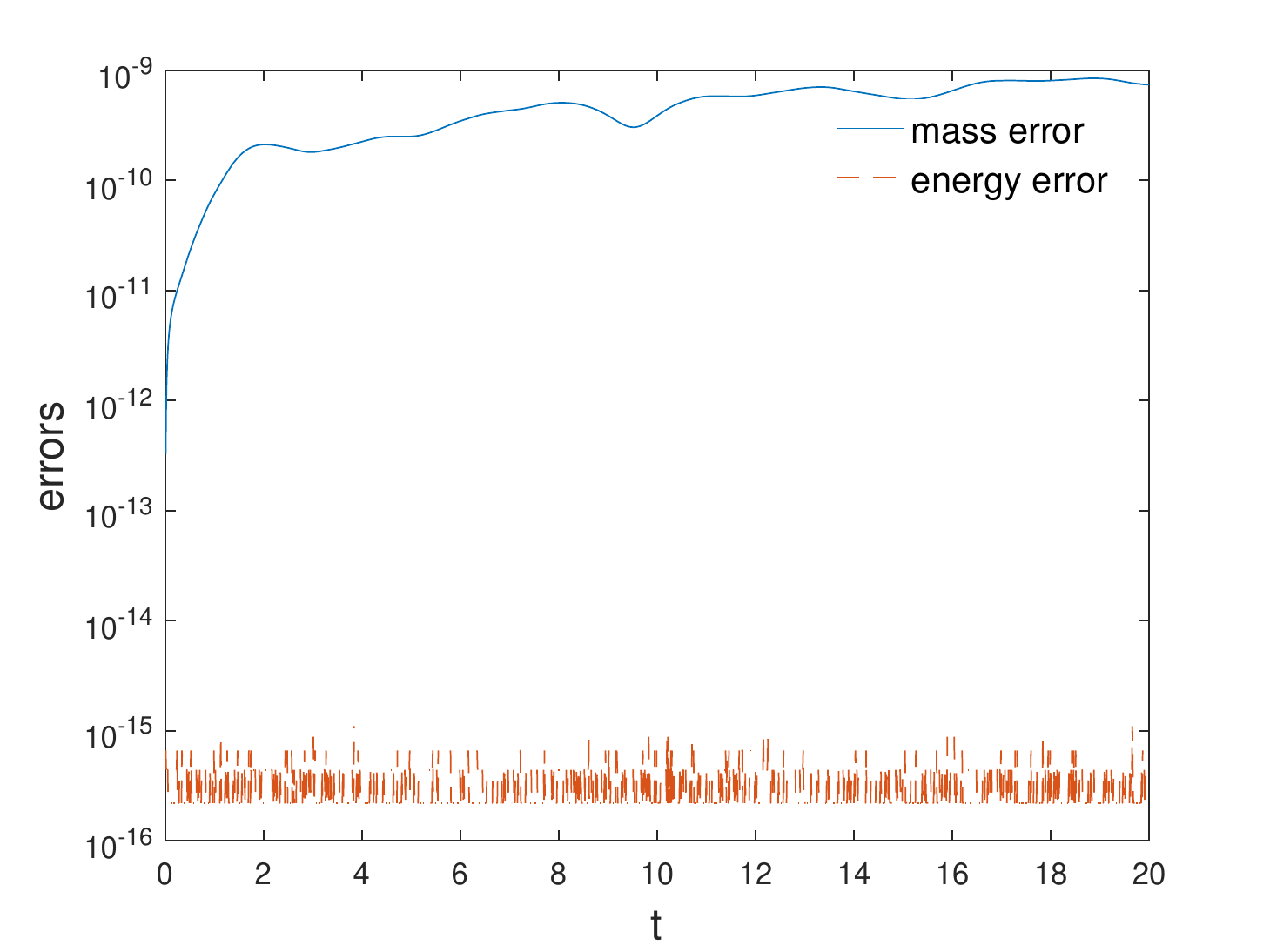}
	\end{minipage}}
	\subfigure[\textbf{EIP-ME}]{
		\begin{minipage}[b]{0.3\linewidth}
			\includegraphics[width=1\linewidth]{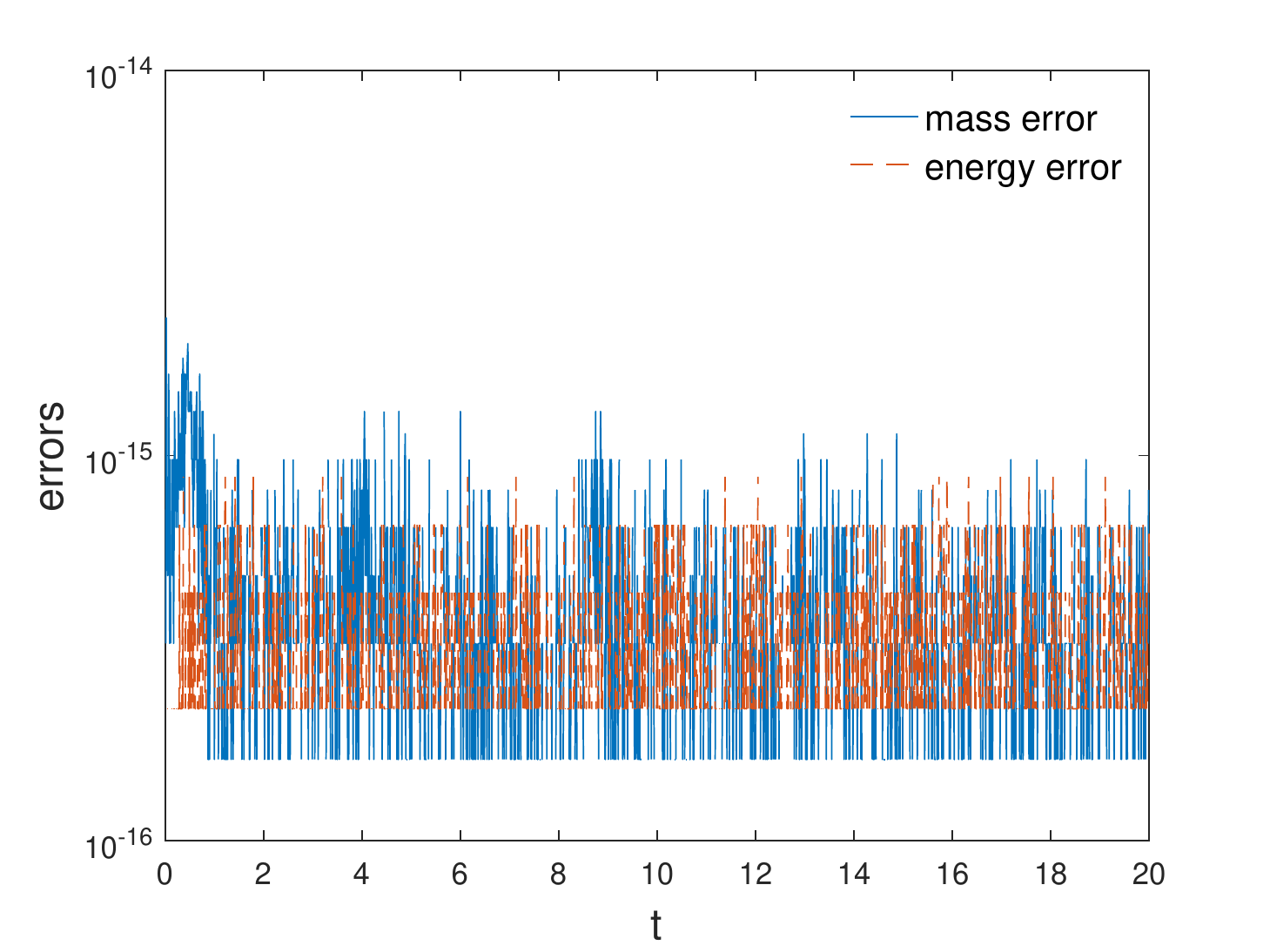}
	\end{minipage}}
	\caption{Errors in two invariants computed by the EIP methods for the dynamics of vortex lines.}\label{ex4-fig6}
\end{figure}

Finally, we consider the EIP methods for the dynamics of vortex lines in a 3D rotating BEC. The ground state solution is obtained similar as the above 2D experiment with computational domain  $\mathcal{D}=[-10,10]\times[-10,10]\times[-15,15]$, and the parameters are chosen as $\Omega=0.7$, $\beta=400$. The partition numbers $J=K=64$ and the time step $h=0.005$. The potential function $V_0$ in the Gross-Pitaevskii operator becomes a quadratic potential as
\[
V_0(\bm x)=\frac{1}{2}(\gamma_x^2 x^2+\gamma_y^2 y^2+\gamma_z^2 z^2),
\]
with $\gamma_x=\gamma_y=1, \gamma_z=1/2$. The $10^{-3}$-isosurface of the modulus of the corresponding ground state solution is drawn in the first plot of Figure.~\ref{ex4-fig4} where four vortex lines can be observed  clearly. By perturbing the angular velocity $\Omega=0.7$ to $0.9$, we can generate dynamics of the four vortex lines by the \textbf{EIP-M} method, and the other methods have similar results. Figure.~\ref{ex4-fig4} also presents the snapshots of $10^{-3}$-isosurfaces at different times where the four vortex lines rotate clockwise around the $z$- axis from above. This can be confirmed again by the slices along the surface $z=0$ in Figure.~\ref{ex4-fig5}. Moreover, as demonstrated in Figure.~\ref{ex4-fig6} the errors in the mass and energy exhibit a similar tendency as that in the 2D dynamics of vortex lattices.

\section{Conclusion}

In this paper, we propose a novel explicit and practically invariants-preserving method for conservative systems, which can be viewed as a further simplification of the standard projection method. We prove that such simplification does not affect the numerical behaviors, that is, the order of accuracy retains the same as the underlying Runge-Kutta method and the invariants are preserved to round-off errors in practice. The detailed implementations are also provided to show the flexibility in preservation of single or multiple invariants, and in the generalization to high-dimensional problems. Extensive numerical experiments are carried out for both ODEs and PDEs to verify the theoretical analysis and demonstrate the efficiency and conservative properties of our method.

Notice that the underlying method is taken as an explicit Runge-Kutta method throughout the entire discussion, however, the proposed invariants-preserving method is rather general and various existing non-conservative but efficient numerical methods can also be brought into our framework to achieve invariants preservation in practice and may essentially improve the numerical performance of the original methods.

\section*{Acknowledgments}
	
This work is supported by the National Key Research and Development Project of China (Grant Nos. 2017YFC0601505, 2018YFC1504205), the National Natural Science Foundation of China (Grant Nos. 11971242,11801269,11771213), the Major Projects of Natural Sciences of University in Jiangsu Province of China (Grant No. 18KJA110003), the Natural Science Foundation of Jiangsu Province (Grant Nos. BK20180413, BK20171480), the Foundation of Jiangsu Key Laboratory for Numerical Simulation of Large Scale Complex Systems (202002).

\appendix
\renewcommand\thesection{\appendixname~\Alph{section}}
\renewcommand\theequation{\Alph{section}.\arabic{equation}}

\section{Initial datum of the solar system}

\begin{table}[H]
	\centering
	\caption{Initial datum of planets in the solar system.}\label{tab-5-1}
	\begin{tabular*}{0.9\textwidth}[h]{@{\extracolsep{\fill}}c c c c} \hline
		Planet  & Initial position & Initial velocity & G*mass  \\ \hline
		& 1.563021412664830e+10 & -5.557001175482630e+04&    \\[1ex]
		Mercury  & 4.327888220902108e+10 & 1.840863017229157e+04  & 2.203209e+13  \\[1ex]
		& 2.102123103174893e+09 & 6.602621285552567e+03 &  \\[1ex]
		\hline
		& -9.030189258080004e+10 & -1.907374632532257e+04 & \\[1ex]
		\multirow{3}{*}{Venus}  	& 5.802615456116644e+10 & -2.963461693326599e+04 &  3.248586e+14 \\[1ex]
		& 6.006513603716755e+09 & 6.946391255404438e+02 &  \\ [1ex]
		\hline
		& -1.018974476358996e+11 & -2.201749257051057e+04 &      \\[1ex]
		Earth  & 1.065689158175689e+11 & -2.071074857788741e+04  & 3.986004e+14  \\[1ex]
		& -3.381951053601424e+06 & 1.575245213712245e+00& \\ [1ex]
		\hline
		& -2.443763125844157e+11 & -3.456935754608896e+03 &    \\[1ex]
		Mars  & 4.473211564076996e+10 & -2.176307370133160e+04& 4.282830e+13   \\[1ex]
		& 6.935657388967808e+09 & -3.711433859326417e-02  &   \\[1ex]
		\hline
		& -2.3516546827532200e+11 & -1.262559929908801e+04 &    \\[1ex]
		Jupiter &  7.421837640432589e+11 & -3.332552395475581e+03  & 1.266865e+17 \\[1ex]
		&  2.179850895804323e+09 & 2.962741332356101e+02  & \\[1ex]
		\hline
		& -1.011712827283427e+12 & 6.507898648442419e+03 &    \\[1ex]
		Saturn  & -1.077496255617324e+12 & -6.640809674126991e+03 & 3.793120e+16    \\[1ex]
		& 5.901251900068215e+10 & -1.434198106014633e+02 & \\[1ex]
		\hline
		& 2.934840841770302e+12  & -1.433852081777671e+03  &    \\[1ex]
		Uranus  & 6.048399137411513e+11  & 6.347897341634990e+03  & 5.793966e+15   \\[1ex]
		& -3.576451387567792e+10 & 4.228261484335974e+01  &  \\[1ex]
		\hline
		& 4.055112581124043e+12 & 2.275119229131818e+03 & -  \\[1ex]
		Neptune  	& -1.914578873112663e+12 & 4.942356914027413e+03 & 6.835107e+15   \\[1ex]
		& -5.400973716179796e+10 & -1.548950389954096e+02  &  \\[1ex]
		\hline
		& 9.514009594170194e+11 & 5.431808363374300e+03 &   \\[1ex]
		Pluto   & -4.776029500570151e+12 & -2.387056445508962e+01  & 8.72400e+11   \\[1ex]
		& 2.358627841705075e+11 & -1.551877289694926e+03 &  \\ \hline
	\end{tabular*}
\end{table}

\end{document}